\documentclass[12pt,reqno]{amsart}
\usepackage[margin=1in]{geometry}
\usepackage{amsmath,amssymb,amsthm,graphicx,amsxtra, setspace}
\usepackage[utf8]{inputenc}
\usepackage{mathrsfs}
\usepackage{hyperref}
\usepackage{xcolor}
\usepackage{upgreek}
\usepackage{mathtools}
\allowdisplaybreaks
\newtheorem{theorem}{Theorem}[section]
\newtheorem{lem}[theorem]{Lemma}
\newtheorem{rem}[theorem]{Remark}

\newtheorem{Def}[theorem]{Definition}

\newtheorem{Ex}[theorem]{Example}
\newtheorem{Ass}[theorem]{Assumption}
\usepackage{latexsym}
\usepackage{hyperref}
\usepackage{graphicx}
\usepackage{epstopdf}
\allowdisplaybreaks

\let\originalleft\left
\let\originalright\right
\renewcommand{\left}{\mathopen{}\mathclose\bgroup\originalleft}
\renewcommand{\right}{\aftergroup\egroup\originalright}

\newcommand{\Addresses}{{% additional braces for segregating \footnotesize
		\footnote{
			\noindent	 \textsuperscript{1,3}Department of Applied Sciences and Engineering, Indian Institute of Technology Roorkee-IIT Roorkee, Roorkee, 247667, India.

			\noindent  \textit{e-mail\textsuperscript{1}:} \texttt{sumit@as.iitr.ac.in.}
			
			\noindent  \textit{e-mail\textsuperscript{3}:} \texttt{jay.dabas@gmail.com.}

			\noindent \textsuperscript{2}Department of Mathematics, Indian Institute of Technology Roorkee-IIT Roorkee, Roorkee, 247667, India.\par\nopagebreak
			\noindent  \textit{e-mail\textsuperscript{2}:} \texttt{maniltmohan@ma.iitr.ac.in, maniltmohan@gmail.com.}

			\noindent \textsuperscript{*}Corresponding author.

			\textit{Key words:} approximate controllability, non-instantaneous impulses,fractional differential equation.
			
			Mathematics Subject Classification (2020): 34K06, 34A12, 37L05, 93B05.

}}}

\begin{document}
	\title[Approximate controllability of fractional systems]{Approximate controllability of fractional order non-instantaneous impulsive functional evolution equations with state-dependent delay in Banach spaces\Addresses}
	\author [S. Arora, M. T. Mohan and J. Dabas]{S. Arora\textsuperscript{1}, Manil T. Mohan\textsuperscript{2}  and J. dabas\textsuperscript{3*}}
	\maketitle{}
	\begin{abstract} The present paper deals with the control problems governed by fractional non-instantaneous impulsive functional evolution equations with state-dependent delay involving Caputo fractional derivatives in Banach spaces. The main objective of this work is to formulate sufficient conditions for the approximate controllability of the considered system in separable reflexive Banach spaces. We have exploited the resolvent operator technique and Schauder's fixed point theorem in the proofs to achieve this goal. The approximate controllability of linear system is discussed in detail, which lacks in the existing literature. We also provide an example to illustrate the efficiency of the developed results.  Moreover, we point out some shortcomings of the existing works in the context of characterization of mild solution and phase space, and approximate controllability of fractional order impulsive systems in Banach spaces. 
	\end{abstract}
	%%%%%%%%%%%%%%%%%%%%%%%%%%%%%%%%%%%%%%%%%%%%%%%%%%%%%%%%%%%%
	\section{Introduction}\label{intro}\setcounter{equation}{0}
	
	Many physical processes such as harvesting, natural disaster, shocks etc, cause abrupt changes in their states at certain time instant. These sudden changes occur for negligible time period and they are estimated in the form of instantaneous impulses. The theory of instantaneous impulsive systems has remarkable applications in several areas of science and engineering, for example, population dynamics, ecology, network control system with scheduling protocol etc., (cf. \cite{NE,AM,TY}, etc). However, certain dynamics of evolution processes in pharmacotherapy cannot be modeled by instantaneous impulsive dynamical systems, for example, in hemodynamical equilibrium of a person, introduction of insulin into the bloodstream and the consequent absorption of  the body are gradual processes and  stay active for a finite time interval. Thus, we cannot describe this situation via instantaneous impulsive systems. Therefore, Hern\'andez et. al. \cite{EHD} introduced a new class of impulses termed as non-instantaneous impulses, which starts at an arbitrary fixed point and stays active on a finite time interval and they established the existence of solutions for such a class of impulsive differential equations. Later, Wang et. al. \cite{JRW,JRY}, extended this model to two general classes of impulsive differential equations, which are very important in the study of dynamics of evolutionary processes in pharmacotherapy. For more details on the theory of non-instantaneous impulsive systems, we refer the interested readers to \cite{YTJ,JMF}, etc, and the references therein. In addition to this, there are various real world phenomena, for example, neural network, inferred grinding models, ecological models, heat conduction in materials with fading memory etc. In these phenomena, the current state of a system is influenced by the past states. The dynamics of such processes are characterized by delay differential equations (finite or infinite), see for example, \cite{XFL,LA,JW}, etc. In the application point of view, the functional evolution systems with state-dependent delays are more prevalent and adequate, see for instance, \cite{AO,CNF}, etc, and the references therein. 
	
	The concept of controllability plays a vital role in the study of control systems. Controllability (exact or approximate) refers that the solution of a control system can steer from an arbitrary initial state to a desired final state by using some control function. In the infinite dimensional setting, in comparison with exactly controllable systems, approximate controllable systems are more extensive and have wide range of applications, cf. \cite{M,TRR,TR,EZ}, etc. In the past two decades, the problem of approximate controllability of various kinds of systems (in Hilbert and Banach spaces) such as impulsive differential equations, functional differential equations, stochastic systems, Sobolev type evolution systems, etc,  is extensively studied  with the help of fixed point approach and produced excellent results, see for instance, \cite{SM,SAM,FUX,AGK,M,ER}, etc. 
	
	On the other hand, fractional differential equations (FDEs), which involve fractional derivatives of the form $\frac{d^{\alpha}}{dt^{\alpha}}$, where $\alpha>0$ is not necessarily an integer, attained great importance due to their ability to model complex phenomena. They naturally appear  in diffusion processes, electrical engineering, viscoelasticity, control theory of dynamical systems, quantum mechanics, biological sciences, electromagnetic theory, signal processing, finance, economics, and many other fields (cf. \cite{HlR,AAK,VM,IPF,TPJ,YAR,MST}, etc).  A comprehensive study on fractional calculus and FDEs are available in \cite{AAK,IPF}, etc.  For the past few decades, FDEs in infinite dimensions seek incredible attention of many researchers and eminent contributions have been made both in theory as well as in applications. Several authors studied the existence and approximate controllability results for fractional order systems in Hilbert spaces, see for instance, \cite{SNS,NIZ,SRY,JWY}, etc, and the references therein. 
	
	The study of approximate controllability of the fractional order control systems in Banach spaces has not got much attention in the literature. In \cite{NMI}, Mahmudov developed sufficient conditions for the approximate controllability of the Sobolev type fractional evolution equations with Caputo derivative in separable reflexive Banach spaces using the Schauder fixed point theorem. Later in \cite{MIN}, he studied the approximate controllability of the  fractional neutral evolution systems by taking infinite delay using Krasnoselkii’s fixed-point theorem. After that Chalishazar et. al. \cite{DNC} extended his work by considering instantaneous impulses and examined the approximate controllability in Banach spaces. 
	
	The articles  \cite{PCY,FZM,RS,RGR}, etc claimed the approximate controllability of the factional order systems in general Banach spaces using resolvent operator condition. But, the resolvent operator defined in these works is valid only if the state space is a Hilbert space, whose dual is identified by the space itself (see, the resolvent operator definition in the expression \eqref{2.1} and Remark \ref{rem2.8}). Moreover, many papers deal with the fractional order impulsive systems with delays, see for instance, \cite{ZT,ZH,YZO}, etc. In these works, the characterization of norm or seminorm defined in the phase space involves uniform norm, but the choice of such a norm or seminorm is not suitable in the impulsive case, for counter examples and more details, we refer the interested readers to \cite{GU} (a detailed discussion on this problem is also available in \cite{SAMT}). Many articles considered the fractional order impulsive systems with non-instantaneous impulses (cf. \cite{RMS,Ad,Zyf,Zy} etc.). In theses works, the concept of the mild solution defined for the considered system is not realistic, a counter example and appropriate definition of the mild solution discussed in \cite{Fe,JRY} (see Definition \ref{def2.6} for the mild solution definition for system under our consideration). One of the main aims of this work is to resolve these issues.
	
	Recently, a few papers have been reported on the approximate controllability of the non-instantaneous impulsive systems with and without delays in Hilbert spaces, (cf. \cite{RMS,SKS,SAJ}, etc). Dhayal et. al. \cite{RMS} formulated the approximate controllability results for a class of fractional order non-instantaneous impulsive stochastic differential equations driven by fractional Brownian motion. In \cite{SKS}, Kumar and Abdal derived sufficient conditions for the approximate controllability of non-instantaneous impulsive fractional semilinear measure driven control systems with infinite delay. Liu and his co-authors, in \cite{SAJ},  investigated the approximate controllability of fractional differential equation with non-instantaneous impulses via iterative learning control scheme. To best of our knowledge, there is no work has been reported on the approximate controllability of the fractional order non-instantaneous impulsive systems with state-dependent delay in Banach spaces.
	
	Motivated from the above facts, in this work, we derive sufficient conditions for the approximate controllability of  fractional order non-instantaneous impulsive functional evolution equations with state-dependent delay in separable reflexive Banach spaces. Moreover, we properly define the resolvent operator in Banach spaces, which plays a crucial role in obtaining the aforementioned results (see the expression \eqref{2.1} below). Proper motivation for the construction of different forms of feedback controls  for the fractional order semilinear systems available in the literature has been justified in this work (see Remarks \ref{rem3.6} and \ref{rem4.4} below). Furthermore, our paper modifies the phase space characterization to incorporate Guedda’s observations in \cite{GU}, by replacing the uniform norm on the phase space by integral norm for the impulsive differential equations (see Example \ref{exm2.8}). 
	
	We consider the following fractional order non-instantaneous impulsive functional evolution equation with state-dependent delay: 
	\begin{equation}\label{1.1}
	\left\{
	\begin{aligned}
	^C\mathrm{D}_{0,t}^{\alpha}x(t)&=\mathrm{A}x(t)+\mathrm{B}u(t)+f(t,x_{\rho(t, x_t)}), \ t\in\bigcup_{k=0}^{m} (\tau_k, t_{k+1}]\subset J=[0,T], \\
	x(t)&=h_k(t, x(t_k^{-})),\ t\in(t_k, \tau_k], \ k=1,\dots,m, \\
	x_{0}&=\psi\in \mathfrak{B},
	\end{aligned}
	\right.
	\end{equation}
	where \begin{itemize} \item $^C\mathrm{D}_{0,t}^{\alpha}$ denotes a derivative in Caputo sense of order $\alpha$ with $\frac{1}{2}<\alpha<1$, \item  the operator $\mathrm{A}:\mathrm{D(A)}\subset \mathbb{X}\to\mathbb{X}$ is an infinitesimal generator of a $\mathrm{C}_0$-semigroup $ \mathcal{T}(t) $ on a separable reflexive Banach space $ \mathbb{X}$ (having a strictly convex dual $\mathbb{X}^*$), \item  the linear operator $\mathrm{B}:\mathbb{U}\to\mathbb{X}$ is bounded with $\left\|\mathrm{B}\right\|_{\mathcal{L}(\mathbb{U},\mathbb{X})}= \tilde{M}$ and the control function $u\in \mathrm{L}^{2}(J;\mathbb{U})$, where $\mathbb{U}$ is a separable Hilbert space, \item  the function $ f:J\times \mathfrak{B}\rightarrow \mathbb{X} $, where $\mathfrak{B}$ is a phase space, which will be specified in the subsequent sections, \item for $k=1,\ldots,m$, the functions $h_k:[t_k, \tau_k]\times\mathbb{X}\to\mathbb{X}$ represent the non-instantaneous impulses and the fixed points $\tau_k$ and $t_k$ satisfy $0=t_0=\tau_0<t_1\le \tau_1\le t_2<\ldots<t_m\le \tau_m\le t_{m+1}=T$,   the values  $x(t_k^+)$ and $x(t_k^-)$ stand for the right and left limit of $x(t)$ at the point $t=t_k,$ respectively and satisfy $x(t_k^-)=x(t_k),$ for $k=1,\ldots,m$, \item the function $x_{t}:(-\infty, 0]\rightarrow\mathbb{X}$ with $x_{t}(\theta)=x(t+\theta)$ belongs to some abstract space $\mathfrak{B}$, and the function $\rho: J \times \mathfrak{B} \rightarrow (-\infty, T]$ is continuous.
	\end{itemize}
	
	The rest of the article is organized as follows: Section \ref{pre} starts with some basic definitions of fractional calculus and then presents the phase space axioms. Moreover, we provide assumptions and important results to establish the existence and approximate controllability results of the system \eqref{1.1}. Section \ref{Linear} deals with the  approximate controllability  of  the fractional order linear control system corresponding to \eqref{1.1}. For the approximate controllability results of the linear system, we first formulate the linear regulator problem to obtain the existence of an optimal control (Theorem \ref{optimal}), and then derive the explicit expression of optimal control in Lemma \ref{lem3.1}. With the help of this optimal control, we investigate the approximate controllability of the linear system \eqref{3.2} in Theorem \ref{lem4.2}. In section \ref{semilinear}, we first show the existence of a mild solution of the system \eqref{1.1}, by invoking Schauder's fixed point theorem. Then, we demonstrate the approximate controllability results of the system \eqref{1.1} in Theorem \ref{thm4.4}. A concrete example is presented in the final section to illustrate the developed results in previous sections.
	%%%%%%%%%%%%%%%%%%%%%%%%%%%%%%%%%%%%%%%%%%%%%%%%%%%%%%%%%%%
	\section{Preliminaries}\label{pre}\setcounter{equation}{0}
	Let us first recall some basic definitions and properties from fractional calculus. For a detailed study, the interested readers are referred to see \cite{AAK,IPF}, etc. 	Let $\mathrm{AC}([a,b];\mathbb{R})$ denote the space of all absolutely continuous functions from $[a,b]$ to $\mathbb{R}$ and $\mathrm{AC}^n([a,b];\mathbb{R}),$ for $n\in\mathbb{N},$ represent the space of all  functions $f:[a,b]\to\mathbb{R}$, which have continuous derivatives up to order $n-1$  with $f^{(n-1)}\in\mathrm{AC}([a,b];\mathbb{R})$. 
	\begin{Def}
		The \emph{fractional integral} of a function $f:[a,b]\to\mathbb{R}$, $a,b\in\mathbb{R}$ with $a<b$, of order $q>0$ is defined as
		\begin{align*}
		I_{a}^{q}f(t):=\frac{1}{\Gamma(q)}\int_{a}^{t}\frac{f(s)}{(t-s)^{1-q}}\mathrm{d}s,\ \mbox{ for a.e. } \  t\in[a,b],
		\end{align*}
		where $f\in\mathrm{L}^1([a,b];\mathbb{R})$ and $\Gamma(\alpha)=\int_{0}^{\infty}t^{\alpha-1}e^{-t}\mathrm{d}t$ is the Euler gamma function.
	\end{Def}
	\begin{Def}%[\cite{MR}]
		The \emph{Riemann-Liouville fractional derivative} of a function $f:[a,b]\to\mathbb{R}$  of order $q>0$ is given as 
		\begin{align*}
		^L\mathrm{D}_{a,t}^{q}f(t):=\frac{1}{\Gamma(n-q)}\frac{d^n}{dt^n}\int_{a}^{t}(t-s)^{n-q-1}f(s)\mathrm{d}s,\ \mbox{ for a.e. }\ t\in[a,b],
		\end{align*}
		with $n-1< q<n,$ and the function $f\in\mathrm{AC}^n([a,b];\mathbb{R})$.
	\end{Def}
	\begin{Def}%[\cite{MR}]
		The \emph{Caputo fractional derivative} of a function $f\in\mathrm{AC}^n([a,b];\mathbb{R})$ of order $q>0$   is defined as 
		\begin{align*}
		^C\mathrm{D}_{a,t}^{q}f(t):=\ ^L\mathrm{D}_{a,t}^{q}\left[f(t)-\sum_{p=1}^{n-1}\frac{f^{(p)}(a)}{p!}(x-a)^p\right],\ \mbox{ for a.e. } \ t\in[a,b].
		\end{align*}
	\end{Def}
	%\begin{Def}
	%	A real valued function $f\in\mathrm{C}_\gamma,\ \gamma\in\mathbb{R}$ if there exist a real number $p>\gamma$, such that $f(t)=t^pg(t)$, where $g\in\mathrm{C}[0, \infty)$ and $f\in\mathrm{C}^n_\gamma$ if and only if $f^{(n)}\in\mathrm{C}_\gamma,\ n\in\mathbb{N}$.
	%\end{Def}
	\begin{rem}[\cite{IPF}]
		If a function $f\in\mathrm{AC}^n([a,b];\mathbb{R})$, then the Caputo fractional derivative of $f$ can be written as
		\begin{align*}
		^C\mathrm{D}_{a,t}^{q}f(t)=\frac{1}{\Gamma(n-q)}\int_{a}^{t}(t-s)^{n-q-1}f^{(n)}(s)\mathrm{d}s,\ \mbox{ for a.e. }\ \ t \in[a,b], \ n-1< q<n.
		\end{align*}
	\end{rem}

	In order to define the mild solution of the system \eqref{1.1}, first we consider the one-dimensional stable probability density function (cf. \cite{Ua} \cite{MPG})
	\begin{align*}
	\phi_{q}(\xi)=\frac{1}{\pi}\sum_{n=1}^{\infty}(-1)^{n-1}\xi^{-nq -1}\frac{\Gamma(nq+1)}{n!}\sin(n\pi q),\ \xi\in(0,\infty), 0<q<1,
	\end{align*}
	and for any $x\in\mathbb{X}$, we define
	\begin{align*}
	\mathcal{T}_{q}(t)x=\int_{0}^{\infty}\!\!\varphi_{q}(\xi)\mathcal{T}(t^{q}\xi)x\mathrm{d}\xi \ \mbox{ and }\
	\widehat{\mathcal{T}}_{q}(t)x=q\int_{0}^{\infty}\!\!\xi\varphi_{q}(\xi)\mathcal{T}(t^{q}\xi)x\mathrm{d}\xi,
	\end{align*}
	where $\varphi_{q}(\xi)=\frac{1}{q}\xi^{-1-\frac{1}{q}}\phi_{q}(\xi^{-\frac{1}{q}})$ and  $ \mathcal{T}(t) $ is the $\mathrm{C}_0$-semigroup. Note that $\varphi_{q}$ is also a probability density function on $(0, \infty)$, so that
	\begin{align*}
	\int_{0}^{\infty}\varphi_{q}(\xi)\mathrm{d}\xi=1\ \mbox{ and }\ \int_{0}^{\infty}\xi\varphi_q(\xi)\mathrm{d}\xi=\frac{1}{\Gamma(1+q)},
	\end{align*}
	where the final expression is the first moment of $\varphi_q$.  
	
	Let us provide some important properties of the operators $\mathcal{T}_{q}(t)$ and $\widehat{\mathcal{T}}_{q}(t),$ for $t\geq 0$.
	\begin{lem}[\cite{YZ}] \label{lem2.5}The operators $\mathcal{T}_{q}(t)$ and $\widehat{\mathcal{T}}_{q}(t)$ have the following properties:
		\begin{enumerate}
			\item [(i)] For any fixed $t\ge0,$ the operators $\mathcal{T}_{q}(t)$ and  $\widehat{\mathcal{T}}_{q}(t)$ are linear and bounded. Moreover
			\begin{align*}
			\left\|\mathcal{T}_{q}(t)\right\|_{\mathcal{L}(\mathbb{X})}\le M \ \mbox{ and }\  \|\widehat{\mathcal{T}}_{q}(t)\|_{\mathcal{L}(\mathbb{X})}\le\frac{M q}{\Gamma(1+q)},
			\end{align*} 
			where $M$ is a constant such that $\left\|\mathcal{T}(t)\right\|_{\mathcal{L(\mathbb{X})}}\leq M$.
			\item [(ii)] The operators $\mathcal{T}_{q}(t)$ and  $\widehat{\mathcal{T}}_{q}(t)$ are strongly continuous for $t\ge0$.
			\item [(iii)] If $\mathcal{T}(t)$ is compact for $t>0$, then the operators $\mathcal{T}_{q}(t)$ and  $\widehat{\mathcal{T}}_{q}(t)$ are also compact for $t>0$. 
		\end{enumerate}
	\end{lem}
	Let us define the set 
	\begin{align*} 
	\mathrm{PC}(J;\mathbb{X})&:=\big\{x:J \rightarrow \mathbb{X} : x\vert_{t\in I_k}\in\mathrm{C}(I_k;\mathbb{X}),\ I_k:=(t_k, t_{k+1}],\ k=0,1,\ldots,m \ \mbox{ and }\ x(t_k^+)\\&\qquad \qquad \mbox{ and }\ x(t_k^-)\ \mbox{ exist for each }\ k=1,\ldots,m, \ \mbox{ and satisfy }\ x(t_k)=x(t_k^-)\big\}, 
	\end{align*}
	endowed with the norm $\left\|x\right\|_{\mathrm{PC}(J;\mathbb{X})}:=\sup\limits_{t\in J}\left\|x\right\|_{\mathbb{X}}$.
	
	We now introduce the concept of mild solution for the system \eqref{1.1} (cf. \cite{JRY}). 
	\begin{Def}[Mild solution]\label{def2.6}
		A function $x(\cdot;\psi,u):(-\infty, T]\to\mathbb{X}$   is said to be a  \emph{mild solution} of \eqref{1.1}, if  $x_0=\psi\in\mathfrak{B}$ and $x\vert_{J}\in\mathrm{PC}(J;\mathbb{X})$ and satisfies the following:
		\begin{equation}\label{2.2}
		x(t)=
		\begin{dcases}
		\mathcal{T}_{\alpha}(t)\psi(0)+\int_{0}^{t}(t-s)^{\alpha-1}\widehat{\mathcal{T}}_{\alpha}(t-s)\left[\mathrm{B}u(s)+f(s,x_{\rho(s, x_s)})\right]\mathrm{d}s,\ t\in[0, t_1],\\
		h_k(t, x(t_k^-)),\ t\in(t_k, \tau_k],\ k=1,\ldots,m,\\
		\mathcal{T}_{\alpha}(t-\tau_k)h_k(\tau_k, x(t_k^-))-\int_{0}^{\tau_k}(\tau_k-s)^{\alpha-1}\widehat{\mathcal{T}}_{\alpha}(\tau_k-s)\left[\mathrm{B}u(s)+f(s,x_{\rho(s, x_s)})\right]\mathrm{d}s\\\quad+\int_{0}^{t}(t-s)^{\alpha-1}\widehat{\mathcal{T}}_{\alpha}(t-s)\left[\mathrm{B}u(s)+f(s,x_{\rho(s, x_s)})\right]\mathrm{d}s,\ t\in(\tau_k,t_{k+1}],\ k=1,\ldots,m.
		\end{dcases}
		\end{equation}
	\end{Def}
	%%%%%%%%%%%%%%%%%%%%%%%%%%%%%%%%%%%%%%%%%%%%%%%%%%%%%%%%%%%%%%%%
	\subsection{Phase space}
	We now provide the definition of phase space $\mathfrak{B}$ introduced in \cite{HY}, and suitably modify to incorporate the impulsive systems (cf. \cite{VOj}). The linear space $\mathfrak{B}$ equipped with the seminorm $\left\|\cdot\right\|_{\mathfrak{B}}$, consisting of all functions from $(-\infty, 0]$ into $\mathbb{X}$  and satisfying the following axioms:
	\begin{enumerate}
		\item [(A1)] If $x: (-\infty, T]\rightarrow \mathbb{X}$ such that $x_{0}\in \mathfrak{B}$ and $x|_{J}\in \mathrm{PC}(J;\mathbb{X})$. Then the following conditions hold:
		\begin{itemize}
			\item [(i)] $x_{t}\in\mathfrak{B}$ for $t\in J$.
			\item [(ii)] $\left\|x_{t}\right\|_{\mathfrak{B}}\leq \Lambda(t)\sup\{\left\|x(s)\right\|_{\mathbb{X}}: 0
			\leq s\leq t\}+\Upsilon(t)\left\|x_{0}\right\|_{\mathfrak{B}},$  for $t\in J$, where $\Lambda, \Upsilon:[0, \infty)\rightarrow [0, \infty)$ are independent of $x$, the function $\Lambda(\cdot)$ is strictly positive and continuous, $\Upsilon(\cdot)$ is locally bounded.
			%\item [(iii)] the map $t\mapsto y_t$  is a $\mathfrak{B}$-valued continuous function on $J.$
		\end{itemize}
		\item [(A2)] The space $\mathfrak{B}$ is complete. 
	\end{enumerate}  
	For any $\psi\in \mathfrak{B}$, the function $\psi_{t}, \ t\leq 0,$ defined as $\psi_{t}(\theta)=\psi(t+\theta),\  \theta \in (-\infty, 0].$ Then for any function $x(\cdot)$ satisfying the axiom (A1) with $x_{0}=\psi$, we can extend the mapping $t\mapsto x_{t}$ by setting $x_{t}=\psi_{t}, \  t\leq0$, to the whole interval $(-\infty, T]$. Moreover, let us  introduce a set 
	\begin{align*}
	\mathcal{Q}(\rho^-)&=\{\rho(s, \varphi):\  \rho(s, \varphi)\leq 0, \mbox{ for } (s, \varphi)\in J \times \mathfrak{B}\}.
	\end{align*} 
	Assume that the function $t\mapsto \psi_{t}$,  defined from $\mathcal{Q}(\rho^-)$ into $\mathfrak{B}$ is continuous, and there exists a continuous and bounded function $\varTheta^{\psi}: \mathcal{Q}(\rho^-)\rightarrow (0, \infty)$ such that  
	\begin{align*}
	\left\|\psi_{t}\right\|_{\mathfrak{B}}&\leq \varTheta^{\psi}(t)\left\|\psi\right\|_{\mathfrak{B}}.
	\end{align*}
	\begin{lem}[\cite{HY}]{\label{lem2.7}}
		Let $x:(-\infty, T]\rightarrow \mathbb{X}$ be a function such that $x_0=\psi$ and $x\vert_J\in\mathrm{PC}(J;\mathbb{X})$. Then 
		\begin{align*}
		\left\|x_s\right\|_{\mathfrak{B}}&\leq H_{1}\left\|\psi\right\|_{\mathfrak{B}}+H_{2}\sup \big\{ \left\|x(\theta)\right \|_{\mathbb{X}}:\theta \in [0, \max \{0, s\}] \big\}, \ s\in \mathcal{Q}(\rho^-) \cup J,
		\end{align*}
		where $$H_{1}= \sup\limits_{t\in \mathcal{Q}(\rho^-)} \Theta^{\psi}(t)  + \sup\limits_{t\in J}\Upsilon(t), \;\; H_{2}= \sup\limits_{t\in J}\Lambda(t).$$ 
	\end{lem}
	\begin{Ex}\label{exm2.8}
		Let us take $\mathfrak{B}=\mathrm{PC}_{r}\times\mathrm{L}^p_h(\mathbb{X}), r\ge0, 1\le p<\infty$, which consists of all functions $\psi:(-\infty,0]\to\mathbb{X}$ such that $\psi\vert_{[-r,0]}\in \mathrm{PC}([-r,0];\mathbb{X}),$ Lebesgue measurable on $(-\infty,-r)$ and $h\|\psi(\cdot)\|^p_{\mathbb{X}}$ is Lebesgue integrable on $(-\infty,-r]$. The seminorm in $\mathfrak{B}$ is defined as 
		\begin{align}
		\label{Bnorm}\left\|\psi\right\|_{\mathcal{P}}:=\int_{-r}^{0}\left\|\psi(\theta)\right\|_\mathbb{X}\mathrm{d}\theta+\left(\int_{-\infty}^{-r}h(\theta)\left\|\psi(\theta)\right\|^p_{\mathbb{X}}\mathrm{d}\theta\right)^{\frac{1}{p}},
		\end{align}
		where the function $h:(-\infty, 0]\to\mathbb{R}^+$ is locally bounded and Lebesgue integrable. Moreover, there exists a locally bounded function $H:(-\infty,0]\to\mathbb{R}^+$ such that $h(t+\theta)\le H(t)h(\theta),$ for all $t\le0$ and $\theta\in(-\infty,0)\backslash \mathcal{O}_t$ where $\mathcal{O}_t\subseteq(-\infty,0)$ is a set with Lebesgue measure zero. 
	\end{Ex}
	%%%%%%%%%%%%%%%%%%%%%%%%%%%%%%%%%%%%%%%%%%%%%%%%%%%%%%%%%%%%%%%%%%%%%%%
	\subsection{Resolvent operator and assumptions}
	To discuss the approximate controllability of the system \eqref{1.1}, we first define the following operators:
	\begin{equation}\label{2.1}
	\left\{
	\begin{aligned}
	L_0^Tu&:=\int^{T}_{0}(T-t)^{\alpha-1}\widehat{\mathcal{T}}_{\alpha}(T-t)\mathrm{B}u(t)\mathrm{d}t,\\
	\Phi_{0}^{T}&:=\int^{T}_{0}(T-t)^{2(\alpha-1)}\widehat{\mathcal{T}}_{\alpha}(T-t)\mathrm{B}\mathrm{B}^*\widehat{\mathcal{T}}_{\alpha}(T-t)^*\mathrm{d}t,\\
	\mathcal{R}(\lambda,\Phi_{0}^{T})&:=(\lambda \mathrm{I}+\Phi_{0}^{T}\mathcal{J})^{-1},\ \lambda > 0,
	\end{aligned}
	\right.
	\end{equation}
	where $\mathrm{B}^{*}$ and $\widehat{\mathcal{T}}_{\alpha}(t)^*$ denote the adjoint operators of $\mathrm{B}$ and $\widehat{\mathcal{T}}_{\alpha}(t),$ respectively. It is immediate that the operator $L_0^T$ is linear and bounded for $\frac{1}{2}<\alpha<1$. Moreover, the map $\mathcal{J} : \mathbb{X} \rightarrow 2^{\mathbb{X}^*}$ stands for the duality mapping, which is defined as 
	\begin{align*}
	\mathcal{J}[x]&=\{x^* \in \mathbb{X}^* : \langle x, x^* \rangle=\left\|x\right\|_{\mathbb{X}}^2= \left\|x^*\right\|_{\mathbb{X}^*}^2 \}, \mbox{ for all } x\in \mathbb{X},
	\end{align*}
	where  $\langle \cdot, \cdot  \rangle $  represents  a duality pairing between $\mathbb{X}$ and $\mathbb{X}^*$. Since the space $\mathbb{X}$ is a reflexive Banach space, then  $\mathbb{X}^*$ becomes strictly convex (see, \cite{AA}), which implies that the mapping $\mathcal{J}$ is bijective, strictly monotonic and demicontinuous, that is, $$x_k\to x\ \mbox{ in }\ \mathbb{X}\ \mbox{ implies }\ \mathcal{J}[x_k] \xrightharpoonup{w} \mathcal{J}[x]\ \mbox{ in } \ \mathbb{X}^*\ \mbox{ as }\ k\to\infty.$$ Moreover, the inverse mapping $\mathcal{J}^{-1}:\mathbb{X}^*\to\mathbb{X}$ is also duality mapping.
	\begin{rem}\label{rem2.8}
		If $\mathbb{X}$ is a separable Hilbert space (identified with its own dual), then the resolvent operator is defined as $\mathcal{R}(\lambda,\Phi_{0}^{T}):=(\lambda \mathrm{I}+\Phi_{0}^{T})^{-1},\ \lambda > 0$.
	\end{rem}
	\begin{lem}[Lemma 2.2 \cite{M}]\label{lem2.9}
		For every $h\in\mathbb{X}$ and $\lambda>0$, the equation
		\begin{align}\label{2.4}\lambda z_{\lambda}+\Phi_{0}^{T}\mathcal{J}[z_{\lambda}]=\lambda h,\end{align}
		has a unique solution  $z_{\lambda}(h)=\lambda(\lambda \mathrm{I}+\Phi_{0}^{T}\mathcal{J})^{-1}(h)=\lambda\mathcal{R}(\lambda,\Phi_{0}^{T})(h)$ and \begin{align}\label{2.5}
		\left\|z_{\lambda}(h)\right\|_{\mathbb{X}}=\left\|\mathcal{J}[z_{\lambda}(h)]\right\|_{\mathbb{X}^*}\leq\left\|h\right\|_{\mathbb{X}}.
		\end{align}
		\begin{proof}
			Since the non-negative operator $\Phi_0^T$ is linear and bounded for $\frac{1}{2}<\alpha<1$, then proceeding similar way as in the proof of Lemma 2.2 \cite{M}, one can obtain the results.
		\end{proof}
	\end{lem}

	\begin{Def}[\cite{FUX}]
		The system \eqref{1.1} is said to be \emph{approximately controllable} on $ J $, for any initial function $\psi\in\mathfrak{B}$, if the closure of reachable set is whole space $\mathbb{X}$, that is, $\overline{\mathfrak{R}(T,\psi)}=\mathbb{X},$ where the reachable set is defined as \begin{align*}\mathfrak{R}(T,\psi) = \{x(T;\psi,u): u(\cdot) \in \mathrm{L}^{2}(J;\mathbb{U})\}.\end{align*}
	\end{Def}
	We impose the following assumptions to investigate the approximate controllability of the system \eqref{1.1}:
	\begin{Ass}\label{as2.1} 
		\begin{enumerate}
			\item [\textit{($H0$)}] For every $h\in\mathbb{X}$, $z_\lambda=z_{\lambda}(h)=\lambda\mathcal{R}(\lambda,\Phi_{0}^{T})(h) \rightarrow 0$ as $\lambda\downarrow 0$ in strong topology, where $z_{\lambda}(h)$ is a solution of the equation \eqref{2.4}.
			\item[\textit{(H1)}] The $\mathrm{C}_0$-semigroup of bounded linear operator $\mathcal{T}(t)$ is compact for $t>0$ with bound $M\geq 1,$ such that $\|\mathcal{T}(t)\|_{\mathcal{L(\mathbb{X})}}\leq M$.
			\item [\textit{(H2)}] 
			\begin{enumerate} 
				\item [(i)] Let $x:(-\infty, T]\rightarrow \mathbb{X}$ be such that $x_0=\psi$ and $x|_{J}\in \mathrm{PC}(J;\mathbb{X}).$ The function $t\mapsto f(t, x_{\rho(t,x_{t})}) $ is strongly measurable on $J$  and the function $f(t,\cdot): \mathfrak{B}\rightarrow \mathbb{X}$ is continuous for a.e. $t\in J$. Also, the map $t\mapsto f(s,x_{t})$ is continuous on $\mathcal{Q}(\rho^-) \cup J,$ for every $s\in J$. 
				\item [(ii)] For each positive integer $r$, there exists a constant $\alpha_1\in[0,\alpha]$ and a function $\gamma_{r}\in \mathrm{L}^{\frac{1}{\alpha_1}}(J;\mathbb{R^{+}})$, such that $$  \sup_{\left\| \psi\right\|_{\mathfrak{B}}\leq r} \left\|f(t, \psi)\right\|_{\mathbb{X}}\leq\gamma_{r}(t), \mbox{ for a.e.} \ t \in J \mbox{ and } \psi\in \mathfrak{B},$$ with
				$$ \liminf_{r \rightarrow \infty } \frac {\left\|\gamma_r\right\|_{\mathrm{L}^{\frac{1}{\alpha_1}}(J;\mathbb{R^+})}}{r} = \beta< \infty. $$ 
			\end{enumerate}
			\item [\textit{(H3)}] The impulses $ h_k:[t_k,\tau_k]\times\mathbb{X}\to\mathbb{X}$, for  $k=1,\dots,m$, are such that 
			\begin{itemize}
				\item [(i)]The impulses $h_k(\cdot,x):[t_k,\tau_k]\to\mathbb{X}$ are continuous for each $x\in\mathbb{X}$. 
				\item [(ii)] Each $h_k(t,\cdot):\mathbb{X}\to\mathbb{X}$ is completely continuous, for all  $t\in[t_k,\tau_k]$.
				\item[(iii)] $\left\|h_k(t,x) \right\|_{\mathbb{X}} \leq l_k, \mbox{ for each } t\in [t_k, \tau_k] \mbox{ and } x \in \mathbb{X},$ where $l_k$'s are positive constants.
			\end{itemize} 
		\end{enumerate}
	\end{Ass}
	The following version of the discrete Gronwall-Bellman lemma (cf. \cite{Ch}) is used in the sequel.  
	\begin{lem}\label{lem2.13}
	If $\{f_n\}_{n=0}^{\infty}, \{g_n\}_{n=0}^{\infty}$ and $\{w_n\}_{n=0}^{\infty}$ are non-negative sequences and $$f_n\le g_n+\sum_{k=0}^{n-1}w_kf_k,\ \text{ for }\ n\geq 0,$$  then $$f_n\le g_n+\sum_{k=0}^{n-1}g_kw_k\exp\left(\sum_{j=k+1}^{n-1}g_j\right),\text{  for }\ n\geq 0.$$
	\end{lem}
	
	%%%%%%%%%%%%%%%%%%%%%%%%%%%%%%%%%%%%%%%%%%%%%%%%%%%%%%%%%%%%%%%%%
	\section{Linear Control Problem} \label{Linear}\setcounter{equation}{0}
	The present section is devoted for discussing the approximate controllability of the fractional order linear control problem corresponding to \eqref{1.1}. To establish this result, we first obtain the existence of an optimal control by minimizing the cost functional given by
	\begin{equation}\label{3.1}
	\mathcal{G}(x,u)=\left\|x(T)-x_{T}\right\|^{2}_{\mathbb{X}}+\lambda\int^{T}_{0}\left\|u(t)\right\|^{2}_{\mathbb{U}}\mathrm{d}t,
	\end{equation}
	where $x(\cdot)$ is the solution of the linear control system:
	\begin{equation}\label{3.2}
	\left\{
	\begin{aligned}
	^CD_{0,t}^{\alpha}x(t)&= \mathrm{A}x(t)+\mathrm{B}u(t),\ t\in J,\\
	x(0)&=\zeta,
	\end{aligned}
	\right.
	\end{equation}
	with the control $u\in \mathbb{U}$, $x_{T}\in \mathbb{X}$ and $\lambda >0$. Since $\mathrm{B}u\in\mathrm{L}^1(J;\mathbb{X})$, the system \eqref{3.2} has a unique mild solution $x\in \mathrm{C}(J;\mathbb{X}) $ given by (see Corollary 2.2, Chapter 4, \cite{P} and Lemma 4.68, Chapter 4, \cite{Y})
	\begin{align*}
	x(t)= \mathcal{T}_{\alpha}(t)\zeta+\int_{0}^{t}(t-s)^{\alpha-1}\widehat{\mathcal{T}}_{\alpha}(t-s)\mathrm{B}u(s)\mathrm{d}s,
	\end{align*}
	for any $u\in\mathscr{U}_{\mathrm{ad}}=\mathrm{L}^2(J;\mathbb{U})$ (admissible control class).  Next, we define the \emph{admissible class} $\mathscr{A}_{\mathrm{ad}}$  for the system \eqref{3.2} as
	\begin{align*}
	\mathscr{A}_{\mathrm{ad}}:=\big\{(x,u) :x\mbox{ is \mbox{the unique mild solution} of }\eqref{3.2}  \mbox{ with the control }u\in\mathscr{U}_{\mathrm{ad}}\big\}.
	\end{align*}
	For any given control $u\in\mathscr{U}_{\mathrm{ad}}$, the system \eqref{3.2} has a unique mild solution, which ensures that the set $\mathscr{A}_{\mathrm{ad}}$ is nonempty. By using the definition of the cost functional, we can formulate the optimal control problem as:
	\begin{align}\label{3.3}
	\min_{ (x,u) \in \mathscr{A}_{\mathrm{ad}}}  \mathcal{G}(x,u).
	\end{align}
	In the next theorem, we show the existence of an optimal pair for the problem \eqref{3.3}.
	\begin{theorem}[Existence of an optimal pair]\label{optimal}
		For a given $\zeta\in\mathbb{X}$ and fixed $\frac{1}{2}<\alpha<1$, there exists a unique optimal pair  $(x^0,u^0)\in\mathscr{A}_{\mathrm{ad}}$ for the problem \eqref{3.3}.
	\end{theorem}
	\begin{proof}
		Let us assume $$\mathcal{G} := \inf \limits _{u \in \mathscr{U}_{\mathrm{ad}}}\mathcal{G}(x,u).$$ Since $0\leq \mathcal{G} < +\infty$, there exists a minimizing sequence $\{u^n\}_{n=1}^{\infty} \in \mathscr{U}_{\mathrm{ad}}$ such that $$\lim_{n\to\infty}\mathcal{G}(x^n,u^n) = \mathcal{G},$$ where $x^n(\cdot)$ is the unique mild solution of the system \eqref{3.2}, corresponding to the control $u^n(\cdot),$  for each $n\in\mathbb{N}$ with $x^n(0)=\zeta$.	Note that $x^n(\cdot)$ satisfies
		\begin{align}\label{3.4}
		x^n(t)&=\mathcal{T}_{\alpha}(t)\zeta+\int_{0}^{t}(t-s)^{\alpha-1}\widehat{\mathcal{T}}_{\alpha}(t-s)\mathrm{B}u^n(s)\mathrm{d}s,
		\end{align} 
		for $t\in J$. Since  $0\in\mathscr{U}_{\mathrm{ad}}$, without loss of generality, we may assume that $\mathcal{G}(x^n,u^n) \leq \mathcal{G}(x,0)$, where $(x,0)\in\mathscr{A}_{\mathrm{ad}}$. Using the definition of $\mathcal{G}(\cdot,\cdot)$, we easily get
		\begin{align}\label{35}
		\left\|x^n(T)-x_{T}\right\|^{2}_{\mathbb{X}}+\lambda\int^{T}_{0}\left\|u^n(t)\right\|^{2}_{\mathbb{U}}\mathrm{d}t\leq \left\|x(T)-x_{T}\right\|^{2}_{\mathbb{X}}\leq 2\left(\|x(T)\|_{\mathbb{X}}^2+\|x_T\|_{\mathbb{X}}^2\right)<+\infty.
		\end{align}
		From the above estimate, it is clear that, there exists a large $L>0$ (independent of $n$), such that 
		\begin{align}\label{3.5}\int_0^T \|u^n(t)\|^2_{\mathbb{U}} \mathrm{d} t \leq L < +\infty .\end{align}
		Using the expression \eqref{3.4}, we compute
		\begin{align*}
		\left\|x^n(t)\right\|_{\mathbb{X}}&\le\left\|\mathcal{T}_{\alpha}(t)\zeta\right\|_{\mathbb{X}}+\left\|\int_{0}^{t}(t-s)^{\alpha-1}\widehat{\mathcal{T}}_{\alpha}(t-s)\mathrm{B}u^n(s)\mathrm{d}s\right\|_{\mathbb{X}}\nonumber\\ &\le\left\|\mathcal{T}_{\alpha}(t)\right\|_{\mathcal{L}(\mathbb{X})}\left\|\zeta\right\|_{\mathbb{X}}+\int_{0}^{t}(t-s)^{\alpha-1}\|\widehat{\mathcal{T}}_{\alpha}(t-s)\|_{\mathcal{L}(\mathbb{X})}\left\|\mathrm{B}\right\|_{\mathcal{L}(\mathbb{U},\mathbb{X})}\left\|u^n(s)\right\|_{\mathbb{U}}\mathrm{d}s\nonumber\\
		&\le M\left\|\zeta\right\|_{\mathbb{X}}+\frac{M\tilde{M}\alpha}{\Gamma(1+\alpha)}\int_{0}^{t}(t-s)^{\alpha-1}\left\|u^n(s)\right\|_{\mathbb{U}}\mathrm{d}s
		\nonumber\\&\le M\left\|\zeta\right\|_{\mathbb{X}}+\frac{M\tilde{M}\alpha}{\Gamma(1+\alpha)}\frac{T^{2\alpha-1}}{2\alpha-1}\left(\int_{0}^{t}\left\|u^n(s)\right\|_{\mathbb{U}}^2\mathrm{d}s\right)^{\frac{1}{2}}\nonumber\\&\le M\left\|\zeta\right\|_{\mathbb{X}}+\frac{M\tilde{M}\alpha}{\Gamma(1+\alpha)}\frac{T^{2\alpha-1}L^{\frac{1}{2}}}{2\alpha-1}<+\infty,
		\end{align*}
		for all  $t\in J$ and $\frac{1}{2}<\alpha<1$. Since we know that $\mathrm{L}^{2}(J;\mathbb{X})$ is reflexive, then by applying the Banach-Alaoglu theorem, we can find a subsequence $\{x^{n_k}\}_{k=1}^{\infty}$ of $\{x^n\}_{n=1}^{\infty}$ such that 
		\begin{align}\label{3.6}
		x^{n_k}\xrightharpoonup{w}x^0\ \mbox{ in }\ \mathrm{L}^{2}(J;\mathbb{X}), \ \mbox{ as }\ k\to\infty. 
		\end{align}
		From the estimate \eqref{3.5}, we also infer that the sequence $\{u^n\}_{n=1}^{\infty}$ is uniformly bounded in the space $\mathrm{L}^2(J;\mathbb{U})$. Further, by using the Banach-Alaoglu theorem, there exists a subsequence, say,  $\{u^{n_k}\}_{k=1}^{\infty}$ of $\{u^n\}_{n=1}^{\infty}$ such that 
		\begin{align*}
		u^{n_k}\xrightharpoonup{w}u^0\ \mbox{ in }\ \mathrm{L}^2(J;\mathbb{U})=\mathscr{U}_{\mathrm{ad}}, \ \mbox{ as }\ k\to\infty. 
		\end{align*}
		Since $\mathrm{B}$ is a bounded linear operator from $\mathbb{U}$ to $\mathbb{X}$, then we have 
		\begin{align}\label{3.7}
		\mathrm{B}	u^{n_k}\xrightharpoonup{w}\mathrm{B}u^0\ \mbox{ in }\ \mathrm{L}^2(J;\mathbb{X}),\ \mbox{ as }\ k\to\infty.
		\end{align}
		Moreover, by using the above  convergences together with the compactness of the operator $(\mathrm{Q}f)(\cdot) =\int_{0}^{\cdot}(\cdot-s)^{\alpha-1}\widehat{\mathcal{T}}_{\alpha}(\cdot-s) f(s)\mathrm{d}s:\mathrm{L}^2(J;\mathbb{X})\rightarrow \mathrm{C}(J;\mathbb{X}) $ (see Lemma \ref{lem2.12} below), we obtain
		\begin{align*}
		\left\|\int_{0}^{t}(t-s)^{\alpha-1}\widehat{\mathcal{T}}_{\alpha}(t-s)\mathrm{B}u^{n_k}(s)\mathrm{d}s-\int_{0}^{t}(t-s)^{\alpha-1}\widehat{\mathcal{T}}_{\alpha}(t-s)\mathrm{B}u(s)\mathrm{d}s\right\|_{\mathbb{X}}\to0,\ \mbox{ as }\  k\to\infty,
		\end{align*}
		for all $t\in J$. We now estimate 
		\begin{align}\label{3.8}
		\left\|x^{n_k}(t)-x^*(t)\right\|_{\mathbb{X}}&=\left\|\int_{0}^{t}(t-s)^{\alpha-1}\widehat{\mathcal{T}}_{\alpha}(t-s)\mathrm{B}u^{n_k}(s)\mathrm{d}s-\int_{0}^{t}(t-s)^{\alpha-1}\widehat{\mathcal{T}}_{\alpha}(t-s)\mathrm{B}u^{0}(s)\mathrm{d}s\right\|_{\mathbb{X}}\nonumber\\&\to 0,\ \mbox{ as } \ k\to\infty, \ \mbox{for all }\ t\in J,
		\end{align}
		where 
		\begin{align*}
		x^*(t)=\mathcal{T}_{\alpha}(t)\zeta+\int_{0}^{t}(t-s)^{\alpha-1}\widehat{\mathcal{T}}_{\alpha}(t-s)\mathrm{B}u^{0}(s)\mathrm{d}s,\ t\in J.
		\end{align*}
		It is clear by the above expression, the function  $x^*\in \mathrm{C}(J;\mathbb{X})$ is the unique mild solution of the equation \eqref{3.2} with the control $u^{0}\in\mathscr{U}_{\mathrm{ad}}$. Since the weak limit is unique, then by combining the convergences \eqref{3.6} and \eqref{3.8}, we obtain $x^*(t)=x^0(t),$ for all $t\in J$. Hence, the function $x^0$ is the unique mild solution of the system \eqref{3.2} with the control $u^{0}\in\mathscr{U}_{\mathrm{ad}}$ and also the whole sequence  $x^n\to x^0\in\mathrm{C}(J;\mathbb{X})$. Consequently, we have  $(x^0,u^0)\in\mathscr{A}_{\mathrm{ad}}$.
		
		It remains to show that the functional $\mathcal{G}(\cdot,\cdot)$ attains its minimum at  $(x^0,u^0)$, that is, \emph{$\mathcal{G}=\mathcal{G}(x^0,u^0)$}. Since the cost functional $\mathcal{G}(\cdot,\cdot)$ given in \eqref{3.1} is continuous and convex (see Proposition III.1.6 and III.1.10,  \cite{EI}) on $\mathrm{L}^2(J;\mathbb{X}) \times \mathrm{L}^2(J;\mathbb{U})$, it follows that $\mathcal{G}(\cdot,\cdot)$ is weakly lower semi-continuous (Proposition II.4.5, \cite{EI}). That is, for a sequence 
		$$(x^n,u^n)\xrightharpoonup{w}(x^0,u^0)\ \mbox{ in }\ \mathrm{L}^2(J;\mathbb{X}) \times  \mathrm{L}^2(J;\mathbb{U}),\ \mbox{ as }\ n\to\infty,$$
		we have 
		\begin{align*}
		\mathcal{G}(x^0,u^0) \leq  \liminf \limits _{n\rightarrow \infty} \mathcal{G}(x^n,u^n).
		\end{align*}
		Hence, we obtain 
		\begin{align*}\mathcal{G} \leq \mathcal{G}(x^0,u^0) \leq  \liminf \limits _{n\rightarrow \infty} \mathcal{G}(x^n,u^n)= \lim \limits _{n\rightarrow \infty} \mathcal{G}(x^n,u^n) = \mathcal{G},\end{align*}
		and thus $(x^0,u^0)$ is a minimizer of the problem \eqref{3.3}. Note that the cost functional given in \eqref{3.1} is convex, the constraint \eqref{3.2} is linear and the class $\mathscr{U}_{\mathrm{ad}}=\mathrm{L}^2(J;\mathbb{U})$ is convex, then the  optimal control obtained above is unique.
	\end{proof}

	In the following lemma, we prove the compactness of the operator $(\mathrm{Q}f)(\cdot) =\int_{0}^{\cdot}(\cdot-s)^{\alpha-1}\widehat{\mathcal{T}}_{\alpha}(\cdot-s) f(s)\mathrm{d}s:\mathrm{L}^2(J;\mathbb{X})\rightarrow \mathrm{C}(J;\mathbb{X}) ,\ \mbox{ for }\ \frac{1}{2}<\alpha<1,$ where we assume that $\mathbb{X}$ is a general Banach space. The case of $\alpha=1$ is available in Lemma 3.2, \cite{JYONG}. 
	\begin{lem}\label{lem2.12}
		Suppose that Assumptions (H1) holds. Let the operator $\mathrm{Q}:\mathrm{L}^{2}(J;\mathbb{X})\rightarrow \mathrm{C}(J;\mathbb{X})$ be defined as
		\begin{align}
		(\mathrm{Q}\psi)(t)= \int^{t}_{0}(t-s)^{\alpha-1}\widehat{\mathcal{T}}_{\alpha}(t-s)\psi(s)\mathrm{d}s, \ t\in J,\ \frac{1}{2}<\alpha<1.
		\end{align}
		Then the operator $\mathrm{Q}$ is compact.
	\end{lem}
	\begin{proof}
		We prove that $\mathrm{Q}$ is a compact operator by using  the infinite-dimensional version of Arzel\'a-Ascoli theorem (see, Theorem 3.7, Chapter 2, \cite{JYONG}). Let a closed and bounded ball $\mathcal{B}_R$ in $\mathrm{L}^{2}(J;\mathbb{X})$ be defined as 
		\begin{align*}
		\mathcal{B}_R=\left\{\psi\in \mathrm{L}^{2}(J;\mathbb{X}):\left\|\psi\right\|_{\mathrm{L}^{2}(J;\mathbb{X})}\leq R\right\}.
		\end{align*}
		For $s_1,s_2\in J$  ($s_1<s_2$) and $\psi\in\mathcal{B}_{R}$, we compute the following:
		\begin{align}
		&\left\|(\mathrm{Q}\psi)(s_2)-(\mathrm{Q}\psi)(s_1)\right\|_{\mathbb{X}}\nonumber\\&\le\left\|\int_{0}^{s_1}\left[(s_2-s)^{\alpha-1}\widehat{\mathcal{T}}_{\alpha}(s_2-s)-(s_1-s)^{\alpha-1}\widehat{\mathcal{T}}_{\alpha}(s_1-s)\right]\psi(s)\mathrm{d}s\right\|_{\mathbb{X}}\nonumber\\&\quad+\left\|\int_{s_1}^{s_2}(s_2-s)^{\alpha-1}\widehat{\mathcal{T}}_{\alpha}(s_2-s)\psi(s)\mathrm{d}s\right\|_{\mathbb{X}}\nonumber\\&\le\left\|\int_{0}^{s_1}\left[(s_2-s)^{\alpha-1}-(s_1-s)^{\alpha-1}\right]\widehat{\mathcal{T}}_{\alpha}(s_2-s)\psi(s)\mathrm{d}s\right\|_{\mathbb{X}}\nonumber\\&\quad+\left\|\int_{0}^{s_1}\!\!\!(s_1-s)^{\alpha-1}\left[\widehat{\mathcal{T}}_{\alpha}(s_2-s)-\widehat{\mathcal{T}}_{\alpha}(s_1-s)\right]\psi(s)\mathrm{d}s\right\|_{\mathbb{X}}+\left\|\int_{s_1}^{s_2}\!\!\!(s_2-s)^{\alpha-1}\widehat{\mathcal{T}}_{\alpha}(s_2-s)\psi(s)\mathrm{d}s\right\|_{\mathbb{X}}\nonumber\\&\le\frac{M\alpha}{\Gamma(1+\alpha)}\int_{0}^{s_1}\left|(s_2-s)^{\alpha-1}-(s_1-s)^{\alpha-1}\right|\left\|\psi(s)\right\|_{\mathbb{X}}\mathrm{d}s\nonumber\\&\quad+\int_{0}^{s_1}(s_1-s)^{\alpha-1}\left\|\widehat{\mathcal{T}}_{\alpha}(s_2-s)-\widehat{\mathcal{T}}_{\alpha}(s_1-s)\right\|_{\mathcal{L}(\mathbb{X})}\left\|\psi(s)\right\|_{\mathbb{X}}\mathrm{d}s\nonumber\\&\quad+\frac{M\alpha}{\Gamma(1+\alpha)}\int_{s_1}^{s_2}(s_2-s)^{\alpha-1}\left\|\psi(s)\right\|_{\mathbb{X}}\mathrm{d}s\nonumber\\&\le\frac{MR\alpha}{\Gamma(1+\alpha)}\left(\int_{0}^{s_1}\left|(s_2-s)^{\alpha-1}-(s_1-s)^{\alpha-1}\right|^2\mathrm{d}s\right)^{\frac{1}{2}}\nonumber\\&\quad+\int_{0}^{s_1}\!\!\!(s_1-s)^{\alpha-1}\left\|\widehat{\mathcal{T}}_{\alpha}(s_2-s)-\widehat{\mathcal{T}}_{\alpha}(s_1-s)\right\|_{\mathcal{L}(\mathbb{X})}\left\|\psi(s)\right\|_{\mathbb{X}}\mathrm{d}s+\frac{MR\alpha}{\Gamma(1+\alpha)}\left(\frac{(s_2-s_1)^{2\alpha-1}}{2\alpha-1}\right)^{\frac{1}{2}}\nonumber.
		\end{align}
		If $s_1=0,$ then by the above estimate, we deduce that 
		\begin{align*}
		\lim_{s_2\to 0^+}\left\|(\mathrm{Q}\psi)(s_2)-(\mathrm{Q}\psi)(s_1)\right\|_{\mathbb{X}}=0,\; \mbox{ unifromly for } \ \psi\in\mathrm{L}^2(J;\mathbb{X}).
		\end{align*} 
		For $0<\epsilon<s_1<T$, we have
		\begin{align}\label{2.8}
		&\left\|(\mathrm{Q}\psi)(s_2)-(\mathrm{Q}\psi)(s_1)\right\|_{\mathbb{X}}\nonumber\\&\le\frac{MR\alpha}{\Gamma(1+\alpha)}\left(\int_{0}^{s_1}\left|(s_2-s)^{\alpha-1}-(s_1-s)^{\alpha-1}\right|^2\mathrm{d}s\right)^{\frac{1}{2}}\nonumber\\&\quad+\int_{0}^{s_1-\epsilon}(s_1-s)^{\alpha-1}\left\|\widehat{\mathcal{T}}_{\alpha}(s_2-s)-\widehat{\mathcal{T}}_{\alpha}(s_1-s)\right\|_{\mathcal{L}(\mathbb{X})}\left\|\psi(s)\right\|_{\mathbb{X}}\mathrm{d}s\nonumber\\&\quad+\int_{s_1-\epsilon}^{s_1}(s_1-s)^{\alpha-1}\left\|\widehat{\mathcal{T}}_{\alpha}(s_2-s)-\widehat{\mathcal{T}}_{\alpha}(s_1-s)\right\|_{\mathcal{L}(\mathbb{X})}\left\|\psi(s)\right\|_{\mathbb{X}}\mathrm{d}s+\frac{MR\alpha}{\Gamma(1+\alpha)}\left(\frac{(s_2-s_1)^{2\alpha-1}}{2\alpha-1}\right)^{\frac{1}{2}}\nonumber\\&\le\frac{MR\alpha}{\Gamma(1+\alpha)}\left(\int_{0}^{s_1}\left|(s_2-s)^{\alpha-1}-(s_1-s)^{\alpha-1}\right|^2\mathrm{d}s\right)^{\frac{1}{2}}\nonumber\\&\quad+\sup_{s\in[0,s_1-\epsilon]}\left\|\widehat{\mathcal{T}}_{\alpha}(s_2-s)-\widehat{\mathcal{T}}_{\alpha}(s_1-s)\right\|_{\mathcal{L}(\mathbb{X})}\int_{0}^{s_1-\epsilon}(s_1-s)^{\alpha-1}\left\|\psi(s)\right\|_{\mathbb{X}}\mathrm{d}s\nonumber\\&\quad+\frac{2M\alpha}{\Gamma(1+\alpha)}\int_{s_1-\epsilon}^{s_1}(s_1-s)^{\alpha-1}\left\|\psi(s)\right\|_{\mathbb{X}}\mathrm{d}s+\frac{MR\alpha}{\Gamma(1+\alpha)}\left(\frac{(s_2-s_1)^{2\alpha-1}}{2\alpha-1}\right)^{\frac{1}{2}}\nonumber\\&\le\frac{MR\alpha}{\Gamma(1+\alpha)}\left(\int_{0}^{s_1}\left|(s_2-s)^{\alpha-1}-(s_1-s)^{\alpha-1}\right|^2\mathrm{d}s\right)^{\frac{1}{2}}\nonumber\\&\quad+R\sup_{s\in[0,s_1-\epsilon]}\left\|\widehat{\mathcal{T}}_{\alpha}(s_2-s)-\widehat{\mathcal{T}}_{\alpha}(s_1-s)\right\|_{\mathcal{L}(\mathbb{X})}\left(\frac{s_1^{2\alpha-1}-\epsilon^{2\alpha-1}}{2\alpha-1}\right)^{\frac{1}{2}}\nonumber\\&\quad+\frac{2MR\alpha}{\Gamma(1+\alpha)}\left(\frac{\epsilon^{2\alpha-1}}{2\alpha-1}\right)^{\frac{1}{2}}+\frac{MR\alpha}{\Gamma(1+\alpha)}\left(\frac{(s_2-s_1)^{2\alpha-1}}{2\alpha-1}\right)^{\frac{1}{2}}.
		\end{align}
		From Lemma \ref{lem2.5}, we know that the operator $\widehat{\mathcal{T}}_{\alpha}(t)$ is compact for $t>0$, which implies that the operator $\widehat{\mathcal{T}}_{\alpha}(t)$ is continuous under the uniform operator topology (see Theorem 3.2, Chapter 2, \cite{P}). Hence, using the arbitrariness of $\epsilon$ and continuity of  $\widehat{\mathcal{T}}_{\alpha}(t)$ in the uniform operator topology, the right hand side of the expression \eqref{2.8} converges to zero as $|s_2-s_1| \rightarrow 0$. Thus, $\mathrm{Q}\mathcal{B}_R$ is equicontinuous on $\mathrm{L}^{2}(J;\mathbb{X})$. 
		
		Next, we show that $\mathrm{V}(t):= \left\{(\mathrm{Q}\psi)(t):\psi\in \mathcal{B}_R\right\},$  for all $t\in J$ is relatively compact. For $t=0$, it is easy to check that the set $\mathrm{V}(t)$ is relatively compact in $\mathbb{X}$.  Let us take $ 0<t\leq T$ be fixed and for given $\eta$ with $ 0<\eta<t$ and any $\delta>0$, we define
		\begin{align*}
		(\mathrm{Q}^{\eta,\delta}\psi)(t)&=\alpha\int_{0}^{t-\eta}\int_{\delta}^{\infty}\xi(t-s)^{\alpha-1}\varphi_{\alpha}(\xi)\mathcal{T}(t^{\alpha}\xi)\psi(s)\mathrm{d}\xi\mathrm{d}s\nonumber\\&=\mathcal{T}(\eta^{\alpha}\delta)\alpha\int_{0}^{t-\eta}\int_{\delta}^{\infty}\xi(t-s)^{\alpha-1}\varphi_{\alpha}(\xi)\mathcal{T}(t^{\alpha}\xi-\eta^{\alpha}\delta)\psi(s)\mathrm{d}\xi\mathrm{d}s.
		\end{align*}
		Since the operator $\mathcal{T}(\cdot)$ is compact,  the set  $\mathrm{V}_{\eta,\delta}(t)=\{(\mathrm{Q}^{\eta,\delta}_\lambda \psi)(t):\psi\in \mathcal{B}_R\}$ is relatively compact in $\mathbb{X}$. Hence, there exist a finite $ x_{i}$'s, for $i=1,\dots, n $ in $ \mathbb{X} $ such that 
		\begin{align*}
		\mathrm{V}_{\eta,\delta}(t) \subset \bigcup_{i=1}^{n}\mathcal{S}(x_i, \varepsilon/2),
		\end{align*}
		for some $\varepsilon>0$, where $\mathcal{S}(x_i, \varepsilon/2)$ is an open ball centered at $x_i$ and of radius $\varepsilon/2$. Let us choose $\delta>0$ and $\eta>0$ such that 
		\begin{align*}
		\left\|(\mathrm{Q}\psi)(t)-(\mathrm{Q}^{\eta,\delta}\psi)(t)\right\|_{\mathbb{X}}&\le\alpha\left\|\int_{0}^{t}\int_{0}^{\delta}\xi(t-s)^{\alpha-1}\varphi_{\alpha}(\xi)\mathcal{T}(t^{\alpha}\xi-\eta^{\alpha}\delta)\psi(s)\mathrm{d}\xi\mathrm{d}s\right\|_{\mathbb{X}}\nonumber\\&\quad+\alpha\left\|\int_{t-\eta}^{t}\int_{\delta}^{\infty}\xi(t-s)^{\alpha-1}\varphi_{\alpha}(\xi)\mathcal{T}(t^{\alpha}\xi-\eta^{\alpha}\delta)\psi(s)\mathrm{d}\xi\mathrm{d}s\right\|_{\mathbb{X}}\nonumber\\&\le\frac{MR\alpha t^{2\alpha-1}}{2\alpha-1}\int_{0}^{\delta}\xi\varphi(\xi)\mathrm{d}\xi+\frac{MR\alpha}{\Gamma(1+\alpha)}\frac{\eta^{2\alpha-1}}{2\alpha-1}\le\frac{\varepsilon}{2}.
		\end{align*}
		Consequently $$ \mathrm{V}(t)\subset \bigcup_{i=1}^{n}\mathcal{S}(x_i, \varepsilon ). $$ Thus, for each $t\in J$, the set $\mathrm{V}(t)$ is relatively compact in $ \mathbb{X}$. Then by invoking the  Arzela-Ascoli theorem, we conclude that the operator $\mathrm{Q}$ is compact.
	\end{proof}
	Since $\mathbb{X}^*$ is strictly convex, the norm $\|\cdot\|_{\mathbb{X}}$ is Gateaux differentiable (cf. Fact 8.12, \cite{MFb}). Furthermore, every separable Banach space admits an equivalent Gateaux differentiable norm (cf. Theorem 8.13, \cite{MFb}). Since $\mathcal{J}$ is single-valued,  the Gateaux derivative of $\phi(x)=\frac{1}{2}\|x\|_{\mathbb{X}}^2$ is the duality map, that is, $$\langle\partial_x\phi(x),y\rangle=\lim_{\varepsilon \to 0}\frac{\phi(x+\varepsilon y)-\phi(x)}{\varepsilon}=\frac{1}{2}\frac{\mathrm{d}}{\mathrm{d}\varepsilon}\|x+\varepsilon y\|_{\mathbb{X}}^2\Big|_{\varepsilon=0}=\langle\mathcal{J}[x],y\rangle,$$ for $y\in\mathbb{X}$, where $\partial_x\phi(x)$ denotes the Gateaux derivative of $\phi$ at $x\in\mathbb{X}$. In fact, since $\mathbb{U}$ is a separable Hilbert space (identified with its own dual), by Theorem 8.24, \cite{MFb}, we infer that $\mathbb{U}$ admits a Fr\'echet differentiable norm.  The explicit expression of the optimal control ${u}$ in the feedback form is obtained in the following lemma: 
	
	\begin{lem}\label{lem3.1}
		Let $u$ be the optimal control satisfying \eqref{3.3} and minimizing the cost functional \eqref{3.1}. Then $u$ is given by 
		\begin{align*}
		u(t)=(T-t)^{\alpha-1}\mathrm{B}^*\widehat{\mathcal{T}}_{\alpha}(T-t)^*\mathcal{J}\left[\mathcal{R}(\lambda,\Phi_{0}^T)p(x(\cdot))\right],\  t\in [0, T),\ \lambda>0,\ \frac{1}{2}<\alpha<1, 
		\end{align*}
		with
		\begin{align*}
		p(x(\cdot))=x_{T}-\mathcal{T}_{\alpha}(T)\zeta.
		\end{align*}
	\end{lem}
	\begin{proof}
		Let us first consider the functional 
		\begin{align*}
		\mathcal{I}(\varepsilon)=\mathcal{G}(x_{u+\varepsilon w},u+\varepsilon w),
		\end{align*}
		where $(x,u)$ is the optimal solution of \eqref{3.3} and $w\in \mathrm{L}^{2}(J;\mathbb{U})$. Also the function  $x_{u+\varepsilon w}$ is the unique mild solution of \eqref{3.2} corresponding to the control $u+\varepsilon w$. Then, it is immediate that 
		\begin{align*}
		x_{u+\varepsilon w}(t)= \mathcal{T}_{\alpha}(t)\zeta+\int_{0}^{t}(t-s)^{\alpha-1}\widehat{\mathcal{T}}_{\alpha}(t-s)\mathrm{B}(u+\varepsilon w)(s)\mathrm{d}s.
		\end{align*}
		It is clear $\varepsilon=0$ is the critical point of $\mathcal{I}(\varepsilon)$. We now evaluate the first variation of the cost functional $\mathcal{G}$ (defined in \eqref{3.1}) as
		\begin{align*}
		\frac{\mathrm{d}}{\mathrm{d}\varepsilon}\mathcal{I}(\varepsilon)\Big|_{\varepsilon=0}&=\frac{\mathrm{d}}{\mathrm{d}\varepsilon}\bigg[\left\|x_{u+\varepsilon w}(T)-x_{T}\right\|^{2}_{\mathbb{X}}+\lambda\int^{T}_{0}\left\|u(t)+\varepsilon w(t)\right\|^{2}_{\mathbb{U}}\mathrm{d}t\bigg]_{\varepsilon=0}\nonumber\\
		&=2\bigg[\langle \mathcal{J}(x_{u+\varepsilon w}(T)-x_{T}), \frac{\mathrm{d}}{\mathrm{d}\varepsilon}(x_{u+\varepsilon w}(T)-x_{T})\rangle\nonumber\\&\qquad +2\lambda\int^{T}_{0}(u(t)+\varepsilon w(t),\frac{\mathrm{d}}{\mathrm{d}\varepsilon}(u(t)+\varepsilon w(t)))\mathrm{d}t\bigg]_{\varepsilon=0}\nonumber\\
		&=2\left\langle\mathcal{J}(x(T)-x_T),\int_0^T(T-t)^{\alpha-1}\widehat{\mathcal{T}}_{\alpha}(T-t)\mathrm{B}w(t)\mathrm{d}t \right\rangle+2\lambda\int_0^T(u(t),w(t))\mathrm{d} t. 
		\end{align*}
		By taking the first variation of the cost functional is zero, we deduce that
		\begin{align}\label{3.9}
		0&=\left\langle\mathcal{J}(x(T)-x_T),\int_0^T(T-t)^{\alpha-1}\widehat{\mathcal{T}}_{\alpha}(T-t)\mathrm{B}w(t)\mathrm{d}t\right\rangle+\lambda\int_0^T(u(t),w(t))\mathrm{d} t\nonumber\\&=\int_0^T(T-t)^{\alpha-1}\left\langle\mathcal{J}(x(T)-x_T),\widehat{\mathcal{T}}_{\alpha}(T-t)\mathrm{B}w(t) \right\rangle\mathrm{d}t+\lambda\int_0^T(u(t),w(t))\mathrm{d} t\nonumber\\&= \int_0^T\left((T-t)^{\alpha-1}\mathrm{B}^*\widehat{\mathcal{T}}_{\alpha}(T-t)^*\mathcal{J}(x(T)-x_T)+\lambda u(t),w(t) \right)\mathrm{d}t,
		\end{align}
		where $(\cdot,\cdot)$ is the inner product in the Hilbert space $\mathbb{U}$. Since $w\in \mathrm{L}^{2}(J;\mathbb{U})$ is an arbitrary element (one can choose $w$ to be $(T-t)^{\alpha-1}\mathrm{B}^*\widehat{\mathcal{T}}_{\alpha}(T-t)^*\mathcal{J}(x(T)-x_T)+\lambda u(t)$), it follows that the optimal control is given by
		\begin{align}\label{3.10}
		u(t)&= -\lambda^{-1}(T-t)^{\alpha-1}\mathrm{B}^*\widehat{\mathcal{T}}_{\alpha}(T-t)^*\mathcal{J}(x(T)-x_T),
		\end{align}
		for a.e. $t\in [0,T]$. Since by the relations \eqref{3.9} and \eqref{3.10}, it is clear that $u\in\mathrm{C}([0,T);\mathbb{U})$. Using the above expression of the control, we find
		%	It also holds for all $t\in[0,T)$, since from the expressions \eqref{3.9} and \eqref{3.10}, it is clear that $u$ is continuous and belongs to $\mathrm{C}([0,T);\mathbb{X})$. Therefore the state system \eqref{3.2} at a final point $T$ with the above control $u$ is given by
		\begin{align}\label{3.11}
		x(T)&=\mathcal{T}_{\alpha}(T)\zeta-\int^{T}_{0}\lambda^{-1}(T-s)^{2(\alpha-1)}\widehat{\mathcal{T}}_{\alpha}(T-s)\mathrm{B}\mathrm{B}^*\widehat{\mathcal{T}}_{\alpha}(T-s)^*\mathcal{J}(x(T)-x_T)\mathrm{d}s\nonumber\\
		&= \mathcal{T}_{\alpha}(T)\zeta-\lambda^{-1}\Phi_{0}^T\mathcal{J}\left[x(T)-x_{T}\right].
		\end{align}
		Let us assume
		\begin{align}\label{3.12}
		p(x(\cdot)):=x_{T}-\mathcal{T}_{\alpha}(T)\zeta.
		\end{align}
		Combining \eqref{3.11} and \eqref{3.12}, we have
		\begin{align}\label{3.13}
		x(T)-x_{T}&=-p(x(\cdot))-\lambda^{-1}\Phi_{0}^T\mathcal{J}\left[x(T)-x_{T}\right].
		\end{align}
		From \eqref{3.13}, one can easily deduce that 
		\begin{align}\label{3.15}
		x(T)-x_T=-\lambda\mathrm{I}(\lambda\mathrm{I}+\Phi_0^T\mathcal{J})^{-1}p(x(\cdot))=-\lambda\mathcal{R}(\lambda,\Phi_0^T)p(x(\cdot)).
		\end{align}
		Finally, from \eqref{3.10}, we get  the expression for optimal control  as
		\begin{align*}
		u(t)=(T-t)^{\alpha-1}\mathrm{B}^*\widehat{\mathcal{T}}_{\alpha}(T-t)^*\mathcal{J}\left[\mathcal{R}(\lambda,\Phi_{0}^T)p(x(\cdot))\right],\ \mbox{ for } \ t\in [0,T),
		\end{align*}
		which completes the proof. 
	\end{proof}
	Next, we examine the approximate controllability of the linear control system \eqref{3.2} through the following lemma.  
	\begin{lem}\label{lem4.2}
		The linear control system \eqref{3.2} is approximately controllable on $J$ if and only if Assumption (H0) holds. 
	\end{lem}
	A proof of the above lemma can be obtained by proceeding similarly as in the proof of  Theorem 3.2, \cite{SM}.
	\begin{rem}\label{rem3.4}
		If Assumption (\textit{H0}) holds, then by Theorem 2.3, \cite{M}, we know that the operator $\Phi_{0}^{T}$ is positive and vice versa. The positivity of $\Phi_{0}^{T}$ is equivalent to $$ \langle x^*, \Phi_{0}^{T}x^*\rangle=0\Rightarrow x^*=0.$$ We know that 
		\begin{align}
		\langle x^*, \Phi_{0}^{T}x^*\rangle =\int_0^T\left\|(T-t)^{\alpha-1}\mathrm{B}^*\widehat{\mathcal{T}}_\alpha(T-t)^*x^*\right\|_{\mathbb{X}^*}^2\mathrm{d}t.
		\end{align}
		By the above fact and  Lemma \ref{lem4.2}, we infer that the approximate controllability of the linear system \eqref{3.2} is equivalent to the condition $$\mathrm{B}^*\widehat{\mathcal{T}}_\alpha(T-t)^*x^*=0,\ 0\le t<T \Rightarrow x^*=0.$$
	\end{rem}

	\begin{rem}\label{rem3.6}
		Instead  \eqref{3.1}, one can consider the following cost functional also: 
		\begin{equation}\label{3.19}
		\mathcal{G}(x,u)=\left\|x(T)-x_{T}\right\|^{2}_{\mathbb{X}}+\lambda\int^{T}_{0}(T-t)^{\alpha-1}\left\|u(t)\right\|^{2}_{\mathbb{U}}\mathrm{d}t. 
		\end{equation}
		For $\frac{1}{2}<\alpha<1$,	the existence of optimal solution for the problem \eqref{3.3} follows similarly as in the proof of Theorem \ref{optimal}. For this, we have to replace \eqref{35} by 
		\begin{align*}
		T^{\alpha-1}\lambda\int^{T}_{0}\left\|u^n(t)\right\|^{2}_{\mathbb{U}}\mathrm{d}t\leq  &\left\|x^n(T)-x_{T}\right\|^{2}_{\mathbb{X}}+\lambda\int^{T}_{0}(T-t)^{\alpha-1}\left\|u^n(t)\right\|^{2}_{\mathbb{U}}\mathrm{d}t\nonumber\\&\leq \left\|x(T)-x_{T}\right\|^{2}_{\mathbb{X}}\leq 2\left(\|x(T)\|_{\mathbb{X}}^2+\|x_T\|_{\mathbb{X}}^2\right)<+\infty,
		\end{align*}
		so that \eqref{3.5} follows easily. Using the cost functional given in \eqref{3.19}, calculations similar to Lemma \ref{lem3.1} yields  the optimal control $u(\cdot)$ as 	\begin{align*}
		u(t)=\mathrm{B}^*\widehat{\mathcal{T}}_{\alpha}(T-t)^*\mathcal{J}\left[\mathcal{R}(\lambda,\Phi_{0}^T)p(x(\cdot))\right],\  t\in [0, T),\ \lambda>0,\ \frac{1}{2}<\alpha<1, 
		\end{align*}
		with
		$	p(x(\cdot))=x_{T}-\mathcal{T}_{\alpha}(T)\zeta,$ and 
		\begin{equation}\label{3.20}
		\Phi_{0}^{T}=\int^{T}_{0}(T-t)^{\alpha-1}\widehat{\mathcal{T}}_{\alpha}(T-t)\mathrm{B}\mathrm{B}^*\widehat{\mathcal{T}}_{\alpha}(T-t)^*\mathrm{d}t. 
		\end{equation}
	\end{rem}
	
	\section{Approximate Controllability of the Semilinear Impulsive System} \label{semilinear}\setcounter{equation}{0}
	The purpose of this section is to investigate the approximate controllability of the fractional order  semilinear impulsive system \eqref{1.1}. In order to acquire sufficient conditions on approximate controllability, we first show that for $\lambda>0$ and $x_T\in\mathbb{X}$, there exists a mild solution of the system \eqref{1.1} with the control function defined as 
	\begin{align}\label{C}
	u^\alpha_{\lambda}(t)&=\sum_{k=0}^{m}u^\alpha_{k,\lambda}(t)\chi_{[\tau_k, t_{k+1})}(t), \ t\in J,\  \frac{1}{2}<\alpha<1,
	\end{align}
	where 
	\begin{align*}
	u^\alpha_{k,\lambda}(t)&=(t_{k+1}-t)^{\alpha-1}\mathrm{B}^*\widehat{\mathcal{T}}_{\alpha}(t_{k+1}-t)^*\mathcal{J}\left[\mathcal{R}(\lambda,\Phi_{\tau_k}^{t_{k+1}})p_k(x(\cdot))\right],
	\end{align*}
	for $t\in [\tau_k, t_{k+1}),k=0,1,\ldots,m$, with
	\begin{align*}
	p_0(x(\cdot))&=\zeta_{0}-\mathcal{T}_{\alpha}(t_1)\psi(0)-\int^{t_1}_{0}(t_1-s)^{\alpha-1}\widehat{\mathcal{T}}_{\alpha}(t_1-s)f(s,\tilde{x}_{\rho(s,\tilde{x}_s)})\mathrm{d}s,\nonumber\\
	p_k(x(\cdot))&=\zeta_{k}-\mathcal{T}_{\alpha}(t_{k+1}-\tau_k)h_k(\tau_k,\tilde{x}(t_k^-))\\&\quad+\int_{0}^{\tau_k}(\tau_k-s)^{\alpha-1}\widehat{\mathcal{T}}_{\alpha}(\tau_k-s)\left[f(s,\tilde{x}_{\rho(s,\tilde{x}_s)})+\mathrm{B}\sum_{j=0}^{k-1}u^\alpha_{j,\lambda}(s)\chi_{[\tau_j, t_{j+1})}(s)\right]\mathrm{d}s \\&\quad-\int_{0}^{t_{k+1}}(t_{k+1}-s)^{\alpha-1}\widehat{\mathcal{T}}_{\alpha}(t_{k+1}-s)f(s,\tilde{x}_{\rho(s,\tilde{x}_s)})\mathrm{d}s	\\&\quad-\int_{0}^{\tau_k}(t_{k+1}-s)^{\alpha-1}\widehat{\mathcal{T}}_{\alpha}(t_{k+1}-s)\mathrm{B}\sum_{j=0}^{k-1}u^\alpha_{j,\lambda}(s)\chi_{[\tau_j, t_{j+1})}(s)\mathrm{d}s, \ k=1,\ldots,m,
	%p_m(x(\cdot))&=x_{T}-\mathcal{T}_{\alpha}(T-\tau_m)g(\tau_m,\tilde{x}(t_m^-))-\int_{\tau_m}^{T}(T-s)^{\alpha-1}\widehat{\mathcal{T}}_{\alpha}(T-s)f(s,\tilde{x}_{\rho(s,\tilde{x}_s)})\mathrm{d}s,
	\end{align*} 
	and  $\tilde{x}:(-\infty,T]\rightarrow\mathbb{X}$  such that $\tilde{x}(t)=\psi(t), \ t\in(-\infty,0]\ \tilde{x}(t)=x(t),\ t\in J=[0,T],$ and $\zeta_{k}\in \mathbb{X}$ for $k=0,1,\ldots,m$.
	\begin{rem}
		Since the operator $\Phi_{\tau_k}^{t_{k+1}},$  for each $k=0,\ldots,m,$ is non-negative, linear and bounded for $\frac{1}{2}<\alpha<1$, Lemma \ref{lem2.9} is also valid for each $\Phi_{\tau_k}^{t_{k+1}},$  for $k=0,\ldots,m$.
	\end{rem}
	The following theorem provides the existence of mild solution of the system \eqref{1.1} with the control \eqref{C}. 

	\begin{theorem}\label{thm4.3}
		Let Assumptions (H1)-(H3) hold true. Then for every $ \lambda>0 $ and  fixed $\zeta_{k}\in\mathbb{X},$ for $k=0,1,\ldots,m$, the system  \eqref{1.1} with the control \eqref{C} has at least one mild solution on $J$,  provided 
		\begin{align}\label{cnd}
	\frac{MH_{2}\alpha\beta}{\Gamma(1+\alpha)}\frac{2T^{\alpha-\alpha_1}}{\mu^{1-\alpha_1}}\left\{1+\frac{(m+1)(m+2)\tilde{R}}{2}+\frac{m(m+1)\tilde{R}^2}{2}\sum_{j=0}^{m-1}e^{\frac{(m+j)(m-j-1)\tilde{R}}{2}}\right\}<1,
		\end{align}
		where $\tilde{R}=\left(\frac{M\tilde{M}\alpha}{\Gamma(1+\alpha)}\right)^{2}\frac{2T^{2\alpha-1}}{\lambda(2\alpha-1)}$ and $\mu=\frac{\alpha-\alpha_1}{1-\alpha_1}$.
	\end{theorem}
	\begin{proof}
		Let us take a set $\mathrm{E}:=\{x\in\mathrm{PC}(J;\mathbb{X}) : x(0)=\psi(0)\}$ with the norm $\left\|\cdot\right\|_{\mathrm{PC}(J;\mathbb{X})}$. For each $r>0$, we consider a set $\mathrm{E}_{r}=\{x\in\mathrm{E} : \left\|x\right\|_{\mathrm{PC}(J;\mathbb{X})}\le r\}$.
		
		For $\lambda>0$, let us define an operator $F_{\lambda}:\mathrm{E}\to\mathrm{E}$ such that
		\begin{eqnarray}\label{2}
		(F_{\lambda}x)(t)=\left\{
		\begin{aligned}
		&\mathcal{T}_{\alpha}(t)\psi(0)+\int_{0}^{t}(t-s)^{\alpha-1}\widehat{\mathcal{T}}_{\alpha}(t-s)\left[\mathrm{B}u^{\alpha}_{\lambda}(s)+f(s,\tilde{x}_{\rho(s, \tilde{x}_s)})\right]\mathrm{d}s,\ t\in[0, t_1],\\
		&	h_k(t, \tilde{x}(t_k^-)),\ t\in(t_k, \tau_k],\ k=1,\ldots,m,\\
		&\mathcal{T}_{\alpha}(t-\tau_k)h_k(\tau_k,\tilde{x}(t_k^-))-\int_{0}^{\tau_k}(\tau_k-s)^{\alpha-1}\widehat{\mathcal{T}}_{\alpha}(\tau_k-s)\left[\mathrm{B}u^{\alpha}_{\lambda}(s)+f(s,\tilde{x}_{\rho(s, \tilde{x}_s)})\right]\mathrm{d}s\\&\quad+\int_{0}^{t}(t-s)^{\alpha-1}\widehat{\mathcal{T}}_{\alpha}(t-s) \left[\mathrm{B}u^{\alpha}_{\lambda}(s)+f(s,\tilde{x}_{
			\rho(s, \tilde{x}_s)})\right]\mathrm{d}s,\ t\in(\tau_k,t_{k+1}],\ k=1,\ldots,m,
		\end{aligned}
		\right.
		\end{eqnarray}
		where $u^{\alpha}_{\lambda}$ is given in \eqref{C}. From the definition of $ F_{\lambda}$, we infer that the system $\eqref{1.1}$ has a mild solution, if the operator $ F_{\lambda}$ has a fixed point. We divide the proof of the fact that the operator $F_{\lambda}$ has a fixed point in the following steps. 
		\vskip 0.1in 
		\noindent\textbf{Step (1): } \emph{$ F_{\lambda}(\mathrm{E}_r)\subset \mathrm{E}_r,$ for some $ r $}. On the contrary, let us suppose that our claim is not true. Then for any $\lambda>0$ and for all $r>0$, there exists $x^r\in \mathrm{E}_r,$ such that $\left\|(F_{\lambda}x^r)(t)\right\|_\mathbb{X}>r,$ for some $t\in J$, where $t$ may depend upon $r$. 
		First, by using Assumption \ref{as2.1}, we estimate 
		\begin{align}\label{4.4}
		\left\|p_0(x(\cdot))\right\|_{\mathbb{X}}&\le\left\|\zeta_0\right\|_{\mathbb{X}}+\left\|\mathcal{T}_{\alpha}(t_1)\psi(0)\right\|_{\mathbb{X}}+\int_{0}^{t_1}(t_1-s)^{\alpha-1}\left\|\widehat{\mathcal{T}}_{\alpha}(t_1-s)f(s,\tilde{x}_{\rho(s,\tilde{x}_s)})\right\|_{\mathbb{X}}\mathrm{d}s\nonumber\\&\le\left\|\zeta_0\right\|_{\mathbb{X}}+M\left\|\psi(0)\right\|_{\mathbb{X}}+\frac{M\alpha}{\Gamma(1+\alpha)}\int_{0}^{t_1}(t_1-s)^{\alpha-1}\gamma_{r'}(s)\mathrm{d}s\nonumber\\&\le\left\|\zeta_0\right\|_{\mathbb{X}}+M\left\|\psi(0)\right\|_{\mathbb{X}}+\frac{M\alpha}{\Gamma(1+\alpha)}\left(\int_{0}^{t_1}(t_1-s)^{\frac{\alpha-1}{1-\alpha_1}}\mathrm{d}s\right)^{1-\alpha_1}\left(\int_{0}^{t_1}(\gamma_{r'}(s))^{\frac{1}{\alpha_1}}\mathrm{d}s\right)^{\alpha_1}\nonumber\\&\le\left\|\zeta_0\right\|_{\mathbb{X}}+M\left\|\psi(0)\right\|_{\mathbb{X}}+\frac{M\alpha}{\Gamma(1+\alpha)}\frac{t_{1}^{\alpha-\alpha_1}}{\mu^{1-\alpha_1}}\left\|\gamma_{r'}\right\|_{\mathcal{L}^{\frac{1}{\alpha_1}}([0,t_1];\mathbb{R^{+}})}\nonumber\\&\le\left\|\zeta_0\right\|_{\mathbb{X}}+M\left\|\psi(0)\right\|_{\mathbb{X}}+\frac{M\alpha}{\Gamma(1+\alpha)}\frac{T^{\alpha-\alpha_1}}{\mu^{1-\alpha_1}}\left\|\gamma_{r'}\right\|_{\mathcal{L}^{\frac{1}{\alpha_1}}(J;\mathbb{R^{+}})}\nonumber\\&\le\left\|\zeta_0\right\|_{\mathbb{X}}+M\left\|\psi(0)\right\|_{\mathbb{X}}+\frac{M\alpha}{\Gamma(1+\alpha)}\frac{2T^{\alpha-\alpha_1}}{\mu^{1-\alpha_1}}\left\|\gamma_{r'}\right\|_{\mathcal{L}^{\frac{1}{\alpha_1}}(J;\mathbb{R^{+}})}=N_0,
		\end{align} 
		where $\mu=\frac{\alpha-\alpha_1}{1-\alpha_1}$ and $r'=H_{1}\left\|\psi\right\|_{\mathfrak{B}}+H_{2}r$.	Further, using Assumption \ref{as2.1}, we estimate 
		\begin{align*}
		&\left\|p_k(x(\cdot))\right\|_{\mathbb{X}}\nonumber\\&\le\left\|\zeta_k\right\|_{\mathbb{X}}+\left\|\mathcal{T}_{\alpha}(t-\tau_k)h_k(\tau_k,\tilde{x}(t_k^-))\right\|_{\mathbb{X}}+\int_{0}^{\tau_k}(\tau_k-s)^{\alpha-1}\left\|\widehat{\mathcal{T}}_{\alpha}(\tau_k-s)f(s,\tilde{x}_{\rho(s,\tilde{x}_s)})\right\|_{\mathbb{X}}\mathrm{d}s\nonumber\\&\quad+\int_{0}^{t_{k+1}}(t_{k+1}-s)^{\alpha-1}\left\|\widehat{\mathcal{T}}_{\alpha}(t_{k+1}-s)f(s,\tilde{x}_{\rho(s,\tilde{x}_s)})\right\|_{\mathbb{X}}\mathrm{d}s\nonumber\\&\quad+\int_{0}^{\tau_k}(\tau_k-s)^{\alpha-1}\left\|\widehat{\mathcal{T}}_{\alpha}(\tau_k-s)\mathrm{B}\sum_{j=0}^{k-1}u^\alpha_{j,\lambda}(s)\chi_{[\tau_j,t_{j+1})}(s)\right\|_{\mathbb{X}}\mathrm{d}s\nonumber\\&\quad+\int_{0}^{\tau_k}(t_{k+1}-s)^{\alpha-1}\left\|\widehat{\mathcal{T}}_{\alpha}(t_{k+1}-s)\mathrm{B}\sum_{j=0}^{k-1}u^\alpha_{j,\lambda}(s)\chi_{[\tau_j, t_{j+1})}(s)\right\|_{\mathbb{X}}\mathrm{d}s\nonumber\\&\le\left\|\zeta_k\right\|_{\mathbb{X}}+Ml_k+\frac{M\alpha}{\Gamma(1+\alpha)}\frac{\tau_k^{\alpha-\alpha_1}}{\mu^{1-\alpha_1}}\left\|\gamma_{r'}\right\|_{\mathcal{L}^{\frac{1}{\alpha_1}}([0,\tau_k];\mathbb{R^{+}})}+\frac{M\alpha}{\Gamma(1+\alpha)}\frac{t_{k+1}^{\alpha-\alpha_1}}{\mu^{1-\alpha_1}}\left\|\gamma_{r'}\right\|_{\mathcal{L}^{\frac{1}{\alpha_1}}([0,t_{k+1}];\mathbb{R^{+}})}\nonumber\\&\quad+\frac{1}{\lambda}\left(\frac{M\tilde{M}\alpha}{\Gamma(1+\alpha)}\right)^{2}\sum_{j=0}^{k-1}\left\|p_j(x(\cdot))\right\|_{\mathbb{X}}\int_{\tau_j}^{t_{j+1}}(\tau_k-s)^{\alpha-1}(t_{j+1}-s)^{\alpha-1}\mathrm{d}s\nonumber\\&\quad+\frac{1}{\lambda}\left(\frac{M\tilde{M}\alpha}{\Gamma(1+\alpha)}\right)^{2}\sum_{j=0}^{k-1}\left\|p_j(x(\cdot))\right\|_{\mathbb{X}}\int_{\tau_j}^{t_{j+1}}(t_{k+1}-s)^{\alpha-1}(t_{j+1}-s)^{\alpha-1}\mathrm{d}s\nonumber\\&\le\left\|\zeta_k\right\|_{\mathbb{X}}+Ml_k+\frac{2M\alpha}{\Gamma(1+\alpha)}\frac{T^{\alpha-\alpha_1}}{\mu^{1-\alpha_1}}\left\|\gamma_{r'}\right\|_{\mathcal{L}^{\frac{1}{\alpha_1}}(J;\mathbb{R^{+}})}\nonumber\\&\quad+\frac{2}{\lambda}\left(\frac{M\tilde{M}\alpha}{\Gamma(1+\alpha)}\right)^{2}\sum_{j=0}^{k-1}\left\|p_j(x(\cdot))\right\|_{\mathbb{X}}\int_{\tau_j}^{t_{j+1}}(t_{j+1}-s)^{2(\alpha-1)}\mathrm{d}s\nonumber\\&\le N_k+\frac{2}{\lambda}\left(\frac{M\tilde{M}\alpha}{\Gamma(1+\alpha)}\right)^{2}\sum_{j=0}^{k-1}\left\|p_j(x(\cdot))\right\|_{\mathbb{X}}\frac{(t_{j+1}-\tau_j)^{2\alpha-1}}{2\alpha-1}\nonumber\\&\le N_k+\frac{2T^{2\alpha-1}}{\lambda(2\alpha-1)}\left(\frac{M\tilde{M}\alpha}{\Gamma(1+\alpha)}\right)^{2}\sum_{j=0}^{k-1}\left\|p_j(x(\cdot))\right\|_{\mathbb{X}}
		\nonumber\\&= N_k+\tilde{R}\sum_{j=0}^{k-1}\left\|p_j(x(\cdot))\right\|_{\mathbb{X}},
		\end{align*}
		where  $\tilde{R}=\left(\frac{M\tilde{M}\alpha}{\Gamma(1+\alpha)}\right)^{2}\frac{2T^{2\alpha-1}}{\lambda(2\alpha-1)}$ and $N_k=\left\|\zeta_k\right\|_{\mathbb{X}}+Ml_k+\frac{2T^{\alpha-\alpha_1}}{\mu^{1-\alpha_1}}\frac{M\alpha}{\Gamma(1+\alpha)}\left\|\gamma_{r'}\right\|_{\mathrm{L}^{\frac{1}{\alpha_1}}(J;\mathbb{R^+})}$ for $k=1,\ldots,m.$ Applying the discrete Gronwall-Bellman lemma (Lemma \ref{lem2.13}), we obtain 
		\begin{align}\label{4.5}
		\left\|p_k(x(\cdot))\right\|_{\mathbb{X}}\le&N_k+\tilde{R}\sum_{j=0}^{k-1}N_je^{\frac{(k+j)(k-j-1)\tilde{R}}{2}}=C_k, \ \mbox{for}\ k=1,\ldots,m.
		\end{align}
		Taking $t\in[0,t_1]$ and using the relation \eqref{2.5}, Lemma \ref{lem2.5} and Assumption \ref{as2.1} \textit{(H1)}-\textit{(H3)}, we compute
		\begin{align}\label{4.21}
		r&<\left\|(F_{\lambda}x^r)(t)\right\|_\mathbb{X}\nonumber\\&=\left\|\mathcal{T}_{\alpha}(t)\psi(0)+\int_{0}^{t}(t-s)^{\alpha-1}\widehat{\mathcal{T}}_{\alpha}(t-s)\left[\mathrm{B}u^{\alpha}_{\lambda}(s)+f(s,\tilde{x}_{\rho(s, \tilde{x}_s)})\right]\mathrm{d}s\right\|_{\mathbb{X}} \nonumber\\&\le\left\|\mathcal{T}_{\alpha}(t)\psi(0)\right\|_{\mathbb{X}}+\int_0^t(t-s)^{\alpha-1}\left\|\widehat{\mathcal{T}}_{\alpha}(t-s)\mathrm{B}u^{\alpha}_{\lambda}(s)\right\|_{\mathbb{X}}\mathrm{d}s\nonumber\\&\quad+\int_0^t(t-s)^{\alpha-1}\left\|\widehat{\mathcal{T}}_{\alpha}(t-s)f(s,\tilde{x}_{\rho(s,\tilde{x}_s)})\right\|_{\mathbb{X}}\mathrm{d}s\nonumber\\&\le  M\left\|\psi(0)\right\|_{\mathbb{X}}+\frac{M\tilde{M}\alpha}{\Gamma(1+\alpha)}\!\!\int_0^t\!\!(t-s)^{\alpha-1}\left\|u^{\alpha}_{\lambda}(s)\right\|_{\mathbb{U}}\mathrm{d}s+\frac{M\alpha}{\Gamma(1+\alpha)}\!\!\int_0^t\!\!(t-s)^{\alpha-1}\left\|f(s,\tilde{x}_{\rho(s,\tilde{x}_s)})\right\|_{\mathbb{X}}\mathrm{d}s\nonumber\\&\le M\left\|\psi(0)\right\|_{\mathbb{X}}+\frac{1}{\lambda}\left(\frac{M\tilde{M}\alpha}{\Gamma(1+\alpha)}\right)^{2}\left\|p_0(x(\cdot))\right\|_{\mathbb{X}}\int_{0}^{t}(t-s)^{\alpha-1}(t_1-s)^{\alpha-1}\mathrm{d}s\nonumber\\&\quad+\frac{M\alpha}{\Gamma(1+\alpha)}\int_0^t(t-s)^{\alpha-1}\gamma_{r'}(s)\mathrm{d}s\nonumber\\&\le M\left\|\psi(0)\right\|_{\mathbb{X}}+\left(\frac{M\tilde{M}\alpha}{\Gamma(1+\alpha)}\right)^{2}\frac{t^{2\alpha-1}}{\lambda(2\alpha-1)}\left\|p_0(x(\cdot))\right\|_{\mathbb{X}}\nonumber\\&\quad+\frac{M\alpha}{\Gamma(1+\alpha)}\left(\int_0^t(t-s)^{\frac{\alpha-1}{1-\alpha_1}}\mathrm{d}s\right)^{1-\alpha_1}\left\|\gamma_{r'}\right\|_{\mathrm{L}^{\frac{1}{\alpha_1}}([0,t];\mathbb{R^{+}})}\nonumber\\&\le M\left\|\psi(0)\right\|_{\mathbb{X}}+\left(\frac{M\tilde{M}\alpha}{\Gamma(1+\alpha)}\right)^{2}\frac{T^{2\alpha-1}}{\lambda(2\alpha-1)}\left\|p_0(x(\cdot))\right\|_{\mathbb{X}}+\frac{M\alpha}{\Gamma(1+\alpha)}\frac{T^{\alpha-\alpha_1}}{\mu^{1-\alpha_1}}\left\|\gamma_{r'}\right\|_{\mathrm{L}^{\frac{1}{\alpha_1}}(J;\mathbb{R^+})}\nonumber\\&\le M\left\|\psi(0)\right\|_{\mathbb{X}}+\left(\frac{M\tilde{M}\alpha}{\Gamma(1+\alpha)}\right)^{2}\frac{2T^{2\alpha-1}}{\lambda(2\alpha-1)}\sum_{j=0}^{m}\left\|p_j(x(\cdot))\right\|_{\mathbb{X}}+\frac{M\alpha}{\Gamma(1+\alpha)}\frac{2T^{\alpha-\alpha_1}}{\mu^{1-\alpha_1}}\left\|\gamma_{r'}\right\|_{\mathrm{L}^{\frac{1}{\alpha_1}}(J;\mathbb{R^+})}\nonumber\\&\le M\left\|\psi(0)\right\|_{\mathbb{X}}+\tilde{R}\sum_{j=0}^{m}C_j+\frac{M\alpha}{\Gamma(1+\alpha)}\frac{2T^{\alpha-\alpha_1}}{\mu^{1-\alpha_1}}\left\|\gamma_{r'}\right\|_{\mathrm{L}^{\frac{1}{\alpha_1}}(J;\mathbb{R^+})},
		\end{align}
		where $C_0=N_0$. For $t\in(t_k,\tau_k],\ k=1,\ldots,m$, we obtain  
		\begin{align}\label{4.22}
		r<\left\|(F_{\lambda}x^r)(t)\right\|_\mathbb{X}&\le\left\|h_k(t, \tilde{x}(t_k^-))\right\|_{\mathbb{X}}\nonumber\\&\le l_k\le l_k+\tilde{R}\sum_{j=0}^{m}C_j+\frac{2T^{\alpha-\alpha_1}}{\mu^{1-\alpha_1}}\frac{M\alpha}{\Gamma(1+\alpha)}\left\|\gamma_{r'}\right\|_{\mathrm{L}^{\frac{1}{\alpha_1}}(J;\mathbb{R^+})}.
		\end{align}
		Taking $t\in(\tau_k,t_{k+1}], \ k=1,\dots,m$, we evaluate
		\begin{align}\label{4.23}
		r&<\left\|(F_{\lambda}x^r)(t)\right\|_\mathbb{X}\nonumber\\&=\bigg\|\mathcal{T}_{\alpha}(t-\tau_k)h_k(\tau_k,\tilde{x}(t_k^-))-\int_{0}^{\tau_k}(\tau_k-s)^{\alpha-1}\widehat{\mathcal{T}}_{\alpha}(\tau_k-s)\left[\mathrm{B}u^{\alpha}_{\lambda}(s)+f(s,\tilde{x}_{\rho(s,\tilde{x}_s)})\right]\mathrm{d}s\nonumber\\&\quad+\int_{0}^{t}(t-s)^{\alpha-1}\widehat {\mathcal{T}}_{\alpha}(t-s)\left[\mathrm{B}u^{\alpha}_{\lambda}(s)+f(s,\tilde{x}_{\rho(s,\tilde{x}_s)})\right]\mathrm{d}s\bigg\|_{\mathbb{X}}\nonumber\\&\le\left\|\mathcal{T}_{\alpha} (t)h_k(\tau_k, \tilde{x}(t_k^-))\right\|_{\mathbb{X}}+\int_{0}^{\tau_k}(\tau_k-s)^{\alpha-1}\left\|\widehat{\mathcal{T}}_{\alpha}(\tau_k-s)\mathrm{B}u^{\alpha}_{\lambda}(s)
		\right\|_{\mathbb{X}}\mathrm{d}s\nonumber\\&\quad+\int_{0}^{\tau_k}(\tau_k-s)^{\alpha-1}\left\|\widehat{\mathcal{T}}_{\alpha}(\tau_k-s)f(s,\tilde{x}_{\rho(s,\tilde{x}_s)})\right\|_{\mathbb{X}}\mathrm{d}s+\int_{0}^t(t-s)^{\alpha-1}\left\|\widehat{\mathcal{T}}_{\alpha}(t-s)\mathrm{B}u^{\alpha}_{\lambda}(s)\right\|_{\mathbb{X}}\mathrm{d}s\nonumber\\&\quad +\int_{0}^t(t-s)^{\alpha-1}\left\|\widehat{\mathcal{T}}_{\alpha}(t-s)f(s,\tilde{x}_{\rho(s, \tilde{x}_s)})\right\|_{\mathbb{X}}\mathrm{d}s
		\nonumber\\&\le Ml_k+\frac{M\tilde{M}\alpha}{\Gamma(1+\alpha)}\int_{0}^{\tau_k}(\tau_k-s)^{\alpha-1}\left\|u^{\alpha}_{\lambda}(s)\right\|_{\mathbb{U}}\mathrm{d}s+\frac{M\alpha}
		{\Gamma(1+\alpha)}\int_{0}^{\tau_k}(\tau_k-s)^{\alpha-1}\left\|f(s,\tilde{x}_{\rho(s, \tilde{x}_s)})\right\|_{\mathbb{X}}\mathrm{d}s\nonumber\\&\quad+\frac{M\tilde{M}\alpha}{\Gamma(1+\alpha)}\int_{0}^{t}(t-s)^{\alpha-1}\left\|u^{\alpha}_{\lambda}(s)\right\|_{\mathbb{U}}\mathrm{d}s+\frac{M\alpha}
		{\Gamma(1+\alpha)}\int_{0}^{t}(t-s)^{\alpha-1}\left\|f(s,\tilde{x}_{\rho(s, \tilde{x}_s)})\right\|_{\mathbb{X}}\mathrm{d}s\nonumber\\&\le Ml_k+\frac{1}{\lambda}\left(\frac{M\tilde{M}\alpha}{\Gamma(1+\alpha)}\right)^{2}\sum_{j=0}^{k-1}\left\|p_j(x(\cdot))\right\|_{\mathbb{X}}\int_{\tau_j}^{t_{j+1}}(\tau_k-s)^{\alpha-1}(t_{j+1}-s)^{\alpha-1}\mathrm{d}s\nonumber\\&\quad+\frac{M\alpha}{\Gamma(1+\alpha)}\left(\int_{0}^{\tau_k}(\tau_k-s)^{\frac{\alpha-1}{1-\alpha_1}}\mathrm{d}s\right)^{1-\alpha_1}\left\|\gamma_{r'}\right\|_{\mathcal{L}^{\frac{1}{\alpha_1}}([0,\tau_k];\mathbb{R^{+}})}\nonumber\\&\quad+\frac{1}{\lambda}\left(\frac{M\tilde{M}\alpha}{\Gamma(1+\alpha)}\right)^{2}\Bigg[\sum_{j=0}^{k-1}\left\|p_j(x(\cdot))\right\|_{\mathbb{X}}\int_{\tau_j}^{t_{j+1}}(t-s)^{\alpha-1}(t_{j+1}-s)^{\alpha-1}\mathrm{d}s\nonumber\\&\qquad+\left\|p_k(x(\cdot))\right\|_{\mathbb{X}}\int_{\tau_k}^{t}(t-s)^{\alpha-1}(t_{k+1}-s)^{\alpha-1}\mathrm{d}s\Bigg]+\frac{M\alpha}{\Gamma(1+\alpha)}\left(\int_{0}^{t}(t-s)^{\frac{\alpha-1}{1-\alpha_1}}\mathrm{d}s\right)^{1-\alpha_1}\left\|\gamma_{r'}\right\|_{\mathcal{L}^{\frac{1}{\alpha_1}}([0,t];\mathbb{R^{+}})}\nonumber\\&\le Ml_k+\frac{2}{\lambda}\left(\frac{M\tilde{M}\alpha}{\Gamma(1+\alpha)}\right)^{2}\sum_{j=0}^{k-1}\left\|p_j(x(\cdot))\right\|_{\mathbb{X}}\int_{\tau_j}^{t_{j+1}}(\tau_k-s)^{\alpha-1}(t_{j+1}-s)^{\alpha-1}\mathrm{d}s\nonumber\\&\quad+\frac{M\alpha}{\Gamma(1+\alpha)}\left(\int_{0}^{\tau_k}(\tau_k-s)^{\frac{\alpha-1}{1-\alpha_1}}\mathrm{d}s\right)^{1-\alpha_1}\left\|\gamma_{r'}\right\|_{\mathcal{L}^{\frac{1}{\alpha_1}}([0,\tau_k];\mathbb{R^{+}})}\nonumber\\&\quad+\frac{1}{\lambda}\left(\frac{M\tilde{M}\alpha}{\Gamma(1+\alpha)}\right)^{2}\left\|p_k(x(\cdot))\right\|_{\mathbb{X}}\int_{\tau_k}^{t}(t-s)^{\alpha-1}(t_{k+1}-s)^{\alpha-1}\mathrm{d}s\nonumber\\&\quad+\frac{M\alpha}{\Gamma(1+\alpha)}\left(\int_{0}^{t}(t-s)^{\frac{\alpha-1}{1-\alpha_1}}\mathrm{d}s\right)^{1-\alpha_1}\left\|\gamma_{r'}\right\|_{\mathcal{L}^{\frac{1}{\alpha_1}}([0,t];\mathbb{R^{+}})}\nonumber\\&\le Ml_k+\frac{2}{\lambda}\left(\frac{M\tilde{M}\alpha}{\Gamma(1+\alpha)}\right)^{2}\sum_{j=0}^{k-1}\left\|p_j(x(\cdot))\right\|_{\mathbb{X}}\int_{\tau_j}^{t_{j+1}}(t_{j+1}-s)^{2(\alpha-1)}\mathrm{d}s\nonumber\\&\quad+\frac{M\alpha}{\Gamma(1+\alpha)}\left(\int_{0}^{\tau_k}(\tau_k-s)^{\frac{\alpha-1}{1-\alpha_1}}\mathrm{d}s\right)^{1-\alpha_1}\left\|\gamma_{r'}\right\|_{\mathcal{L}^{\frac{1}{\alpha_1}}([0,\tau_k];\mathbb{R^{+}})}\nonumber\\&\quad+\frac{1}{\lambda}\left(\frac{M\tilde{M}\alpha}{\Gamma(1+\alpha)}\right)^{2}\left\|p_k(x(\cdot))\right\|_{\mathbb{X}}\int_{\tau_k}^{t}(t-s)^{2(\alpha-1)}\mathrm{d}s+\frac{M\alpha}{\Gamma(1+\alpha)}\left(\int_{0}^{t}(t-s)^{\frac{\alpha-1}{1-\alpha_1}}\mathrm{d}s\right)^{1-\alpha_1}\left\|\gamma_{r'}\right\|_{\mathcal{L}^{\frac{1}{\alpha_1}}([0,t];\mathbb{R^{+}})}\nonumber\\&\le Ml_k+\frac{2}{\lambda}\left(\frac{M\tilde{M}\alpha}{\Gamma(1+\alpha)}\right)^{2}\sum_{j=0}^{k-1}\left\|p_j(x(\cdot))\right\|_{\mathbb{X}}\frac{(t_{j+1}-\tau_j)^{2\alpha-1}}{2\alpha-1}\nonumber\\&\quad+\frac{M\alpha}{\Gamma(1+\alpha)}\left(\int_{0}^{\tau_k}(\tau_k-s)^{\frac{\alpha-1}{1-\alpha_1}}\mathrm{d}s\right)^{1-\alpha_1}\left\|\gamma_{r'}\right\|_{\mathcal{L}^{\frac{1}{\alpha_1}}([0,\tau_k];\mathbb{R^{+}})}\nonumber\\&\quad+\frac{1}{\lambda}\left(\frac{M\tilde{M}\alpha}{\Gamma(1+\alpha)}\right)^{2}\left\|p_k(x(\cdot))\right\|_{\mathbb{X}}\frac{(t-\tau_k)^{2\alpha-1}}{2\alpha-1}+\frac{M\alpha}{\Gamma(1+\alpha)}\left(\int_{0}^{t}(t-s)^{\frac{\alpha-1}{1-\alpha_1}}\mathrm{d}s\right)^{1-\alpha_1}\left\|\gamma_{r'}\right\|_{\mathcal{L}^{\frac{1}{\alpha_1}}([0,t];\mathbb{R^{+}})}\nonumber\\&\le Ml_k+\frac{2T^{2\alpha-1}}{\lambda(2\alpha-1)}\left(\frac{M\tilde{M}\alpha}{\Gamma(1+\alpha)}\right)^{2}\sum_{j=0}^{k-1}\left\|p_j(x(\cdot))\right\|_{\mathbb{X}}\nonumber\\&\quad+\frac{M\alpha}{\Gamma(1+\alpha)}\left(\int_{0}^{\tau_k}(\tau_k-s)^{\frac{\alpha-1}{1-\alpha_1}}\mathrm{d}s\right)^{1-\alpha_1}\left\|\gamma_{r'}\right\|_{\mathcal{L}^{\frac{1}{\alpha_1}}([0,\tau_k];\mathbb{R^{+}})}+\frac{2T^{2\alpha-1}}{\lambda(2\alpha-1)}\left(\frac{M\tilde{M}\alpha}{\Gamma(1+\alpha)}\right)^{2}\left\|p_k(x(\cdot))\right\|_{\mathbb{X}}\nonumber\\&\quad+\frac{M\alpha}{\Gamma(1+\alpha)}\left(\int_{0}^{t}(t-s)^{\frac{\alpha-1}{1-\alpha_1}}\mathrm{d}s\right)^{1-\alpha_1}\left\|\gamma_{r'}\right\|_{\mathcal{L}^{\frac{1}{\alpha_1}}([0,t];\mathbb{R^{+}})}\nonumber\\&\le Ml_k+\tilde{R}\sum_{j=0}^{k}C_j+\frac{2T^{\alpha-\alpha_1}}{\mu^{1-\alpha_1}}\frac{M\alpha}{\Gamma(1+\alpha)}\left\|\gamma_{r'}\right\|_{\mathrm{L}^{\frac{1}{\alpha_1}}(J;\mathbb{R^+})}\nonumber\\&\le Ml_k+\tilde{R}\sum_{j=0}^{m}C_j+\frac{2T^{\alpha-\alpha_1}}{\mu^{1-\alpha_1}}\frac{M\alpha}{\Gamma(1+\alpha)}\left\|\gamma_{r'}\right\|_{\mathrm{L}^{\frac{1}{\alpha_1}}(J;\mathbb{R^+})}.
		\end{align}
		Using Assumption \ref{as2.1} (\textit{H2})(ii), we easily obtain
		\begin{align*}
		\liminf_{r \rightarrow \infty }\frac {\left\|\gamma_{r'}\right\|_{\mathrm{L}^{\frac{1}{\alpha_1}}(J;\mathbb{R^+})}}{r}&=\liminf_{r \rightarrow \infty } \left (\frac {\left\|\gamma_{r'}\right\|_{\mathrm{L}^{\frac{1}{\alpha_1}}(J;\mathbb{R^+})}}{r'}\times\frac{r'}{r}\right)=H_2\beta.
		\end{align*}
		Thus, dividing by $r$ in expressions \eqref{4.21}, \eqref{4.22}, \eqref{4.23} and then passing $r\to\infty$, we obtain
		\begin{align*}
		\frac{MH_{2}\alpha\beta}{\Gamma(1+\alpha)}\frac{2T^{\alpha-\alpha_1}}{\mu^{1-\alpha_1}}\left\{1+\frac{(m+1)(m+2)\tilde{R}}{2}+\frac{m(m+1)\tilde{R}^2}{2}\sum_{j=0}^{m-1}e^{\frac{(m+j)(m-j-1)\tilde{R}}{2}}\right\}>1,
		\end{align*}
		which is  a contradiction to \eqref{cnd}. Hence, for some  $ r>0$, $F_{\lambda}(\mathrm{E}_{r})\subset \mathrm{E}_{r}.$
		\vskip 0.1in 
		\noindent\textbf{Step (2): } \emph{The operator $ F_{\lambda}$ is continuous}. To achieve this goal, we consider a sequence $\{{x}^n\}^\infty_{n=1}\subseteq \mathrm{E}_r$ such that ${x}^n\rightarrow {x}\mbox{ in }{\mathrm{E}_r},$ that is,
		$$\lim\limits_{n\rightarrow \infty}\left\|x^n-x\right\|_{\mathrm{PC}(J;\mathbb{X})}=0.$$
		From Lemma \ref{lem2.7}, we infer that
		\begin{align*}
		\left\|\tilde{x_{s}^n}- \tilde{x_{s}}\right\|_{\mathfrak{B}}&\leq H_{2}\sup\limits_{\theta\in J}\left\|\tilde{x^{n}}(\theta)-\tilde{x}(\theta)\right\|_{\mathbb{X}}=H_{2}\left\|x^{n}-x\right\|_{\mathrm{PC}(J;\mathbb{X})}\rightarrow 0 \ \mbox{ as } \ n\rightarrow\infty,
		\end{align*}
		for all $s\in\mathcal{Q}(\rho^-)\cup J$. Since $\rho(s, \tilde{x_s^k})\in \mathcal{Q}(\rho^-)\cup J,$ for all $k\in\mathbb{N}$, then we conclude that
		\begin{align*}
		\left\|\tilde{x^n}_{\rho(s,\tilde{x_s^k})}-\tilde{x}_{\rho(s,\tilde{x_s^k})}\right\|_{\mathfrak{B}}\rightarrow 0 \ \mbox{ as } \ n\rightarrow\infty, \ \mbox{ for all }\ s\in J\  \mbox{ and }\ k\in\mathbb{N}.
		\end{align*}
		In particular, we choose $k=n$ and use the above convergence together with Assumption \ref{as2.1} \textit{(H1)} to  obtain
		\begin{align}\label{4.25}
		\left\|f(s, \tilde{x^n}_{\rho(s,\tilde{x_s^n})})-f(s,\tilde{x}_{\rho(s,\tilde{x_s})})\right\|_{\mathbb{X}}&\leq\left\|f(s, \tilde{x^n}_{\rho(s,\tilde{x_s^n})})-f(s,\tilde{x}_{\rho(s,\tilde{x_s^n})})\right\|_{\mathbb{X}}\nonumber\\&\quad+\left\|f(s, \tilde{x}_{\rho(s,\tilde{x_s^n})})-f(s,\tilde{x}_{\rho(s,\tilde{x_s})})\right\|_{\mathbb{X}}\nonumber\\&\to 0\ \mbox{ as }\ n\to\infty, \mbox{ uniformly for } \ s\in J. 
		\end{align}
		From the above convergence and the dominated convergence theorem, we evaluate
		\begin{align}\label{4.26}
		\left\|p_0(x^n(\cdot))-p_0(x(\cdot))\right\|_{\mathbb{X}}&\le\left\|\int^{t_1}_{0}(t_1-s)^{\alpha-1}\widehat{\mathcal{T}}_{\alpha}(t_1-s)\left[f(s,\tilde{x^n}_{\rho(s,\tilde{x^n_s})})-f(s,\tilde{x}_{\rho(s,\tilde{x_s})})\right]\mathrm{d}s\right\|_{\mathbb{X}}\nonumber\\&\le\frac{M\alpha}{\Gamma(1+\alpha)}\int^{t_1}_{0}(t_1-s)^{\alpha-1}\left\|f(s,\tilde{x^n}_{\rho(s,\tilde{x^n_s})})-f(s,\tilde{x}_{\rho(s,\tilde{x_s})})\right\|_{\mathbb{X}}\mathrm{d}s\nonumber\\&\to 0\ \mbox{ as }\ n\to\infty.
		\end{align}
		Using the convergences \eqref{4.26} and the relation \eqref{2.5}, we calculate
		\begin{align*}
		\left\|\mathcal{R}(\lambda,\Phi_{0}^{t_{1}})p_0(x^{n}(\cdot))-\mathcal{R}(\lambda,\Phi_{0}^{t_{1}})p_0(x(\cdot))\right\|_{\mathbb{X}}&=\frac{1}{\lambda}\left\|\lambda\mathcal{R}(\lambda,\Phi_{\tau_k}^{t_{k+1}})\left(p_0(x^{n}(\cdot))-p_0(x(\cdot))\right)\right\|_{\mathbb{X}}\nonumber\\&\leq\frac{1}{\lambda}\left\|p_0(x^{n}(\cdot))-p_0(x(\cdot))\right\|_{\mathbb{X}}\nonumber\\&\to 0 \ \mbox{ as } \ n \to \infty.
		\end{align*}
		Since the mapping $\mathcal{J}:\mathbb{X}\to\mathbb{X}^{*}$  is  demicontinuous, it is easy to obtain
		\begin{align}\label{4.2}
		\mathcal{J}\left[\mathcal{R}(\lambda,\Phi_{\tau_k}^{t_{k+1}})p_0(x^{n}(\cdot))\right]\xrightharpoonup{w}\mathcal{J}\left[\mathcal{R}(\lambda,\Phi_{\tau_k}^{t_{k+1}})p_0(x(\cdot))\right] \ \mbox{ as } \ n\to\infty  \ \mbox{ in }\ \mathbb{X}^{*}.
		\end{align}
		From Lemma \ref{lem2.5}, we infer that the operator $\widehat{\mathcal{T}}_{\alpha}(t)$ is compact for $t>0$. Therefore, the operator $\widehat{\mathcal{T}}_{\alpha}(t)^*$ is also compact for $t>0$.   Hence, by using the compactness of this operator together with the weak convergence \eqref{4.2}, one can obtain
		\begin{align}\label{4.1}
		&\left\|u^{n,\alpha}_{0,\lambda}(t)-u_{0,\lambda}^{\alpha}(t)\right||_{\mathbb{U}}\nonumber\\&\le\left\|(t_{1}-t)^{\alpha-1}\mathrm{B}^*\widehat{\mathcal{T}}_{\alpha}(t_{1}-t)^*\left[\mathcal{J}\left[\mathcal{R}(\lambda,\Phi_{0}^{t_{1}})p_0(x^{n}(\cdot))\right]-\mathcal{J}\left[\mathcal{R}(\lambda,\Phi_{0}^{t_{1}})p_0k(x(\cdot))\right]\right]\right\|_{\mathbb{U}}\nonumber\\&\le(t_{1}-t)^{\alpha-1}\tilde{M}\left\|\widehat{\mathcal{T}}_{\alpha}(t_{1}-t)^*\left[\mathcal{J}\left[\mathcal{R}(\lambda,\Phi_{0}^{t_{1}})p_0(x^{n}(\cdot))\right]-\mathcal{J}\left[\mathcal{R}(\lambda,\Phi_{0}^{t_{1}})p_0(x(\cdot))\right]\right]\right\|_{\mathbb{X}}\nonumber\\&\to 0\ \mbox{ as }\ n\to\infty, \ \mbox{ for  each }\ t\in[0,t_{1}).
		\end{align}
		Similarly, for $k=1$, we compute
		\begin{align}\label{4.27}
		\left\|p_1(x^n(\cdot))-p_1(x(\cdot))\right\|_{\mathbb{X}}&\le\left\|\mathcal{T}_{\alpha}(t_{2}-\tau_1)\left[h_1(\tau_1,\tilde{x^n}(t_1^-))-h_1(\tau_1,\tilde{x}(t_1^-))\right]\right\|_{\mathbb{X}}\nonumber\\&\quad+\left\|\int_{0}^{\tau_{1}}(\tau_{1}-s)^{\alpha-1}\widehat{\mathcal{T}}_{\alpha}(\tau_{1}-s)\left[f(s,\tilde{x^n}_{\rho(s,\tilde{x^n_s})})-f(s,\tilde{x}_{\rho(s,\tilde{x}_s)})\right]\mathrm{d}s\right\|_{\mathbb{X}}\nonumber\\&\quad+\left\|\int_{0}^{t_{2}}(t_{2}-s)^{\alpha-1}\widehat{\mathcal{T}}_{\alpha}(t_{2}-s)\left[f(s,\tilde{x^n}_{\rho(s,\tilde{x^n_s})})-f(s,\tilde{x}_{\rho(s,\tilde{x}_s)})\right]\mathrm{d}s\right\|_{\mathbb{X}}\nonumber\\&\quad+\left\|\int_{0}^{\tau_1}(\tau_1-s)^{\alpha-1}\widehat{\mathcal{T}}_{\alpha}(\tau_1-s)\mathrm{B}\left[u^{n,\alpha}_{0,\lambda}(s)-u^{\alpha}_{0,\lambda}(s)\right]\mathrm{d}s\right\|_{\mathbb{X}} \nonumber\\&\quad+\left\|\int_{0}^{\tau_1}(t_2-s)^{\alpha-1}\widehat{\mathcal{T}}_{\alpha}(t_2-s)\mathrm{B}\left[u^{n,\alpha}_{0,\lambda}(s)-u^{\alpha}_{0,\lambda}(s)\right]\mathrm{d}s\right\|_{\mathbb{X}} \nonumber\\&\le M\left\|h_1(\tau_1,\tilde{x^n}(t_1^-))-h_1(\tau_1,\tilde{x}(t_1^-))\right\|_{\mathbb{X}}\nonumber\\&\quad+\frac{M\alpha}{\Gamma(1+\alpha)}\int_{0}^{\tau_{1}}(\tau_{1}-s)^{\alpha-1}\left\|f(s,\tilde{x^n}_{\rho(s,\tilde{x^n_s})})-f(s,\tilde{x}_{\rho(s,\tilde{x}_s)})\right\|_{\mathbb{X}}\mathrm{d}s\nonumber\\&\quad+\frac{M\alpha}{\Gamma(1+\alpha)}\int_{0}^{t_{2}}(t_{2}-s)^{\alpha-1}\left\|f(s,\tilde{x^n}_{\rho(s,\tilde{x^n_s})})-f(s,\tilde{x}_{\rho(s,\tilde{x}_s)})\right\|_{\mathbb{X}}\mathrm{d}s\nonumber\\&\quad+\frac{M\tilde{M}\alpha}{\Gamma(1+\alpha)}\int_{0}^{t_1}(\tau_1-s)^{\alpha-1}\left\|u^{n,\alpha}_{0,\lambda}(s)-u^{\alpha}_{0,\lambda}(s)\right\|_{\mathbb{U}}\mathrm{d}s\nonumber\\&\quad+\frac{M\tilde{M}\alpha}{\Gamma(1+\alpha)}\int_{0}^{t_1}(t_2-s)^{\alpha-1}\left\|u^{n,\alpha}_{0,\lambda}(s)-u^{\alpha}_{0,\lambda}(s)\right\|_{\mathbb{U}}\mathrm{d}s\nonumber\\&\to0\ \mbox{ as }\ n\to\infty, 
		\end{align}
		where we used Assumption \ref{as2.1} (\textit{H3}),  convergences \eqref{4.25}, \eqref{4.1} and the dominated convergence theorem. Moreover, similar to the convergence \eqref{4.1}, one can obtain
		\begin{align*}
		&\left\|u^{n,\alpha}_{1,\lambda}(t)-u_{1,\lambda}^{\alpha}(t)\right||_{\mathbb{U}}\to 0\ \mbox{ as }\ n\to\infty, \ \mbox{ for  each }\ t\in[\tau_1,t_{2}).
		\end{align*}
		Further, applying a similar analogy as above  for $k=2,\ldots,m,$ one can compute 
		\begin{align*}
		&\left\|u^{n,\alpha}_{k,\lambda}(t)-u_{k,\lambda}^{\alpha}(t)\right||_{\mathbb{U}}\to 0\ \mbox{ as }\ n\to\infty, \mbox{ for  each }\ t\in[\tau_k,t_{k+1}), k=2,\ldots,m.
		\end{align*} 
		Therefore, we have
		\begin{align}\label{4.29}
		\left\|u^{n,\alpha}_{\lambda}(t)-u_{\lambda}^{\alpha}(t)\right||_{\mathbb{U}}\to 0\ \mbox{ as }\ n\to\infty,\ \mbox{ for  each }\ t\in[\tau_k,t_{k+1}), k=0,1,\ldots,m.
		\end{align}
		Using the convergences \eqref{4.25}, \eqref{4.29} and the dominated convergence theorem, we arrive at 
		\begin{align*}
		\left\|(F_{\lambda}x^n)(t)-(F_{\lambda}x)(t)\right\|_{\mathbb{X}}&\le\int_{0}^{t}(t-s)^{\alpha-1}\left\|\widehat{\mathcal{T}}_{\alpha}(t-s)\mathrm{B}\left[u^{n,\alpha}_{\lambda}(s)-u_{\lambda}^{\alpha}(s)\right]\right\|_{\mathbb{X}}\mathrm{d}s\nonumber\\&\quad+\int_{0}^{t}(t-s)^{\alpha-1}\left\|\widehat{\mathcal{T}}_{\alpha}(t-s)\left[f(s,\tilde{x^n}_{\rho(s,\tilde{x^n_s})})-f(s,\tilde{x}_{\rho(s,\tilde{x}_s)})\right]\right\|_{\mathbb{X}}\mathrm{d}s\nonumber\\&\le\frac{M\tilde{M}\alpha}{\Gamma(1+\alpha)}\int_{0}^{t}(t-s)^{\alpha-1}\left\|u^{n,\alpha}_{\lambda}(s)-u_{\lambda}^{\alpha}(s)\right\|_{\mathbb{U}}\mathrm{d}s\nonumber\\&\quad+\frac{M\alpha}{\Gamma(1+\alpha)}\int_{0}^{t}(t-s)^{\alpha-1}\left\|f(s,\tilde{x^n}_{\rho(s, \tilde{x^n_s})})-f(s,\tilde{x}_{\rho(s, \tilde{x}_s)})\right\|_{\mathbb{X}}\mathrm{d}s\nonumber\\&\to 0 \ \mbox{ as }\ n\to\infty, \ \mbox{ for }\ t\in[0,t_1].
		\end{align*}
		Similarly, for  $t\in(\tau_k,t_{k+1}],\ k=1,\ldots,m$, we deduce that 
		\begin{align*}
		\left\|(F_{\lambda}x^n)(t)-(F_{\lambda}x)(t)\right\|_{\mathbb{X}}&\le\left\|\mathcal{T}_{\alpha}(t-\tau_k)\left[h_k(\tau_k,\tilde{x^n}(t_k^-))-h_k(\tau_k,\tilde{x}(t_k^-))\right]\right\|_{\mathbb{X}}\nonumber\\&\quad+\int_{0}^{\tau_k}(\tau_k-s)^{\alpha-1}\left\|\widehat{\mathcal{T}}_{\alpha}(\tau_k-s)\left[f(s,\tilde{x^n}_{\rho(s, \tilde{x^n_s})})-f(s,\tilde{x}_{\rho(s, \tilde{x}_s)})\right]\right\|_{\mathbb{X}}\mathrm{d}s\nonumber\\&\quad+\int_{0}^{\tau_k}(\tau_1-s)^{\alpha-1}\left\|\widehat{\mathcal{T}}_{\alpha}(\tau_1-s)\mathrm{B}\left[u^{n,\alpha}_{\lambda}(s)-u_{\lambda}^{\alpha}(s)\right]\right\|_{\mathbb{X}}\mathrm{d}s\nonumber\\&\quad+\int_{0}^{t}(t-s)^{\alpha-1}\left\|\widehat{\mathcal{T}}_{\alpha}(t-s)\mathrm{B}\left[u^{n,\alpha}_{\lambda}(s)-u_{\lambda}^{\alpha}(s)\right]\right\|_{\mathbb{X}}\mathrm{d}s\nonumber\\&\quad+\int_{0}^{t}(t-s)^{\alpha-1}\left\|\widehat{\mathcal{T}}_{\alpha}(t-s)\left[f(s,\tilde{x^n}_{\rho(s, \tilde{x^n_s})})-f(s,\tilde{x}_{\rho(s, \tilde{x}_s)})\right]\right\|_{\mathbb{X}}\mathrm{d}s\nonumber\\&\le M\left\|h_k(\tau_k,\tilde{x^n}(t_k^-))-h_k(\tau_k,\tilde{x}(t_k^-))\right\|_{\mathbb{X}}\nonumber\\&\quad+\frac{M\alpha}{\Gamma(1+\alpha)}\int_{0}^{\tau_k}(\tau_k-s)^{\alpha-1}\left\|f(s,\tilde{x^n}_{\rho(s, \tilde{x^n_s})})-f(s,\tilde{x}_{\rho(s, \tilde{x}_s)})\right\|_{\mathbb{X}}\mathrm{d}s\nonumber\\&\quad+ \frac{M\tilde{M}\alpha}{\Gamma(1+\alpha)}\int_{0}^{\tau_k}(\tau_k-s)^{\alpha-1}\left\|u^{n,\alpha}_{\lambda}(s)-u_{\lambda}^{\alpha}(s)\right\|_{\mathbb{U}}\mathrm{d}s\nonumber\\&\quad+\frac{M\tilde{M}\alpha}{\Gamma(1+\alpha)}\int_{0}^{t}(t-s)^{\alpha-1}\left\|u^{n,\alpha}_{\lambda}(s)-u_{\lambda}^{\alpha}(s)\right\|_{\mathbb{U}}\mathrm{d}s\nonumber\\&\quad+\frac{M\alpha}{\Gamma(1+\alpha)}\int_{0}^{t}(t-s)^{\alpha-1}\left\|f(s,\tilde{x^n}_{\rho(s, \tilde{x^n_s})})-f(s,\tilde{x}_{\rho(s, \tilde{x}_s)})\right\|_{\mathbb{X}}\mathrm{d}s\nonumber\\&\to 0 \ \mbox{ as }\ n\to\infty.
		\end{align*}
		Moreover, for $t\in(t_k, \tau_k],\ k=1,\ldots,m$, applying Assumption \ref{as2.1} (\textit{H3}), we obtain
		\begin{align*}
		\left\|(F_{\lambda}x^n)(t)-(F_{\lambda}x)(t)\right\|_{\mathbb{X}}&\le\left\|h_k(t,\tilde{x^n}(t_k^-))-h_k(t,\tilde{x}(t_k^-))\right\|_{\mathbb{X}}\to 0 \ \mbox{ as }\ n\to\infty.
		\end{align*}
		Hence, it follows that $F_{\lambda}$ is continuous.
		\vskip 0.1in 
		\noindent\textbf{Step (3): } \emph{$ F_{\lambda}$ is a compact operator.} In order to prove this claim, we use the  well-known Ascoli-Arzela theorem. According to the infinite-dimensional version of the Ascoli-Arzela theorem (see, Theorem 3.7, Chapter 2, \cite{JYONG}), it is enough to show that 
		\begin{itemize}
			\item [(i)] the image of $\mathrm{E}_r$ under $F_{\lambda}$ is uniformly bounded (which is proved in Step I),
			\item [(ii)] the image of $\mathrm{E}_r$ under $F_{\lambda}$ is equicontinuous,
			\item [(iii)] for an arbitrary $t\in J$, the set $\mathrm{V}(t)=\{(F_\lambda x)(t):x\in \mathrm{E}_r\}$ is relatively compact.
		\end{itemize}

		First, we claim that the image of $\mathrm{E}_r$ under $F_{\lambda}$ is equicontinuous. For $s_1,s_2\in[0,t_1]$ such that $s_1<s_2$ and $x\in\mathrm{E}_r$, we compute
		\begin{align}
		&\left\|(F_{\lambda}x)(s_2)-(F_{\lambda}x)(s_1)\right\|_{\mathbb{X}}\nonumber\\&\le\left\|\left[\mathcal{T}_{\alpha}(s_2)-\mathcal{T}_{\alpha}(s_1)\right]\psi(0)\right\|_{\mathbb{X}}+\left\|\int_{s_1}^{s_2}(s_2-s)^{\alpha-1}\widehat{\mathcal{T}}_{\alpha}(s_2-s)f(s,\tilde{x}_{\rho(s,\tilde{x}_s)})\mathrm{d}s\right\|_{\mathbb{X}}\nonumber\\&\quad+\left\|\int_{s_1}^{s_2}(s_2-s)^{\alpha-1}\widehat{\mathcal{T}}_{\alpha}(s_2-s)\mathrm{B}u^{\alpha}_{\lambda}(s)\mathrm{d}s\right\|_{\mathbb{X}}\nonumber\\&\quad+\left\|\int_{0}^{s_1}\left[(s_2-s)^{\alpha-1}-(s_1-s)^{\alpha-1}\right]\widehat{\mathcal{T}}_{\alpha}(s_2-s)\mathrm{B}u^{\alpha}_{\lambda}(s)\mathrm{d}s\right\|_{\mathbb{X}}\nonumber\\&\quad+\left\|\int_{0}^{s_1}(s_1-s)^{\alpha-1}\left[\widehat{\mathcal{T}}_{\alpha}(s_2-s)-\widehat{\mathcal{T}}_{\alpha}(s_1-s)\right]\mathrm{B}u^{\alpha}_{\lambda}(s)\mathrm{d}s\right\|_{\mathbb{X}}\nonumber\\&\quad+\left\|\int_{0}^{s_1}\left[(s_2-s)^{\alpha-1}-(s_1-s)^{\alpha-1}\right]\widehat{\mathcal{T}}_{\alpha}(s_2-s)f(s,\tilde{x}_{\rho(s,\tilde{x}_s)})\mathrm{d}s\right\|_{\mathbb{X}}\nonumber\\&\quad+\left\|\int_{0}^{s_1}(s_1-s)^{\alpha-1}\left[\widehat{\mathcal{T}}_{\alpha}(s_2-s)-\widehat{\mathcal{T}}_{\alpha}(s_1-s)\right]f(s,\tilde{x}_{\rho(s,\tilde{x}_s)})\mathrm{d}s\right\|_{\mathbb{X}}\nonumber\\&\le\left\|\mathcal{T}_{\alpha}(s_2)-\mathcal{T}_{\alpha}(s_1)\right\|_{\mathcal{L}(\mathbb{X})}\left\|\psi(0)\right\|_{\mathbb{X}}+\frac{M\alpha}{\Gamma(1+\alpha)}\int_{s_1}^{s_2}(s_2-s)^{\alpha-1}\gamma_{r'}(s)\mathrm{d}s\nonumber\\&\quad+\frac{N_0}{\lambda}\left(\frac{M\tilde{M}\alpha}{\Gamma(1+\alpha)}\right)^{2}\int_{s_1}^{s_2}(s_2-s)^{\alpha-1}(t_1-s)^{\alpha-1}\mathrm{d}s\nonumber\\&\quad+\frac{N_0}{\lambda}\left(\frac{M\tilde{M}\alpha}{\Gamma(1+\alpha)}\right)^{2}\int_{0}^{s_1}\left|(s_2-s)^{\alpha-1}-(s_1-s)^{\alpha-1}\right|(t_1-s)^{\alpha-1}\mathrm{d}s\nonumber\\&\quad+\frac{M\tilde{M}^{2}N_{0}\alpha}{\lambda\Gamma(1+\alpha)}\int_{0}^{s_1}(s_1-s)^{\alpha-1}\left\|\widehat{\mathcal{T}}_{\alpha}(s_2-s)-\widehat{\mathcal{T}}_{\alpha}(s_1-s)\right\|_{\mathcal{L}(\mathbb{X})}(t_1-s)^{\alpha-1}\mathrm{d}s\nonumber\\&\quad+\frac{M\alpha}{\Gamma(1+\alpha)}\int_{0}^{s_1}\left|(s_2-s)^{\alpha-1}-(s_1-s)^{\alpha-1}\right|\gamma_{r'}(s)\mathrm{d}s\nonumber\\&\quad+\int_{0}^{s_1}(s_1-s)^{\alpha-1}\left\|\widehat{\mathcal{T}}_{\alpha}(s_2-s)-\widehat{\mathcal{T}}_{\alpha}(s_1-s)\right\|_{\mathcal{L}(\mathbb{X})}\gamma_{r'}(s)\mathrm{d}s\nonumber\\&\le\left\|\mathcal{T}_{\alpha}(s_2)-\mathcal{T}_{\alpha}(s_1)\right\|_{\mathcal{L}(\mathbb{X})}\left\|\psi(0)\right\|_{\mathbb{X}}+\frac{M\alpha}{\Gamma(1+\alpha)}\frac{(s_2-s_1)^{\alpha-\alpha_1}}{\mu^{1-\alpha_1}}\left(\int_{s_1}^{s_2}(\gamma_{r'}(s))^{\frac{1}{\alpha_1}}\mathrm{d}s\right)^{\alpha_1}\nonumber\\&\quad+\frac{N_{0}}{\lambda}\left(\frac{M\tilde{M}\alpha}{\Gamma(1+\alpha)}\right)^{2}\frac{(s_2-s_1)^{2\alpha-1}}{2\alpha-1}\nonumber\\&\quad+\frac{N_{0}}{\lambda}\left(\frac{M\tilde{M}\alpha}{\Gamma(1+\alpha)}\right)^{2}\int_{0}^{s_1}\left|(s_2-s)^{\alpha-1}-(s_1-s)^{\alpha-1}\right|(t_1-s)^{\alpha-1}\mathrm{d}s\nonumber\\&\quad+\frac{M\tilde{M}^{2}N_{0}\alpha}{\lambda\Gamma(1+\alpha)}\int_{0}^{s_1}(s_1-s)^{\alpha-1}\left\|\widehat{\mathcal{T}}_{\alpha}(s_2-s)-\widehat{\mathcal{T}}_{\alpha}(s_1-s)\right\|_{\mathcal{L}(\mathbb{X})}(t_1-s)^{\alpha-1}\mathrm{d}s\nonumber\\&\quad+\frac{M\alpha}{\Gamma(1+\alpha)}\int_{0}^{s_1}\left|(s_2-s)^{\alpha-1}-(s_1-s)^{\alpha-1}\right|\gamma_{r'}(s)\mathrm{d}s\nonumber\\&\quad+\int_{0}^{s_1}(s_1-s)^{\alpha-1}\left\|\widehat{\mathcal{T}}_{\alpha}(s_2-s)-\widehat{\mathcal{T}}_{\alpha}(s_1-s)\right\|_{\mathcal{L}(\mathbb{X})}\gamma_{r'}(s)\mathrm{d}s.
		\end{align}
		If $s_1=0,$ then from the above estimate, we deduce that 
		\begin{align*}
		\lim_{s_2\to0^+}\left\|(F_{\lambda}x)(s_2)-(F_{\lambda}x)(s_1)\right\|_{\mathbb{X}}=0\ \mbox{ unifromly for } \ x\in\mathrm{E}_r.
		\end{align*} 
		For $0<\nu<s_1<t_1$, we have
		\begin{align}\label{4.31}
		&\left\|(F_{\lambda}x)(s_2)-(F_{\lambda}x)(s_1)\right\|_{\mathbb{X}}\nonumber\\&\le\left\|\mathcal{T}_{\alpha}(s_2)-\mathcal{T}_{\alpha}(s_1)\right\|_{\mathcal{L}(\mathbb{X})}\left\|\psi(0)\right\|_{\mathbb{X}}+\frac{M\alpha}{\Gamma(1+\alpha)}\frac{(s_2-s_1)^{\alpha-\alpha_1}}{\mu^{1-\alpha_1}}\left(\int_{s_1}^{s_2}(\gamma_{r'}(s))^{\frac{1}{\alpha_1}}\mathrm{d}s\right)^{\alpha_1}\nonumber\\&\quad+\frac{N_{0}}{\lambda}\left(\frac{M\tilde{M}\alpha}{\Gamma(1+\alpha)}\right)^{2}\frac{(s_2-s_1)^{2\alpha-1}}{2\alpha-1}\nonumber\\&\quad+\frac{N_{0}}{\lambda}\left(\frac{M\tilde{M}\alpha}{\Gamma(1+\alpha)}\right)^{2}\int_{0}^{s_1}\left|(s_2-s)^{\alpha-1}-(s_1-s)^{\alpha-1}\right|(t_1-s)^{\alpha-1}\mathrm{d}s\nonumber\\&\quad+\frac{M\tilde{M}^{2}N_{0}\alpha}{\lambda\Gamma(1+\alpha)}\int_{0}^{s_1-\nu}(s_1-s)^{\alpha-1}\left\|\widehat{\mathcal{T}}_{\alpha}(s_2-s)-\widehat{\mathcal{T}}_{\alpha}(s_1-s)\right\|_{\mathcal{L}(\mathbb{X})}(t_1-s)^{\alpha-1}\mathrm{d}s\nonumber\\&\quad+\frac{M\tilde{M}^{2}N_{0}\alpha}{\lambda\Gamma(1+\alpha)}\int_{s_1-\nu}^{s_1}(s_1-s)^{\alpha-1}\left\|\widehat{\mathcal{T}}_{\alpha}(s_2-s)-\widehat{\mathcal{T}}_{\alpha}(s_1-s)\right\|_{\mathcal{L}(\mathbb{X})}(t_1-s)^{\alpha-1}\mathrm{d}s\nonumber\\&\quad+\frac{M\alpha}{\Gamma(1+\alpha)}\int_{0}^{s_1}\left|(s_2-s)^{\alpha-1}-(s_1-s)^{\alpha-1}\right|\gamma_{r'}(s)\mathrm{d}s\nonumber\\&\quad+\int_{0}^{s_1-\nu}(s_1-s)^{\alpha-1}\left\|\widehat{\mathcal{T}}_{\alpha}(s_2-s)-\widehat{\mathcal{T}}_{\alpha}(s_1-s)\right\|_{\mathcal{L}(\mathbb{X})}\gamma_{r'}(s)\mathrm{d}s\nonumber\\&\quad+\int_{s_1-\nu}^{s_1}(s_1-s)^{\alpha-1}\left\|\widehat{\mathcal{T}}_{\alpha}(s_2-s)-\widehat{\mathcal{T}}_{\alpha}(s_1-s)\right\|_{\mathcal{L}(\mathbb{X})}\gamma_{r'}(s)\mathrm{d}s\nonumber\\&\le\left\|\mathcal{T}_{\alpha}(s_2)-\mathcal{T}_{\alpha}(s_1)\right\|_{\mathcal{L}(\mathbb{X})}\left\|\psi(0)\right\|_{\mathbb{X}}+\frac{M\alpha}{\Gamma(1+\alpha)}\frac{(s_2-s_1)^{\alpha-\alpha_1}}{\mu^{1-\alpha_1}}\left(\int_{s_1}^{s_2}(\gamma_{r'}(s))^{\frac{1}{\alpha_1}}\mathrm{d}s\right)^{\alpha_1}\nonumber\\&\quad+\frac{N_{0}}{\lambda}\left(\frac{M\tilde{M}\alpha}{\Gamma(1+\alpha)}\right)^{2}\frac{(s_2-s_1)^{2\alpha-1}}{2\alpha-1}\nonumber\\&\quad+\frac{N_{0}}{\lambda}\left(\frac{M\tilde{M}\alpha}{\Gamma(1+\alpha)}\right)^{2}\int_{0}^{s_1}\left|(s_2-s)^{\alpha-1}-(s_1-s)^{\alpha-1}\right|(t_1-s)^{\alpha-1}\mathrm{d}s\nonumber\\&\quad+\frac{M\tilde{M}^{2}N_{0}\alpha}{\lambda\Gamma(1+\alpha)}\sup_{t\in[0,s_1-\nu]}\left\|\widehat{\mathcal{T}}_{\alpha}(s_2-t)-\widehat{\mathcal{T}}_{\alpha}(s_1-t)\right\|_{\mathcal{L}(\mathbb{X})}\int_{0}^{s_1-\nu}(s_1-s)^{\alpha-1}(t_1-s)^{\alpha-1}\mathrm{d}s\nonumber\\&\quad+\frac{2}{\lambda}\left(\frac{M\tilde{M}^{2}\alpha}{\Gamma(1+\alpha)}\right)^{2}\frac{\nu^{2\alpha-1}}{2\alpha-1}+\frac{M\alpha}{\Gamma(1+\alpha)}\int_{0}^{s_1}\left|(s_2-s)^{\alpha-1}-(s_1-s)^{\alpha-1}\right|\gamma_{r'}(s)\mathrm{d}s\nonumber\\&\quad+\sup_{t\in[0,s_1-\nu]}\left\|\widehat{\mathcal{T}}_{\alpha}(s_2-t)-\widehat{\mathcal{T}}_{\alpha}(s_1-t)\right\|_{\mathcal{L}(\mathbb{X})}\int_{0}^{s_1-\nu}(s_1-s)^{\alpha-1}\gamma_{r'}(s)\mathrm{d}s\nonumber\\&\quad+\frac{2M\alpha}{\Gamma(1+\alpha)}\frac{\nu^{\alpha-\alpha_1}}{\mu^{1-\alpha_1}}\left(\int_{s_1-\nu}^{s_1}(\gamma_{r'}(s))^{\frac{1}{\alpha_1}}\mathrm{d}s\right)^{\alpha_1}.
		\end{align}
		Similarly for $s_1,s_2\in(\tau_k,t_{k+1}], \ k=1,\ldots,m$ with $s_1<s_2$ and $x\in\mathrm{E}_r$, one can estimate 
		\begin{align}\label{4.32}
		&\left\|(F_{\lambda}x)(s_2)-(F_{\lambda}x)(s_1)\right\|_{\mathbb{X}}\nonumber\\&\le\left\|\mathcal{T}_{\alpha}(s_2-\tau_k)-\mathcal{T}_{\alpha}(s_1-\tau_k)\right\|_{\mathcal{L}(\mathbb{X})}l_k+\frac{M\alpha}{\Gamma(1+\alpha)}\frac{(s_2-s_1)^{\alpha-\alpha_1}}{\mu^{1-\alpha_1}}\left(\int_{s_1}^{s_2}(\gamma_{r'}(s))^{\frac{1}{\alpha_1}}\mathrm{d}s\right)^{\alpha_1}\nonumber\\&\quad+\frac{C_k}{\lambda}\left(\frac{M\tilde{M}\alpha}{\Gamma(1+\alpha)}\right)^{2}\frac{(s_2-s_1)^{2\alpha-1}}{2\alpha-1}+\frac{M\tilde{M}\alpha}{\Gamma(1+\alpha)}\int_{0}^{s_1}\left|(s_2-s)^{\alpha-1}-(s_1-s)^{\alpha-1}\right|\left\|u_{\lambda}^{\alpha}(s)\right\|_{\mathbb{U}}\mathrm{d}s\nonumber\\&\quad+\tilde{M}\sup_{t\in[0,s_1-\nu]}\left\|\widehat{\mathcal{T}}_{\alpha}(s_2-t)-\widehat{\mathcal{T}}_{\alpha}(s_1-t)\right\|_{\mathcal{L}(\mathbb{X})}\int_{0}^{s_1-\nu}(s_1-s)^{\alpha-1}\left\|u_{\lambda}^{\alpha}(s)\right\|_{\mathbb{U}}\mathrm{d}s\nonumber\\&\quad+\frac{2C_k}{\lambda}\left(\frac{M\tilde{M}^{2}\alpha}{\Gamma(1+\alpha)}\right)^{2}\frac{\nu^{2\alpha-1}}{2\alpha-1}+\frac{M\alpha}{\Gamma(1+\alpha)}\int_{0}^{s_1}\left|(s_2-s)^{\alpha-1}-(s_1-s)^{\alpha-1}\right|\gamma_{r'}(s)\mathrm{d}s\nonumber\\&\quad+\sup_{t\in[0,s_1-\nu]}\left\|\widehat{\mathcal{T}}_{\alpha}(s_2-t)-\widehat{\mathcal{T}}_{\alpha}(s_1-t)\right\|_{\mathcal{L}(\mathbb{X})}\int_{0}^{s_1-\nu}(s_1-s)^{\alpha-1}\gamma_{r'}(s)\mathrm{d}s\nonumber\\&\quad+\frac{2M\alpha}{\Gamma(1+\alpha)}\frac{\nu^{\alpha-\alpha_1}}{\mu^{1-\alpha_1}}\left(\int_{s_1-\nu}^{s_1}(\gamma_{r'}(s))^{\frac{1}{\alpha_1}}\mathrm{d}s\right)^{\alpha_1}\nonumber\\&\le\left\|\mathcal{T}_{\alpha}(s_2-\tau_k)-\mathcal{T}_{\alpha}(s_1-\tau_k)\right\|_{\mathcal{L}(\mathbb{X})}l_k+\frac{M\alpha}{\Gamma(1+\alpha)}\frac{(s_2-s_1)^{\alpha-\alpha_1}}{\mu^{1-\alpha_1}}\left(\int_{s_1}^{s_2}(\gamma_{r'}(s))^{\frac{1}{\alpha_1}}\mathrm{d}s\right)^{\alpha_1}\nonumber\\&\quad+\frac{C_k}{\lambda}\left(\frac{M\tilde{M}\alpha}{\Gamma(1+\alpha)}\right)^{2}\frac{(s_2-s_1)^{2\alpha-1}}{2\alpha-1}\nonumber\\&\quad+\sum_{j=0}^{k-1}\frac{C_j}{\lambda}\left(\frac{M\tilde{M}\alpha}{\Gamma(1+\alpha)}\right)^{2}\int_{\tau_j}^{t_{j+1}}\left|(s_2-s)^{\alpha-1}-(s_1-s)^{\alpha-1}\right|(t_{j+1}-s)^{\alpha-1}\mathrm{d}s\nonumber\\&\quad+\frac{C_k}{\lambda}\left(\frac{M\tilde{M}\alpha}{\Gamma(1+\alpha)}\right)^{2}\int_{\tau_k}^{s_1}\left|(s_2-s)^{\alpha-1}-(s_1-s)^{\alpha-1}\right|(t_{k+1}-s)^{\alpha-1}\mathrm{d}s\nonumber\\&\quad+\sup_{t\in[0,s_1-\nu]}\left\|\widehat{\mathcal{T}}_{\alpha}(s_2-t)-\widehat{\mathcal{T}}_{\alpha}(s_1-t)\right\|_{\mathcal{L}(\mathbb{X})}\sum_{j=0}^{k-1}\frac{C_j}{\lambda}\frac{M\tilde{M}^{2}\alpha}{\Gamma(1+\alpha)}\int_{\tau_j}^{t_{j+1}}(s_1-s)^{\alpha-1}(t_{j+1}-s)^{\alpha-1}\mathrm{d}s\nonumber\\&\quad+\sup_{t\in[0,s_1-\nu]}\left\|\widehat{\mathcal{T}}_{\alpha}(s_2-t)-\widehat{\mathcal{T}}_{\alpha}(s_1-t)\right\|_{\mathcal{L}(\mathbb{X})}\frac{C_k}{\lambda}\frac{M\tilde{M}^{2}\alpha}{\Gamma(1+\alpha)}\int_{\tau_k}^{s_1-\nu}(s_1-s)^{\alpha-1}(t_{k+1}-s)^{\alpha-1}\mathrm{d}s\nonumber\\&\quad+\frac{2C_k}{\lambda}\left(\frac{M\tilde{M}^{2}\alpha}{\Gamma(1+\alpha)}\right)^{2}\frac{\nu^{2\alpha-1}}{2\alpha-1}+\frac{M\alpha}{\Gamma(1+\alpha)}\int_{0}^{s_1}\left|(s_2-s)^{\alpha-1}-(s_1-s)^{\alpha-1}\right|\gamma_{r'}(s)\mathrm{d}s\nonumber\\&\quad+\sup_{t\in[0,s_1-\nu]}\left\|\widehat{\mathcal{T}}_{\alpha}(s_2-t)-\widehat{\mathcal{T}}_{\alpha}(s_1-t)\right\|_{\mathcal{L}(\mathbb{X})}\int_{0}^{s_1-\nu}(s_1-s)^{\alpha-1}\gamma_{r'}(s)\mathrm{d}s\nonumber\\&\quad+\frac{2M\alpha}{\Gamma(1+\alpha)}\frac{\nu^{\alpha-\alpha_1}}{\mu^{1-\alpha_1}}\left(\int_{s_1-\nu}^{s_1}(\gamma_{r'}(s))^{\frac{1}{\alpha_1}}\mathrm{d}s\right)^{\alpha_1}.
		\end{align}
		Moreover, for $s_1,s_2\in(t_k, \tau_k], \ k=1,\ldots,m$ with $s_1<s_2$ and $x\in\mathrm{E}_r$, we have 
		\begin{align}\label{4.33}
		\left\|(F_{\lambda}x)(s_2)-(F_{\lambda}x)(s_1)\right\|_{\mathbb{X}}&\le\left\|h_k(s_2,\tilde{x}(t_k^-))-h_k(s_1,\tilde{x}(t_k^-))\right\|_{\mathbb{X}}.
		\end{align}
		Similar to the estimate \eqref{2.8}, one can easily get that the right hand side of the expressions \eqref{4.31} \eqref{4.32}  and \eqref{4.33} converge to zero as for the arbitrariness of $\nu$ and $|s_2-s_1| \rightarrow 0$.  Thus, the image of $\mathrm{E}_r$ under $F_{\lambda}$ is equicontinuous.
		
		Next,  we show that the set $\mathrm{V}(t)=\{(F_\lambda x)(t):x\in \mathrm{E}_r\}$ is relatively compact for each $t\in J$. For $t=0$, it is easy to check the set  $\mathrm{V}(t)$  is relatively compact in $\mathrm{E}_r$.  Let us take $ 0<t\leq t_1$ be fixed and for given $\eta$ with $ 0<\eta<t$ and any $\delta>0$, we define
		\begin{align*}
	&	(F_{\lambda}^{\eta,\delta}x)(t)\nonumber\\&=\int_{\delta}^{\infty}\varphi_{\alpha}(\xi)\mathcal{T}(t^{\alpha}\xi)\psi(0)\mathrm{d}\xi +\alpha\int_{0}^{t-\eta}\int_{\delta}^{\infty}\xi(t-s)^{\alpha-1}\varphi_{\alpha}(\xi)\mathcal{T}((t-s)^{\alpha}\xi)\mathrm{B}u^{\alpha}_{\lambda}(s)\mathrm{d}\xi\mathrm{d}s\nonumber\\&\quad+\alpha\int_{0}^{t-\eta}\int_{\delta}^{\infty}\xi(t-s)^{\alpha-1}\varphi_{\alpha}(\xi)\mathcal{T}((t-s)^{\alpha}\xi)f(s,\tilde{x}_{\rho(s,\tilde{x}_s)})\mathrm{d}\xi\mathrm{d}s\nonumber\\&=\mathcal{T}(\eta^{\alpha}\delta)\bigg[\int_{\delta}^{\infty}\varphi_{\alpha}(\xi)\mathcal{T}(t^{\alpha}\xi-\eta^{\alpha}\delta)\psi(0)\mathrm{d}\xi\nonumber\\&\qquad\qquad+\alpha\int_{0}^{t-\eta}\int_{\delta}^{\infty}\xi(t-s)^{\alpha-1}\varphi_{\alpha}(\xi)\mathcal{T}((t-s)^{\alpha}\xi-\eta^{\alpha}\delta)\mathrm{B}u^{\alpha}_{\lambda}(s)\mathrm{d}\xi\mathrm{d}s\nonumber\\&\qquad\qquad+\alpha\int_{0}^{t-\eta}\int_{\delta}^{\infty}\xi(t-s)^{\alpha-1}\varphi_{\alpha}(\xi)\mathcal{T}((t-s)^{\alpha}\xi-\eta^{\alpha}\delta)f(s,\tilde{x}_{\rho(s,\tilde{x}_s)})\mathrm{d}\xi\mathrm{d}s\bigg]\nonumber\\&=\mathcal{T}(\eta^{\alpha}\delta)y(t,\eta,\delta),
		\end{align*}
		where $y(\cdot,\cdot,\cdot)$ is the term appearing inside the parenthesis. Using Assumption \ref{as2.1}, one can calculate
		\begin{align*}
		\left\|y(t,\eta,\delta)\right\|_{\mathbb{X}}&\le\int_{\delta}^{\infty}\varphi_{\alpha}(\xi)\left\|\mathcal{T}(t^{\alpha}\xi-\eta^{\alpha}\delta)\psi(0)\right\|_{\mathbb{X}}\mathrm{d}\xi\nonumber\\&\quad+\alpha\int_{0}^{t-\eta}\int_{\delta}^{\infty}\xi(t-s)^{\alpha-1}\varphi_{\alpha}(\xi)\left\|\mathcal{T}((t-s)^{\alpha}\xi-\eta^{\alpha}\delta)\mathrm{B}u^{\alpha}_{\lambda}(s)\right\|_{\mathbb{X}}\mathrm{d}\xi\mathrm{d}s\nonumber\\&\quad+\alpha\int_{0}^{t-\eta}\int_{\delta}^{\infty}\xi(t-s)^{\alpha-1}\varphi_{\alpha}(\xi)\left\|\mathcal{T}((t-s)^{\alpha}\xi-\eta^{\alpha}\delta)f(s,\tilde{x}_{\rho(s,\tilde{x}_s)})\right\|_{\mathbb{X}}\mathrm{d}\xi\mathrm{d}s\nonumber\\&\le
		M\left\|\psi(0)\right\|_{\mathbb{X}}+\frac{N_0}{\lambda}\left(\frac{M\tilde{M}\alpha}{\Gamma(1+\alpha)}\right)^{2}\frac{t^{2\alpha-1}-\eta^{2\alpha-1}}{\lambda(2\alpha-1)}\nonumber\\&\quad\times\frac{M\alpha}{\Gamma(1+\alpha)}\left(\frac{t^\mu-\eta^\mu}{\mu}\right)^{1-\alpha_1}\left(\int_{0}^{t-\eta}(\gamma_{r'}(s))^{\frac{1}{\alpha_1}}\right)^{\alpha_1}<+\infty, \ \mbox{ for }\ t\in[0,t_1].
		\end{align*}
		The compactness of the operator $\mathcal{T}(\cdot)$  implies that the set  $\mathrm{V}_{\eta,\delta}(t)=\{(F^{\eta,\delta}_\lambda x)(t):x\in \mathrm{E}_r\}$ is relatively compact in $\mathbb{X}$. Hence, there exist a finite $ x_{i}$'s, for $i=1,\dots, n $ in $ \mathbb{X} $ such that 
		\begin{align*}
		\mathrm{V}_{\eta,\delta}(t) \subset \bigcup_{i=1}^{n}\mathcal{S}(x_i, \varepsilon/2),
		\end{align*}
		for some $\varepsilon>0$. Let us choose $\delta>0$ and $\eta>0$ such that 
		\begin{align*}
		&\left\|(F_{\lambda}x)(t)-(F_{\lambda}^{\eta,\delta}x)(t)\right\|_{\mathbb{X}}\nonumber\\&\le\left\|\int_{0}^{\delta}\varphi_{\alpha}(\xi)\mathcal{T}(t^{\alpha}\xi)\psi(0)\mathrm{d}\xi\right\|_{\mathbb{X}}+\alpha\left\|\int_{0}^{t}\int_{0}^{\delta}\xi(t-s)^{\alpha-1}\varphi_{\alpha}(\xi)\mathcal{T}((t-s)^{\alpha}\xi)\mathrm{B}u^{\alpha}_{\lambda}(s)\mathrm{d}\xi\mathrm{d}s\right\|_{\mathbb{X}}\nonumber\\&\quad+\alpha\left\|\int_{t-\eta}^{t}\int_{\delta}^{\infty}\xi(t-s)^{\alpha-1}\varphi_{\alpha}(\xi)\mathcal{T}((t-s)^{\alpha}\xi)\mathrm{B}u^{\alpha}_{\lambda}(s)\mathrm{d}\xi\mathrm{d}s\right\|_{\mathbb{X}}\nonumber\\&\quad+\alpha\left\|\int_{0}^{t}\int_{0}^{\delta}\xi(t-s)^{\alpha-1}\varphi_{\alpha}(\xi)\mathcal{T}((t-s)^{\alpha}\xi)f(s,\tilde{x}_{\rho(s,\tilde{x}_s)})\mathrm{d}\xi\mathrm{d}s\right\|_{\mathbb{X}}\nonumber\\&\quad+\alpha\left\|\int_{t-\eta}^{t}\int_{\delta}^{\infty}\xi(t-s)^{\alpha-1}\varphi_{\alpha}(\xi)\mathcal{T}((t-s)^{\alpha}\xi)f(s,\tilde{x}_{\rho(s,\tilde{x}_s)})\mathrm{d}\xi\mathrm{d}s\right\|_{\mathbb{X}}\nonumber\\&\le M\left\|\psi(0)\right\|_{\mathbb{X}}\int_{0}^{\delta}\varphi_{\alpha}(\xi)\mathrm{d}\xi+\frac{M^2\tilde{M}^2N_{0}\alpha t^{2\alpha-1}}{\lambda(2\alpha-1)(\Gamma(\alpha+1))}\int_{0}^{\delta}\xi\varphi_{\alpha}(\xi)\mathrm{d}\xi
		\nonumber\\&\quad+\frac{N_0}{\lambda}\left(\frac{M\tilde{M}\alpha}{\Gamma(1+\alpha)}\right)^{2}\int_{t-\eta}^{t}(t-s)^{\alpha-1}(t_1-s)^{\alpha-1}\mathrm{d}s\nonumber\\&\quad+\frac{M\alpha t^{\alpha-\alpha_1}}{\mu^{1-\alpha_1}}\left\|\gamma_{r'}\right\|_{\mathrm{L}^{\frac{1}{\alpha_1}}(J;\mathbb{R^+})}\int_{0}^{\delta}\xi\varphi_{\alpha}(\xi)\mathrm{d}\xi+\frac{M\alpha}{(\Gamma(1+\alpha))^2}\int_{t-\eta}^{t}(t-s)^{\alpha-1}\gamma_{r'}(s)\mathrm{d}s\nonumber\\&\le\frac{\varepsilon}{2}.
		\end{align*}
		Consequently $$\mathrm{V}(t)\subset \bigcup_{i=1}^{n}\mathcal{S}(x_i, \varepsilon ).$$
		Thus, for each $t\in [0,t_1]$, the set $\mathrm{V}(t)$ is relatively compact in $ \mathbb{X}$. Next, we take $ t\in(\tau_k,t_{k+1}],$ for $k=1,\ldots,m$ and for given $\eta$ with $ 0<\eta<\min\{t-\tau_k, \tau_k\}$ and any $\delta>0$, we define
		\begin{align*}
		&(F_{\lambda}^{\eta,\delta}x)(t)\nonumber\\&=\int_{\delta}^{\infty}\varphi_{\alpha}(\xi)\mathcal{T}((t-\tau_k)^{\alpha}\xi)h_k(\tau_k,\tilde{x}(t_k^-))\mathrm{d}\xi \nonumber\\&\quad-\alpha\int_{0}^{\tau_k-\eta}\int_{\delta}^{\infty}\xi(\tau_k-s)^{\alpha-1}\varphi_{\alpha}(\xi)\mathcal{T}((\tau_k-s)^{\alpha}\xi)\left[\mathrm{B}u^{\alpha}_{\lambda}(s)+(s,\tilde{x}_{\rho(s,\tilde{x}_s)})\right]\mathrm{d}\xi\mathrm{d}s\nonumber\\&\quad+\alpha\int_{0}^{t-\eta}\int_{\delta}^{\infty}\xi(t-s)^{\alpha-1}\varphi_{\alpha}(\xi)\mathcal{T}((t-s)^{\alpha}\xi)\left[\mathrm{B}u^{\alpha}_{\lambda}(s)+(s,\tilde{x}_{\rho(s,\tilde{x}_s)})\right]\mathrm{d}\xi\mathrm{d}s\nonumber\\&=\mathcal{T}(\eta^{\alpha}\delta)\bigg[\int_{\delta}^{\infty}\varphi_{\alpha}(\xi)\mathcal{T}((t-\tau_k)^{\alpha}\xi-\eta^{\alpha}\delta)h_k(\tau_k,\tilde{x}(t_k^-))\mathrm{d}\xi\nonumber\\&\quad+\alpha\int_{0}^{\tau_k-\eta}\int_{\delta}^{\infty}\xi(\tau_k-s)^{\alpha-1}\varphi_{\alpha}(\xi)\mathcal{T}((\tau_k-s)^{\alpha}\xi-\eta^{\alpha}\delta)\left[\mathrm{B}u^{\alpha}_{\lambda}(s)+(s,\tilde{x}_{\rho(s,\tilde{x}_s)})\right]\mathrm{d}\xi\mathrm{d}s\nonumber\\&\quad+\alpha\int_{0}^{t-\eta}\int_{\delta}^{\infty}\xi(t-s)^{\alpha-1}\varphi_{\alpha}(\xi)\mathcal{T}((t-s)^{\alpha}\xi-\eta^{\alpha}\delta)\left[\mathrm{B}u^{\alpha}_{\lambda}(s)+(s,\tilde{x}_{\rho(s,\tilde{x}_s)})\right]\mathrm{d}\xi\mathrm{d}s\bigg].
		\end{align*}
		Proceeding similarly for the case $t\in[0,t_1]$, one can prove that the set $\mathrm{V}(t),$ for $t\in(\tau_k,t_{k+1}],\ k=1,\ldots,m$ is relatively compact in $ \mathbb{X}$. Moreover for $t\in(t_k,\tau_k], k=1\ldots,m$, the fact that the set $\mathrm{V}(t)$ is relative compact follows by the compactness of the impulses $h_k,$ for $k=1,\ldots,m$. Therefore, the set $\mathrm{V}(t)=\{(F_\lambda x)(t):x\in \mathrm{E}_r\},$ for each $t\in J$ is relatively compact in $\mathbb{X}$.
		
		Hence, by invoking the Arzela-Ascoli theorem, we conclude that the operator $ F_{\lambda}$  is compact. Then \emph{Schauder's fixed point theorem} yields that the operator $F_{\lambda}$ has a fixed point in $\mathrm{E}_{r}$, which is a mild solution of the system \eqref{1.1}.
	\end{proof}
	
	In order to prove the approximate controllability of the system \eqref{1.1}, we replace the assumption (\textit{H2}) by the following stronger assumption: 
	\begin{enumerate}\label{as}
		\item [\textit{(H4)}] The function $ f: J \times \mathfrak{B}_{\alpha} \rightarrow \mathbb{X} $ satisfies the assumption (\textit{$H2$})(i) and there exists a function $ \gamma\in \mathrm{L}^{\frac{1}{\alpha_1}}(J;\mathbb{R}^+)$ with $\alpha_1\in[0,\frac{1}{2})$ such that $$ \|f(t,\psi)\|_{\mathbb{X}}\leq \gamma(t),\ \text{ for all }\  (t,\psi) \in J \times  \mathfrak{B}_{\alpha}. $$ 
	\end{enumerate}
	\begin{theorem}\label{thm4.4}
		Suppose that Assumptions (H0)-(H1), (H3)-(H4) and the condition \eqref{cnd} of Theorem \ref{thm4.3} are satisfied. Then the system \eqref{1.1} is approximately controllable.
	\end{theorem}
	\begin{proof}
		By using Theorem \ref{thm4.3}, we infer that for every $\lambda>0$ and $\zeta_k\in \mathbb{X}$ for $k=0,1,\ldots m$, there exists a mild solution $x^{\lambda}\in\mathrm{E}_r$ such that
		\begin{equation}\label{M}
		x^{\lambda}(t)=\begin{dcases}
		\mathcal{T}_{\alpha}(t)\psi(0)+\int_{0}^{t}(t-s)^{\alpha-1}\widehat{\mathcal{T}}_{\alpha}(t-s)\left[\mathrm{B}u^{\alpha}_{\lambda}(s)+f(s,\tilde{x^\lambda}_{\rho(s, \tilde{x^\lambda}_s)})\right]\mathrm{d}s,t\in[0, t_1],\\
		h_k(t, \tilde{x^\lambda}(t_k^-)),t\in(t_k, \tau_k],\ k=1,\ldots,m,\\
		\mathcal{T}_{\alpha}(t-\tau_k)h_k(\tau_k,\tilde{x^\lambda}(t_k^-))-\int_{0}^{\tau_k}(\tau_k-s)^{\alpha-1}\widehat{\mathcal{T}}_{\alpha}(\tau_k-s)\left[\mathrm{B}u^{\alpha}_{\lambda}(s)+f(s,\tilde{x^\lambda}_{\rho(s,\tilde{x^\lambda}_s)})\right]\mathrm{d}s\\\quad+\int_{0}^{t}(t-s)^{\alpha-1}\widehat{\mathcal{T}}_{\alpha}(t-s)\left[\mathrm{B}u^{\alpha}_{\lambda}(s)+f(s,\tilde{x^\lambda}_{\rho(s, \tilde{x^\lambda}_s)})\right]\mathrm{d}s,\\ \qquad\qquad\qquad  t\in(\tau_k,t_{k+1}],\ k=1,\ldots,m,
		\end{dcases}
		\end{equation}
		with the control defined in \eqref{C}. Next, we estimate
		\begin{align}\label{4.35}
		x^{\lambda}(T)&=\mathcal{T}_{\alpha}(T-\tau_m)h_{m}(\tau_m,\tilde{x^\lambda}(t_m^-))\nonumber\\&\quad-\int_{0}^{\tau_m}(\tau_m-s)^{\alpha-1}\widehat{\mathcal{T}}_{\alpha}(\tau_m-s)\left[\mathrm{B}u^{\alpha}_{\lambda}(s)+f(s,\tilde{x^\lambda}_{\rho(s,\tilde{x^\lambda}_s)})\right]\mathrm{d}s\nonumber\\&\quad+\int_{0}^{T}(T-s)^{\alpha-1}\widehat{\mathcal{T}}_{\alpha}(T-s)\left[\mathrm{B}u^{\alpha}_{\lambda}(s)+f(s,\tilde{x^\lambda}_{\rho(s,\tilde{x^\lambda}_s)})\right]\mathrm{d}s\nonumber\\&=\mathcal{T}_{\alpha}(T-\tau_m)h_{m}(\tau_m,\tilde{x^\lambda}(t_m^-))\nonumber\\&\quad-\int_{0}^{\tau_m}(\tau_m-s)^{\alpha-1}\widehat{\mathcal{T}}_{\alpha}(\tau_m-s)\left[\mathrm{B}u^{\alpha}_{\lambda}(s)+f(s,\tilde{x^\lambda}_{\rho(s,\tilde{x^\lambda}_s)})\right]\mathrm{d}s\nonumber\\&\quad+\int_{0}^{T}(T-s)^{\alpha-1}\widehat{\mathcal{T}}_{\alpha}(T-s)f(s,\tilde{x^\lambda}_{\rho(s,\tilde{x^\lambda}_s)})+\int_{0}^{\tau_m}(T-s)^{\alpha-1}\widehat{\mathcal{T}}_{\alpha}(T-s)\mathrm{B}u^{\alpha}_{\lambda}(s)\mathrm{d}s\nonumber\\&\quad+\int_{\tau_m}^{T}(T-s)^{2(\alpha-1)}\widehat{\mathcal{T}}_{\alpha}(T-s)\mathrm{B}\mathrm{B}^*\widehat{\mathcal{T}}_{\alpha}(T-s)^*\mathcal{J}\left[\mathcal{R}(\lambda,\Phi_{\tau_m}^{T})p_m(x^\lambda(\cdot))\right]\mathrm{d}s\nonumber\\&=\mathcal{T}_{\alpha}(T-\tau_m)h_{m}(\tau_m,\tilde{x^\lambda}(t_m^-))\nonumber\\&\quad-\int_{0}^{\tau_m}(\tau_m-s)^{\alpha-1}\widehat{\mathcal{T}}_{\alpha}(\tau_m-s)\left[\mathrm{B}u^{\alpha}_{\lambda}(s)+f(s,\tilde{x^\lambda}_{\rho(s,\tilde{x^\lambda}_s)})\right]\mathrm{d}s\nonumber\\&\quad+\int_{0}^{T}(T-s)^{\alpha-1}\widehat{\mathcal{T}}_{\alpha}(T-s)f(s,\tilde{x^\lambda}_{\rho(s,\tilde{x^\lambda}_s)})+\int_{0}^{\tau_m}(T-s)^{\alpha-1}\widehat{\mathcal{T}}_{\alpha}(T-s)\mathrm{B}u^{\alpha}_{\lambda}(s)\mathrm{d}s\nonumber\\&\quad+\Phi_{\tau_m}^{T}\mathcal{J}\left[\mathcal{R}(\lambda,\Phi_{\tau_m}^{T})p_m(x^\lambda(\cdot))\right]\nonumber\\&=\zeta_m-\lambda\mathcal{R}(\lambda,\Phi_{\tau_m}^{T})p_m(x^\lambda(\cdot)).
		\end{align}
		Moreover,  that fact that $ x^{\lambda}\in\mathrm{E}_r $ is bounded implies that the  sequence $ x^{\lambda}(t),$  for each $t\in J$ is bounded in $\mathbb{X}$. Then by the Banach-Alaoglu theorem, we can find a subsequence, still denoted as $ x^{\lambda},$ such that 
		\begin{align*}
		x^\lambda(t)\xrightharpoonup{w}z(t) \ \mbox{ in }\ \mathbb{X} \ \ \mbox{as}\  \ \lambda\to0^+,\ t\in J.
		\end{align*}
		Using the condition (\textit{$H3$}) of Assumption \ref{as2.1}, we obtain
		\begin{align}\label{4.19}
		h_m(t,x^\lambda(t_m^-))\to h_m(t,z(t_m^-)) \ \mbox{ in }\ \mathbb{X} \ \ \mbox{as}\  \ \lambda\to0^+.
		\end{align}
		Furthermore, by using Assumption  \textit{(H4)}, we get
		\begin{align}
		\int_{s_1}^{s_2}\left\|f(s,\tilde{x^{\lambda}}_{\rho(s,\tilde{x^{\lambda}}_s)})\right\|_{\mathbb{X}}^{2}\mathrm{d}s&\le \int_{s_1}^{s_2}\gamma^2(s)\mathrm{d} s\leq \left(\int_{s_1}^{s_2}\gamma^{\frac{1}{\alpha_1}}(s)\mathrm{d}s\right)^{2\alpha_1}(s_1-s_2)^{1-2\alpha_1}<+\infty, \nonumber
		\end{align}
		for any $ s_1,s_2\in[0,T]$ with $s_1<s_2$. Therefore, the sequence $ \{f(\cdot, \tilde{x^{\lambda}}_{\rho(s,\tilde{x^{\lambda}}_s)}): \lambda >0\} $  in $ \mathrm{L}^2([s_1,s_2]; \mathbb{X})$ is bounded. By an application of the Banach-Alaoglu theorem, we can find a subsequence still denoted as $ \{f(\cdot, \tilde{x^{\lambda}}_{\rho(s,\tilde{x^{\lambda}}_s)}): \lambda > 0 \}$ such that 
		\begin{align}\label{4.36}
		f(\cdot, \tilde{x^{\lambda}}_{\rho(s,\tilde{x^{\lambda}}_s)})\xrightharpoonup{w}f(\cdot) \ \mbox{ in }\ \mathrm{L}^2([s_1,s_2];\mathbb{X}).
		\end{align}
		We now calculate 
		\begin{align}
		&\int_{0}^{\tau_m}\left\|u^{\alpha}_{\lambda}(s)\right\|^2_{\mathbb{U}}\mathrm{d}s\nonumber\\&=\int_{0}^{t_1}\left\|u^{\alpha}_{\lambda}(s)\right\|^2_{\mathbb{U}}\mathrm{d}s+\int_{t_1}^{\tau_1}\left\|u^{\alpha}_{\lambda}(s)\right\|^2_{\mathbb{U}}\mathrm{d}s+\int_{\tau_1}^{t_2}\left\|u^{\alpha}_{\lambda}(s)\right\|^2_{\mathbb{U}}\mathrm{d}s+\cdots+\int_{t_m}^{\tau_m}\left\|u^{\alpha}_{\lambda}(s)\right\|^2_{\mathbb{U}}\mathrm{d}s\nonumber\\&=\int_{0}^{t_1}\left\|u^{\alpha}_{0,\lambda}(s)\right\|^2_{\mathbb{U}}\mathrm{d}s+\int_{\tau_1}^{t_2}\left\|u^{\alpha}_{1,\lambda}(s)\right\|^2_{\mathbb{U}}\mathrm{d}s+\cdots+\int_{\tau_{m-1}}^{t_m}\left\|u^{\alpha}_{m-1,\lambda}(s)\right\|^2_{\mathbb{U}}\mathrm{d}s\nonumber\\&\le\left(\frac{M\tilde{M}\alpha}{\lambda\Gamma(1+\alpha)}\right)^{2}\frac{T^{2\alpha-1}}{2\alpha-1}\sum_{j=0}^{m-1}C_j^2=C,
		\end{align}
		where $C_k,$ for $k=1,\ldots,m-1$ are the same as given in \eqref{4.5} and $C_0=N_0$ given in \eqref{4.4}. Moreover, the above estimate ensures that the sequence $\{u^\alpha_{\lambda}(\cdot): \lambda >0\}$  in $ \mathrm{L}^2([0,\tau_m]; \mathbb{U})$ is bounded. Further, by the Banach-Alaoglu theorem, we can find a subsequence, still denoted as $\{u^\alpha_{\lambda}(\cdot): \lambda >0\}$ such that 
		\begin{align}\label{4}
		u^\alpha_{\lambda}(\cdot)\xrightharpoonup{w}u^\alpha(\cdot) \ \mbox{ in }\ \mathrm{L}^2([0,\tau_m];\mathbb{U}).
		\end{align}
		Next, we compute
		\begin{align}\label{4.37}
		\left\|p_{m}(x^{\lambda}(\cdot))-\omega\right\|_{\mathbb{X}}&\le\left\|\mathcal{T}_{\alpha}(T-\tau_m)(h_k(\tau_m,\tilde{x^\lambda}(t_m^-))-h_k(\tau_m,z(t_m^-)))\right\|_{\mathbb{X}}\nonumber\\&\quad+\left\|\int_{0}^{\tau_m}(\tau_m-s)^{\alpha-1}\widehat{\mathcal{T}}_{\alpha}(\tau_m-s)\mathrm{B}\left[u^\alpha_{\lambda}(s)-u^{\alpha}(s)\right]\mathrm{d}s\right\|_{\mathbb{X}}\nonumber\\&\quad+\left\|\int_{0}^{\tau_m}(\tau_m-s)^{\alpha-1}\widehat{\mathcal{T}}_{\alpha}(\tau_m-s)\left[f(s,\tilde{x^\lambda}_{\rho(s,\tilde{x^\lambda}_s)})-f(s)\right]\mathrm{d}s\right\|_{\mathbb{X}}\nonumber\\&\quad+\left\|\int_{0}^{\tau_m}(T-s)^{\alpha-1}\widehat{\mathcal{T}}_{\alpha}(T-s)\mathrm{B}\left[u^\alpha_{\lambda}(s)-u^{\alpha}(s)\right]\mathrm{d}s\right\|_{\mathbb{X}}\nonumber\\&\quad+\left\|\int_{0}^{T}(T-s)^{\alpha-1}\widehat{\mathcal{T}}_{\alpha}(T-s)\left[f(s,\tilde{x^\lambda}_{\rho(s,\tilde{x^\lambda}_s)})-f(s)\right]\mathrm{d}s\right\|_{\mathbb{X}}\nonumber\\&\le\left\|\mathcal{T}_{\alpha}(T-\tau_m)(h_k(\tau_m,x^\lambda(t_m^-))-h_k(\tau_m,z(t_m^-)))\right\|_{\mathbb{X}}\nonumber\\&\quad+\left\|\int_{0}^{\tau_m}(\tau_m-s)^{\alpha-1}\widehat{\mathcal{T}}_{\alpha}(\tau_m-s)\mathrm{B}\left[u^\alpha_{\lambda}(s)-u^{\alpha}(s)\right]\mathrm{d}s\right\|_{\mathbb{X}}\nonumber\\&\quad+\left\|\int_{0}^{\tau_m}(\tau_m-s)^{\alpha-1}\widehat{\mathcal{T}}_{\alpha}(\tau_m-s)\left[f(s,\tilde{x^\lambda}_{\rho(s,\tilde{x^\lambda}_s)})-f(s)\right]\mathrm{d}s\right\|_{\mathbb{X}}\nonumber\\&\quad+\frac{T^{2\alpha-1}-(T-\tau_m)^{2\alpha-1}}{2\alpha-1}\left(\int_{0}^{\tau_m}\left\|\widehat{\mathcal{T}}_{\alpha}(T-s)\mathrm{B}\left[u^\alpha_{\lambda}(s)-u^{\alpha}(s)\right]\right\|^2_{\mathbb{X}}\mathrm{d}s\right)^{\frac{1}{2}}\nonumber\\&\quad+\left\|\int_{0}^{T}(T-s)^{\alpha-1}\widehat{\mathcal{T}}_{\alpha}(T-s)\left[f(s,\tilde{x^\lambda}_{\rho(s,\tilde{x^\lambda}_s)})-f(s)\right]\mathrm{d}s\right\|_{\mathbb{X}}\nonumber\\&\quad\to 0\ \mbox{as}\ \lambda\to0^+, 
		\end{align}
		where 
		\begin{align*}
		\omega &=\zeta_m-\mathcal{T}_{\alpha}(T-\tau_m)h_m(\tau_m,z(t_m^-))+\int_{0}^{\tau_m}(\tau_m-s)^{\alpha-1}\widehat{\mathcal{T}}_{\alpha}(\tau_m-s)\left[\mathrm{B}u^{\alpha}(s)+f(s)\right]\mathrm{d}s\nonumber\\&\quad-\int_{0}^{\tau_m}(T-s)^{\alpha-1}\widehat{\mathcal{T}}_{\alpha}(T-s)\mathrm{B}u^{\alpha}(s)\mathrm{d}s-\int_{0}^{T}(T-s)^{\alpha-1}\widehat{\mathcal{T}}_{\alpha}(T-s)f(s)\mathrm{d}s.
		\end{align*}
		Here, we used the convergences \eqref{4.19},\eqref{4.36},\eqref{4}, the dominated convergence theorem and the compactness of the operator $f(\cdot)\to\int_{0}^{\cdot}(\cdot-s)^{\alpha-1}\widehat{\mathcal{T}}_{\alpha}(\cdot-s) f(s)\mathrm{d}s:\mathrm{L}^2(J;\mathbb{X})\rightarrow \mathrm{C}(J;\mathbb{X})$, (see Lemma \ref{lem2.12}).
			Finally, by using the equality \eqref{4.35}, we evaluate 
		\begin{align}
		\left\|x^{\lambda}(T)-\zeta_m\right\|_{\mathbb{X}}&\le\left\|\lambda\mathcal{R}(\lambda,\Phi_{\tau_m}^{T})p_m(x(\cdot))\right\|_{\mathbb{X}}\nonumber\\&\le\left\|\lambda\mathcal{R}(\lambda,\Phi_{\tau_m}^{T})(p_m(x(\cdot))-\omega)\right\|_{\mathbb{X}}+\left\|\lambda\mathcal{R}(\lambda,\Phi_{\tau_m}^{T})\omega\right\|_{\mathbb{X}}\nonumber\\&\le\left\|\lambda\mathcal{R}(\lambda,\Phi_{\tau_m}^{T})\right\|_{\mathcal{L}(\mathbb{X})}\left\|p_m(x(\cdot))-\omega\right\|_{\mathbb{X}}+\left\|\lambda\mathcal{R}(\lambda,\Phi_{\tau_m}^{T})\omega\right\|_{\mathbb{X}}.
		\end{align}
		Using the above inequality, \eqref{4.37} and Assumption \ref{as2.1} (\textit{H0}), we obtain
		\begin{align*}
		\left\|x^{\lambda}(T)-\zeta_m\right\|_{\mathbb{X}}\to0,\ \mbox{ as }\ \lambda\to0^+,
		\end{align*}
		which ensures that the system \eqref{1.1} is approximately controllable on $J$.
	\end{proof}
	\begin{rem}\label{rem4.4}
		The works \cite{PCY,NIZ,SRY}, etc considered a different kind of control for the fractional order semilinear problems. If one follows Remark \ref{rem3.6}, the controllability operator defined in \eqref{2.1} changes to \eqref{3.20}
		and  $u^\alpha_{k,\lambda}(\cdot)$ appearing in the control defined in \eqref{C} takes the form 
		\begin{align*}
		u^\alpha_{k,\lambda}(t)&=\mathrm{B}^*\widehat{\mathcal{T}}_{\alpha}(t_{k+1}-t)^*\mathcal{J}\left[\mathcal{R}(\lambda,\Phi_{\tau_k}^{t_{k+1}})p_k(x(\cdot))\right],
		\end{align*}
		for $t\in [\tau_k, t_{k+1}),k=0,1,\ldots,m$. The control provided in \eqref{C} is motivated from the linear regulator problem with the cost functional defined in \eqref{3.1}.  The proof of Theorems \ref{thm4.3} and \ref{thm4.4} follows in a similar way with some obvious modifications in the calculations. 
	\end{rem}
	\section{Application}\label{application}\setcounter{equation}{0}
	In this section, we discuss a concrete example to verify the results developed in previous sections. 
	
	\begin{Ex}\label{ex1} Let us take the following fractional heat equation with non-instantaneous impulses and delay:
		\begin{equation}\label{ex}
		\left\{
		\begin{aligned}
		\frac{\partial^\alpha z(t,\xi)}{\partial t^\alpha}&=\frac{\partial^2z(t,\xi)}{\partial \xi^2}+\eta(t,\xi)+\int_{-\infty}^{t}b(s-t)z(s-\sigma(\|z(t)\|),\xi)\mathrm{d}s, \\&\qquad \qquad \ t\in\bigcup_{k=0}^{m} (\tau_k, t_{k+1}]\subset J=[0,T], \ \xi\in[0,\pi], \\
		z(t,\xi)&=h_k(t,z(t_k^-,\xi)),\ t\in(t_k,\tau_k],\ k=1,\ldots, m,\ \xi\in[0,\pi],\\
		z(t,0)&=0=z(t,\pi), \  t\in [0, 1], \\
		z(\theta,\xi)&=\psi(\theta,\xi), \ \xi\in[0,\pi], \ \theta\leq0.
		\end{aligned}
		\right.
		\end{equation}
		where the function $\eta:[0,1]\times[0,\pi]\to[0,\pi]$  is continuous in $t$ and the functions $\sigma:[0,\infty)\to[0,\infty)$ are also continuous.
	\end{Ex}
	\vskip 0.1 cm
	\noindent\textbf{Step 1:} \emph{$\mathrm{C}_0$-semigroup and phase space:} 
	Let $\mathbb{X}_p= \mathrm{L}^{p}([0,\pi];\mathbb{R})$ with $p\in[2,\infty)$, and $\mathbb{U}=\mathrm{L}^{2}([0,\pi];\mathbb{R})$. Note that  $\mathbb{X}_p$ is separable and reflexive with strictly convex dual $\mathbb{X}_p^*=\mathrm{L}^{\frac{p}{p-1}}([0,\pi];\mathbb{R})$ and $\mathbb{U}$ is separable. We define the linear operator  $\mathrm{A}_p:\mathrm{D}(\mathrm{A}_p)\subset\mathbb{X}_p\to\mathbb{X}_p$  as
	\begin{align*}
	\mathrm{A}_pg(\xi)= g''(\xi),
	\end{align*}
	where $\mathrm{D}(\mathrm{A}_p)= \mathrm{W}^{2,p}([0,\pi];\mathbb{R})\cap\mathrm{W}_0^{1,p}([0,\pi];\mathbb{R})$. Since we know that  $\mathrm{C}_0^{\infty}([0,\pi];\mathbb{R})\subset\mathrm{D}(\mathrm{A}_p)$ and hence $\mathrm{D}(\mathrm{A}_p)$ is dense in $\mathbb{X}_p$ and one can easily verify that the operator $\mathrm{A}_p$ is closed. Next, we consider the following Sturm-Liouville system:
	\begin{equation}\label{59}
	\left\{
	\begin{aligned}
	\left(\lambda\mathrm{I}-\mathrm{A}_p\right)g(\xi)&=l(\xi), \ 0<\xi<\pi,\\
	g(0)=g(\pi)&=0.
	\end{aligned}
	\right.
	\end{equation}
	One can easily  rewrite the above system as 
	\begin{align}\label{511}
	\left(\lambda\mathrm{I}-\Delta\right)g(\xi)&=l(\xi),
	\end{align}
	where $\Delta g(\xi)=g''(\xi)$. Multiplying both sides of \eqref{511} by $g|g|^{p-2}$ and then integrating over $[0,\pi]$, we obtain
	\begin{align}\label{536}
	\lambda\int_0^{\pi}|g(\xi)|^p\mathrm{d}\xi&+(p-1)\int_0^{\pi}|g(\xi)|^{p-2}|f'(\xi)|^2\mathrm{d}\xi=\int_0^{\pi}l(\xi)g(\xi)|g(\xi)|^{p-2}\mathrm{d}\xi.
	\end{align}
	%Using the estimate 
	%\begin{align*}
	%	\|(|g|^{p/2})'\|_{\mathrm{L}^2}=\frac{p}{2}\||g|^{\frac{p-2}{2}}g'\|_{\mathrm{L}^2},
	%\end{align*}
	%and the Poincar\'e inequality $\|g\|_{\mathrm{L}^2}\leq \|g'\|_{\mathrm{L}^2}$ (since the Poincar\'e constant is $1$), we get 
	%\begin{align}\label{537}
	%	\|g\|_{\mathrm{L}^p}^p=\||g|^{p/2}\|_{\mathrm{L}^2}^2\leq\|(|g|^{p/2})'\|_{\mathrm{L}^2}^2=\frac{p^2}{4}\||g|^{\frac{p-2}{2}}g'\|_{\mathrm{L}^2}^2.
	%\end{align}
	%Substituting \eqref{537} in \eqref{536}, we find 
	%\begin{align}
	%	\int_0^{\pi}|g(\xi)|^p\mathrm{d}\xi&+\frac{4(p-1)}{p^2}\int_0^{\pi}|g(\xi)|^{p}\mathrm{d}\xi\leq\int_0^{\pi}|l(\xi)||g(\xi)|^{p-1}\mathrm{d}\xi.
	%\end{align}
	Applying H\"older's inequality, we get
	\begin{align*}
	\lambda\int_0^{\pi}|g(\xi)|^p\mathrm{d}\xi&\leq \left(\int_0^{\pi}|g(\xi)|^p\mathrm{d}\xi\right)^{\frac{p-1}{p}}\left(\int_0^{\pi}|l(\xi)|^p\mathrm{d}\xi\right)^{\frac{1}{p}}.
	\end{align*}
	Thus, we have 
	\begin{align*}
	\|\mathcal{R}(\lambda,\mathrm{A}_p)l\|_{\mathrm{L}^p}=\|g\|_{\mathrm{L}^p}\leq\frac{1}{\lambda}\|l\|_{\mathrm{L}^p},
	\end{align*}
	so that we obtain  
	\begin{align}
	\|\mathcal{R}(\lambda,\mathrm{A}_p\|_{\mathcal{L}(\mathrm{L}^p)}\leq\frac{1}{\lambda}.
	\end{align}
	Hence, by applying the Hille-Yosida theorem, we obtain that the operator $\mathrm{A}_p$ generate a strongly continuous semigroup $\{\mathcal{T}_p(t):t\ge0\}$ of bounded linear operators. 
	
	Moreover, the infinitesimal generator $\mathrm{A}_p$ and the semigroup $\mathcal{T}_p(t)$ can be written as
	\begin{align}
	\mathrm{A}_pg&= \sum_{n=1}^{\infty}-n^{2}\langle g, w_{n} \rangle  w_{n},\ g\in \mathrm{D}(\mathrm{A}_p),\nonumber\\
	\mathcal{T}_p(t)g&= \sum_{n=1}^{\infty}\exp\left(-n^2t\right)\langle g, w_{n} \rangle  w_{n},\ g\in\mathbb{X}_p,
	\end{align}
	where, $w_n(\xi)=\sqrt{\frac{2}{\pi}}\sin(n\xi)$ are the normalized eigenfunctions corresponding to the eigenvalues $\lambda_n=-n^2\ (n\in\mathbb{N})$ of the operator $\mathrm{A}_p$ and $\langle g,w_n\rangle :=\int_0^{\pi}g(\xi)w_n(\xi)\mathrm{d}\xi$. Further, the resolvent operator $\mathcal{R}(\lambda,\mathrm{A}_p)$ is compact (see \cite{MTM} for more details). Therefore, the generated semigroup $\mathcal{T}_p(t)$ is compact for $t>0$. Thus, the condition (\textit{H1}) of Assumption \ref{as2.1} holds.
	
	We now define the following operators
	\begin{align}
	\mathcal{T}_{\alpha,p}(t)g&=\int_{0}^{\infty}\varphi_{\alpha}(\xi)\mathcal{T}_p(t^{\alpha}\xi)g(\xi)\mathrm{d}\xi=\int_{0}^{\infty}\varphi_{\alpha}(\xi)\sum_{n=1}^{\infty}\exp\left(-n^2t^\alpha\xi\right)\langle g, w_{n} \rangle  w_{n}(\xi)\mathrm{d}\xi.\\
	\widehat{\mathcal{T}}_{\alpha,p}(t)g&=\alpha\int_{0}^{\infty}\xi\varphi_{\alpha}(\xi)\mathcal{T}_p(t^{\alpha}\xi)g(\xi)\mathrm{d}\xi=\alpha\int_{0}^{\infty}\xi\varphi_{\alpha}(\xi)\sum_{n=1}^{\infty}\exp\left(-n^2t^\alpha\xi\right)\langle g, w_{n} \rangle  w_{n}(\xi)\mathrm{d}\xi, \label{5.7}
	\end{align}
	for all $g\in\mathbb{X}_p$.  
	
	Let us take $\mathfrak{B}=\mathrm{PC}_{0}\times\mathrm{L}^1_h(\mathbb{X})$ with  $h(\theta)=e^{\nu\theta}$, for some $\nu>0$ (see Example \ref{exm2.8}). Proceeding similar arguments as in section 5,  \cite{SSM}, one can verify that the space $\mathfrak{B}=\mathrm{PC}_{0}\times\mathrm{L}^1_h(\mathbb{X})$  is a phase space, which satisfies the axioms (A1) and (A2) with $\Lambda(t) =\int_{-t}^{-0}h(\theta)\mathrm{d}\theta$ and $\Upsilon(t)= H(-t)$.
	We define  $K:=\sup\limits_{\theta\in(-\infty,0]}\frac{|b(-\theta)|}{h(\theta)}$.
	\vskip 0.1 cm 
	\noindent\textbf{Step 2:} \emph{Abstract formulation and approximate controllability.}
	Let us define $$x(t)(\xi):=z(t,\xi),\ \mbox{ for }\ t\in J\ \mbox{ and }\ \xi\in[0,\pi],$$ and the bounded linear operator $\mathrm{B}:\mathbb{U}\to\mathbb{X}_p$ as  $$\mathrm{B}u(t)(\xi):=\eta(t,\xi)=\int_{0}^{\pi}K(\zeta,\xi)u(t)(\zeta)\mathrm{d}\zeta, \ t\in J,\ \xi\in [0,\pi],$$ where $K\in\mathrm{C}([0,\pi]\times[0,\pi];\mathbb{R})$ with $K(\zeta,\xi)=K(\xi,\zeta),$ for all $\zeta,\xi\in [0,\pi]$. We assume that the operator $\mathrm{B}$ is one-one.
	Let us estimate 
	\begin{align*}
	\left\|\mathrm{B}u(t)\right\|_{\mathbb{X}_p}^p=\int_{0}^{\pi}\left|\int_{0}^{\pi}K(\zeta,\xi)u(t)(\zeta)\mathrm{d}\zeta\right|^p\mathrm{d}\xi.
	\end{align*}
	Applying the Cauchy-Schwarz inequality, we have
	\begin{align*}
	\left\|\mathrm{B}u(t)\right\|_{\mathbb{X}_p}^p&\le\int_{0}^{\pi}\left[\left(\int_{0}^{\pi}|K(\zeta,\xi)|^2\mathrm{d}\zeta\right)^{\frac{1}{2}}\left(\int_{0}^{\pi}|u(t)(\zeta)|^2\mathrm{d}\zeta\right)^{\frac{1}{2}}\right]^{p}\mathrm{d}\xi\\&=\left(\int_{0}^{\pi}|u(t)(\zeta)|^2\mathrm{d}\zeta\right)^{\frac{p}{2}}\int_{0}^{\pi}\left(\int_{0}^{\pi}|K(\zeta,\xi)|^2\mathrm{d}\zeta\right)^{\frac{p}{2}}\mathrm{d}\xi.
	\end{align*}
	Since the kernel $K(\cdot,\cdot)$ is continuous, we arrive at
	\begin{align*}
	\left\|\mathrm{B}u(t)\right\|_{\mathbb{X}_p}\le C\left\|u(t)\right\|_{\mathbb{U}},
	\end{align*}
	so that we get 
	$	\left\|\mathrm{B}\right\|_{\mathcal{L}(\mathbb{U};\mathbb{X}_p)}\le C.$
	Hence, the operator $\mathrm{B}$ is bounded. Moreover, the symmetry of the kernel  implies that the operator $\mathrm{B}=\mathrm{B}^*$ (self adjoint). For example, one can take $K(\xi,\zeta)=1+\xi^2+\zeta^2,\ \mbox{for all}\ \xi, \zeta\in [0,\pi]$.
	The function $\psi:(-\infty,0]\rightarrow\mathbb{X}$ is given as
	\begin{align}
	\nonumber \psi(t)(\xi)=\psi(t,\xi),\ \xi\in[0,\pi].
	\end{align}	 
	Next, the functions $f, \rho:J\times \mathfrak{B}\to\mathbb{X}$ are defined as
	\begin{align}
	\nonumber f(t,\psi)\xi&:=\int_{-\infty}^{0}b(-\theta)\psi(\theta,\xi)\mathrm{d}\theta,\\
	\nonumber\rho(t,\psi):&=t-\sigma(\|\psi(0)\|_{\mathbb{X}}),
	\end{align}	
	for $\xi\in[0,\pi]$. Clearly, $f$ is continuous and uniformly bounded by $K$. These facts guarantee that the function $f$ satisfied the condition \textit{$(H2)$} of Assumption \ref{as2.1} and the condition \textit{$(H4)$}.
	
	Moreover, the impulse functions $h_k:[t_k,\tau_k]\times\mathbb{X}\to\mathbb{X},$ for $k=1,\ldots,m,$ are defined as 
	\begin{align*}
	h_k(t,x)\xi:=\int_{0}^{\pi}\rho_k(t,\xi,z)\cos^2(x(t_k^-)z)\mathrm{d}z, \ \mbox{ for }\ t\in(t_k,\tau_k],
	\end{align*}
	where, $\rho_k\in\mathrm{C}([0,1]\times[0,\pi]^2;\mathbb{R})$. It is easy to verify that the impulses $h_k,$ for $k=1,\ldots,m,$ satisfy the condition \textit{$(H3)$} of Assumption \ref{as2.1}.

	The system \eqref{ex} can be transformed into the abstract form \eqref{1.1} by using the above substitutions and it satisfies  Assumption \ref{as2.1} \textit{$(H1)$-$(H3)$} and Assumption (\textit{$H4$}). Moreover, it remains to verify that the associated linear system of the equation \eqref{1.1} is approximately controllable. In order to prove this, we consider
	$$(T-t)^{\alpha-1}\mathrm{B}^*\widehat{\mathcal{T}}_{\alpha,p}(T-t)^*x^*=0,\ \mbox{ for any}\ x^*\in\mathbb{X}^*,\ 0\le t<T.$$ Therefore, we have $$\mathrm{B}^*\widehat{\mathcal{T}}_{\alpha,p}(T-t)^*x^*=0,\ \mbox{ for any}\ x^*\in\mathbb{X}^*,\ 0\le t<T,$$ and since the operator $\mathrm{B}^*=\mathrm{B}$ is one-one, then we obtain
	$$\widehat{\mathcal{T}}_{\alpha,p}(T-t)^*x^*=0,\ \mbox{ for all } \ t\in [0,T).$$
	Further, we have
	\begin{align*}
	\widehat{\mathcal{T}}_{\alpha,p}(T-t)^*x^*&=\int_{0}^{\infty}\xi\varphi_{\alpha}(\xi)\sum_{n=1}^{\infty}\exp\left(-n^2(T-t)^\alpha\xi\right)\langle x^*, w_{n} \rangle  w_{n}\mathrm{d}\xi=0\nonumber\\&\implies\langle x^*, w_{n} \rangle =0, \ \mbox{ for all }\  n\in\mathbb{N},
	\end{align*}
which implies $x^*=0$.	Hence, the approximately controllability of the linear system follows by using Lemma \ref{lem4.2} and Remark \ref{rem3.4}.  Finally, by invoking Theorem \ref{thm4.4}, we deduce that the semilinear system \eqref{1.1} (equivalent to the system \eqref{ex}) is approximately controllable.
	
	%			\begin{lem}[Lemma 2.3, \cite{SR}]\label{lema2.1}
	%				Let $x:(-\infty,T]\rightarrow \mathbb{X}$, be a function such that $x_0=\psi$ and $x|_J\in \mathcal{PC}$. Then
	%				$$\norm{x_s}_\mathfrak{B}\le H_1\norm{\psi}_\mathcal{P}+H_2\sup\{\norm{x(\theta)}_\mathbb{X}:\theta \in [0,\max\{0,s\}]\},\ s\in\mathcal{Q}(\rho^-)\cup J,$$
	%				where $$H_1=\sup_{t \in \mathcal{Q}(\beta-)}\varTheta^\psi(t)+\sup_{t\in J} \Upsilon(t),\ H_2=\sup_{t\in J}\Lambda(t). $$
	%			\end{lem}

	\medskip\noindent
	{\bf Acknowledgments:} The first author would like to thank Council of Scientific and Industrial Research, New Delhi, Government of India (File No. 09/143(0931)/2013 EMR-I), for financial support to carry out his research work and Department of Mathematics, Indian Institute of Technology Roorkee (IIT Roorkee), for providing stimulating scientific environment and resources. M. T. Mohan would  like to thank the Department of Science and Technology (DST), Govt of India for Innovation in Science Pursuit for Inspired Research (INSPIRE) Faculty Award (IFA17-MA110). 
	
\end{document}